%% file: Kuijlaars-Molag-arXiv-version.tex
\documentclass[a4paper, fleqn]{article}
\usepackage[english]{babel}
\usepackage{graphicx}
\usepackage{amsmath, amssymb, amsthm, latexsym} 
\usepackage{float,fullpage}
\usepackage{amsthm}
\usepackage{mathtools}
\usepackage{xparse}
\usepackage{makeidx}
\usepackage{lipsum}

\usepackage{tikz}
\usepackage{comment}

\newtheorem{theorem}{Theorem}[section]
\newtheorem{lemma}[theorem]{Lemma}
\newtheorem{proposition}[theorem]{Proposition}
\newtheorem{corollary}[theorem]{Corollary}

\theoremstyle{definition}
\newtheorem{definition}[theorem]{Definition}
\newtheorem{remark}[theorem]{Remark}

\newtheorem{rhproblem}[theorem]{Riemann-Hilbert Problem}

\providecommand\phantomsection{}

\numberwithin{equation}{section}

\DeclareMathOperator{\supp}{supp}

\DeclareMathOperator{\IM}{Im}
\DeclareMathOperator{\RE}{Re}
\renewcommand{\Im}{\IM}
\renewcommand{\Re}{\RE}
\newcommand{\ds}{\displaystyle}

\DeclarePairedDelimiterX\MeijerM[3]{\lparen}{\rparen}{\begin{matrix}#1 \\ #2\end{matrix}\delimsize\vert\,#3}

\newcommand\MeijerG[7]{ G^{\,#1,#2}_{#3,#4}\MeijerM*{#5}{#6}{#7}}

\makeindex

\begin{document}

\clearpage
\thispagestyle{empty}

\title{The local universality
of Muttalib-Borodin biorthogonal ensembles with 
parameter $\theta = \frac{1}{2}$}
\author{A. B. J. Kuijlaars and L. D. Molag}
\maketitle
\begin{center}
Katholieke Universiteit Leuven, Department of Mathematics,\\
Celestijnenlaan 200B box 2400, 3001 Leuven, Belgium.\\
E-mail: arno.kuijlaars@kuleuven.be, leslie.molag@kuleuven.be
\end{center}

\begin{abstract} 
The Muttalib-Borodin biorthogonal ensemble is a probability density function for $n$ particles on the positive real line that depends on a parameter $\theta$ and an external field $V$. For $\theta=\frac{1}{2}$ we find the large $n$ behavior of the associated correlation kernel with only few restrictions on $V$. The idea is to relate the ensemble to a type II multiple orthogonal polynomial ensemble that can in turn be related to a $3\times 3$ Riemann-Hilbert problem which we then solve with the Deift-Zhou steepest descent method. The main ingredient is the construction of the local parametrix at the origin, with the help of Meijer G-functions, and its matching condition with a global parametrix. We will present a new iterative technique to obtain the matching condition, which we expect to be applicable in more general situations as well. 
\end{abstract}

\tableofcontents

\section{Introduction and statement of results}
\label{sec:intro}

\subsection{The Muttalib-Borodin ensemble}

The Muttalib-Borodin biorthogonal ensemble with parameter $\theta>0$ and weight
function $w$ is the following probability density function for 
the position of $n$ particles 
on the positive half line $[0,\infty)$
\begin{align}
\label{MBe}
\frac{1}{Z_{n}} \prod_{j<k} (x_{k}-x_{j})(x_{k}^{\theta}-x_{j}^{\theta}) \prod_{j=1}^{n} w(x_{j}).
\end{align}
We will consider an $n$-dependent weight function 
\begin{align}
\label{wx}
w(x)= x^{\alpha} e^{-n V(x)} 
\end{align}
with $\alpha > -1$ and an external field $V$ that has enough increase 
at infinity.

The model is named after Muttalib \cite{Mu}, who introduced it as a 
simplified model for disordered conductors in the metallic regime and
Borodin \cite{Bo}, who obtained profound mathematical results for specific
weights, in particular for Laguerre and Jacobi weights. 
Due to its relation to eigenvalue distributions of random matrix models
\cite{Ch, FoWa, KuSt} the model has attracted considerable attention in recent years. 

In the large $n$ limit the particles have an almost sure limiting 
empirical measure $\mu^*$ which minimizes a corresponding equilibrium problem
\cite{BlLeToWi,Bu,EiSoSt}. The equilibrium problem was studied in detail in \cite{ClRo}.
In \cite{Ku} it was pointed out that there is an equivalent vector
equilibrium problem for vectors of $q+r-1$ measures, when 
$\theta=\frac{q}{r}$ with $q,r\in\mathbb{N}$ and $q<r$. 

The Muttalib-Borodin ensemble is a determinantal point process (in fact a biorthogonal ensemble \cite{Bo}) and as such it has a corresponding correlation kernel $K^{\alpha,\theta}_{V, n}$. The kernel is given by
\begin{align}
\label{KVn}
	K^{\alpha,\theta}_{V,n}(x,y)
	= w(y) \sum_{j=0}^{n-1} p_j(x) q_j(y^{\theta}) 
\end{align}
where  $p_j$, $q_j$ are polynomials of degree $j$  satisfying
for $j,k=0,1, \ldots, n-1$,
\begin{align} 
\label{biortho} 
	\int_0^{\infty} p_j(x) q_k(x^{\theta}) w(x) dx = \delta_{j,k} 
	\end{align}
see \cite{Mu}. Such polynomials exist uniquely (up to multiplicative constants) 
under certain conditions, as was shown in \cite{ClRo}.

Borodin \cite{Bo} computed the hard edge scaling limit for the Laguerre
case, namely if $V(x) = x$, then
\begin{align}
\label{KVnlimit}
	\lim_{n \to \infty} \frac{1}{n^{1+ 1/\theta}}	
	K^{\alpha, \theta}_{V,n} \left(\frac{x}{n^{1+ 1/\theta}}, 
	\frac{y}{n^{1+1/\theta}}\right)
	= \mathbb K^{(\alpha, \theta)}(x,y) 
\end{align}
with limiting kernel
\begin{align} 
\label{KVnlimit2}
\mathbb K^{(\alpha,\theta)}(x,y) = 
	 \theta y^{\alpha} \int_0^1 J_{\frac{\alpha+1}{\theta},\frac{1}{\theta}}(ux)
	 J_{\alpha+1,\theta}\left((uy)^{\theta}\right) u^{\alpha} du  
\end{align}
where 
\[ J_{a,b}(x) = \sum_{j=0}^{\infty} \frac{(-x)^j}{j! \Gamma(a+bj)}. \]
is Wright's generalization of the Bessel function.
For $\theta = 1$, the limit \eqref{KVnlimit} reduces to the well-known
Bessel kernel in the theory of random matrices. 

Several other expressions are known for the kernel \eqref{KVnlimit2}.
For example, Zhang \cite[Theorem 1.2]{Zh} gives a double contour integral,
and if $\theta$ or $1/\theta$ is an integer, then  \eqref{KVnlimit}
can be expressed in terms of Meijer G-functions \cite{KuSt}.
These so-called Meijer G-kernels appear in singular value
distributions for products of random matrices
\cite{AkIpKi,AkKiWe} and in their hard edge scaling limit
 \cite{AkSt1, AkSt2,FoLi, FoWa, KiKuSt, KuZh}. 
 See  \cite{AkIp} for a survey paper.
The bulk and soft edge scaling limits for singular values
of Ginibre random matrices are the  
usual sine and Airy kernels \cite{LiWaZh}, 
and these classical limits were also established 
for the Muttalib-Borodin model in the Laguerre case \cite{Zh}.

It is natural to expect that the limit \eqref{KVnlimit} is not restricted
to the case $V(x) = x$, but holds for much more general external fields. 
In this paper we consider $\theta = \frac{1}{2}$ and we
show that the hard edge scaling limit \eqref{KVnlimit}
indeed holds for a large class of external fields $V$.

\subsection{Statement of main result}

The main assumption on the external field $V$ concerns the behavior
of an equilibrium measure, that we describe first. 
We use
\[ I(\mu) = \iint \log \frac{1}{|x-y|} d\mu(x) d\mu(y),
	\qquad I(\mu,\nu) = \iint \log \frac{1}{|x-y|} d\mu(x) d\nu(y), \]
to denote the logarithmic energy of the measure $\mu$,
and the mutual energy of two measures $\mu$ and $\nu$, respectively.
We also write 
\[ I_{\theta}(\mu) = \iint \log \frac{1}{|x^{\theta}-y^{\theta}|} d\mu(x) d\mu(y). \]
The equilibrium problem that is relevant for \eqref{MBe}
with weight \eqref{wx} is the following. Minimize
\begin{equation} \label{IVtheta}
	I_{V,\theta}(\mu) = 
	\frac{1}{2} I(\mu) + \frac{1}{2} I_{\theta}(\mu) + \int V(x) d\mu(x) 
	\end{equation}
among all probability measures $\mu$ on $[0,\infty)$. We assume that
$V$ is continuous and 
\[ \lim_{x \to +\infty} \frac{V(x)}{\log x} = +\infty. \]
Then existence and
uniqueness of the minimizer can be shown with usual methods from logarithmic
potential theory \cite{De,SaTo} as the functional $I_{V,\theta}$ is
lower semi-continuous and strictly convex on the set of all 
probability measures with finite logarithmic energy. The minimizer $\mu^*_{V,\theta}$ 
is called the $\theta$-equilibrium measure in the presence of the 
external field $V$. If $V$ and $\theta$ are clear from the context
then we simply say  equilibrium measure.

It is known that the equilibrium measure has a compact support
and it is characterized by the Euler-Lagrange variational conditions, see e.g.\ \cite{ClRo},
\begin{align} \label{ELcon1}
	\int \log |x-s| d\mu_{V,\theta}^*(s) + \int \log|x^{\theta} - s^{\theta}| d\mu_{V,\theta}^*(s)
		\begin{cases} = V(x) + \ell,  & x \in \supp(\mu), \\
	  \leq V(x) +  \ell, &  x \in [0,\infty),
	  \end{cases}
\end{align}
for some $\ell \in \mathbb R$. The equilibrium problem can be analyzed
exactly in case $V(x) = x$, see \cite[section 4.5.1]{ClRo} where this was
done for $\theta > 1$ and in particular $\theta = 2$, 
and see Proposition \ref{density} below for $\theta = \frac{1}{2}$. 
The equilibrium measure has
a support $[0,q]$ for some $q > 0$ with a density that has
a square root behavior at $q$ and behaves
like $s^{-1/(1+\theta)}$ as $s \to 0+$.

This is a generic behavior in a wider class of external fields, 
as was shown by Claeys and Romano \cite[Theorem 1.8]{ClRo},
see also Proposition \ref{density} below,
and our results apply to external fields that are one-cut $\theta$-regular in
the following sense.

\begin{definition} \label{ocr}
We call the external field $V$ one-cut $\theta$-regular if
the $\theta$ equilibrium measure $\mu^*_{V,\theta}$ in external
field $V$ is supported on one interval $[0,q]$ for some $q > 0$
with a density that is positive on $(0,q)$ and satisfies
\begin{align} \label{EPreg} 
	\frac{d\mu^*_{V,\theta}(s)}{ds} = 
	\begin{cases} c_{0,V} s^{-\frac{1}{\theta+1}}
	\left(1 + o(1) \right), 
	& \text{ as } s \to 0+, \\[5pt]
	c_{1,V} (q-s)^{1/2} \left(1+ o(1)\right),
	& \text{ as } s \to q-,
	\end{cases} \end{align}
	with positive constants $c_{0,V}, c_{1,V} >0$,
and, in addition, the inequality in \eqref{ELcon1}  is strict for  $x > q$.
\end{definition}
The conditions of Definition \ref{ocr} hold generically for
external fields $V$ that do not push the equilibrium measure
away from $0$. See Proposition \ref{prop:ocr} below for a sufficient
condition on $V$ in case $\theta = \frac{1}{2}$.
For $V(x) = x$ we are able to determine the density of
$\mu^*_{V,\frac{1}{2}}$ explicitly in Proposition \ref{density}. From this
formula the conditions \eqref{EPreg} can be verified directly.

The main result of our paper is the following theorem.

\begin{theorem}\label{mainThm} 
Let $\alpha>-1$, $\theta = \frac{1}{2}$,  and let $V: [0,\infty) \to \mathbb R$ 
be a one-cut $\theta$-regular external field in the sense of Definition \ref{ocr},
which is real analytic on $[0,\infty)$. 
Then for $x,y\in (0,\infty)$ we have
\begin{align} \label{Kalth}
\lim_{n\to\infty} \frac{1}{(c_V n)^3} K^{\alpha,\frac{1}{2}}_{V,n}
\left(\frac{x}{(c_V n)^3},\frac{y}{(c_V n)^3}\right) =
 \mathbb K^{(\alpha, \frac{1}{2})}(x,y)
\end{align}
uniformly on compact subsets of $(0,\infty)$, where $c_{V} = \frac{2\pi}{\sqrt{3}} c_{0,V}$,
and $c_{0,V} > 0$ is the constant from \eqref{EPreg}.
\end{theorem}
In the bulk and at the soft edge we find the usual sine and Airy kernels as the
scaling limits of $K_{V,n}^{\alpha,\frac{1}{2}}$ as $n \to \infty$,
but we do not give the details for that.

\begin{remark}
The one-cut assumption is not essential, but is
included for convenience only. Theorem \ref{mainThm} remains valid if the equilbrium measure 
is supported on several
intervals, provided that one of the intervals is of the form $[0,q]$ and the density
of the equilibrium measure has the behavior \eqref{EPreg} as $s \to 0+$. Similarly, we may
relax the regularity assumption. Theorem \ref{mainThm} also holds if the density of the
equilibrium measure vanishes somewhere in the interior of its support, or if there is
a higher order vanishing at a non-zero endpoint of the support, or if the inequality \eqref{ELcon1} is
an equality at a finite number of points outside of the support. These singular cases are also present
in the classical case $\theta=1$ that is connected with orthogonal polynomials
and $2 \times 2$ matrix valued RH problems \cite{DeKrML1}.
Each of these singular cases can be treated by means of special local parametrices, and this
would work in a similar way in our situation.  
\end{remark}

In section \ref{sec:prelim} we discuss the main ingredients
for the proof of Theorem \ref{mainThm}. These include
the equilibrium measure and the interpretation
of the biorthogonality property \eqref{biortho} as multiple orthogonality in
case $\theta =1/r$ with an integer $r \in \mathbb N$.
The multiple orthogonality leads to a Riemann-Hilbert (RH) problem
\cite{VAGeKu} of size $(r+1) \times (r+1)$. 
To analyze the large $n$ limit we are going
to apply the Deift-Zhou steepest descent method \cite{DeZh} that was
first applied to orthogonal polynomials in the influential papers \cite{DeKrML1,DeKrML2}. The Meijer G-functions 
appear in the construction of a local parametrix at $0$.
The model RH problem for Meijer G-functions is discussed
in section \ref{sec:RHproblem}.

The steepest descent analysis is in sections \ref{sec:first}--\ref{sec:main}.
The main technical  difficulty is the matching of the local
parametrix with the global parametrix. We deal with this issue
in section \ref{sec:local}. The matching is a  serious problem
for larger size RH problems as has been observed in other
situations \cite{BeBo, DeKuZh, DeKu, KuMFWi}.  As in \cite{BeBo}
we develop an iterative technique to improve the matching step by step. 
The proof of Theorem \ref{mainThm} is in the final
section \ref{sec:main}.

\begin{remark}
	We conjecture that  Theorem \ref{mainThm} holds for any $\theta > 0$,
	with the exponent $3$ replaced by $1 +1/\theta$ as in \eqref{KVnlimit}.
	Namely, for real analytic one-cut $\theta$-regular
	external fields $V$ we expect that	there is a constant $c_V > 0$ such that
	\begin{align} \label{Kalthgeneral}
	\lim_{n\to\infty} \frac{1}{(c_V n)^{1+1/\theta}} K^{\alpha,\theta}_{V,n}
	\left(\frac{x}{(c_V n)^{1+1/\theta}},\frac{y}{(c_V n)^{1+1/\theta}}\right) =
	\mathbb K^{(\alpha, \theta)}(x,y).
	\end{align}
	It seems likely that the techniques of this paper extend to the case $\theta = 1/r$
	with $r \in \mathbb N$, $r \geq 3$. As already mentioned, this case leads to
	a RH problem of size $(r+1) \times (r+1)$. The main ingredients that we are
	going to use (see next section) are available for any integer $r$. The challenge
	will be to establish the matching condition which becomes technically more involved.
	
	The case $\theta = 2$ is related to $\theta = 1/2$ by a change of variables. 
	Putting $x_j = y_j^{1/2}$ for $j=1, \ldots, n$ in \eqref{MBe} with $\theta = 2$
	gives a Muttalib-Borodin ensemble with $\theta = 1/2$ and external field
	$V(\sqrt{x})$. Theorem~\ref{mainThm} applies provided that $x \mapsto V(\sqrt{x})$
	is real analytic on $[0,\infty)$. 
	This in turn means that \eqref{Kalthgeneral} with $\theta = 2$  holds provided that 
	$V$ is real analytic and even (and one-cut $\theta$-regular).  
	The extension to arbitrary real analytic external fields remains open.
	
	We use the assumption that $V$ is real analytic at $0$ in the proof of Lemma \ref{Pansatzlemma},
	and in particular for the identity \eqref{phi1jumpDelta2}.
\end{remark}

Throughout the paper we use the principal branches of fractional
powers, i.e., with a branch cut on the negative real axis and positive on the positive real axis.

\section{Preliminaries}
\label{sec:prelim}

We discuss the main ingredients for the proof of Theorem \ref{mainThm}.

\subsection{Multiple orthogonal polynomials}
The first step in the proof is the observation that the polynomials
$p_j$ that appear in the kernel \eqref{KVn} and that satisfy the
biorthogonality \eqref{biortho} can be viewed as multiple orthogonal
polynomials in case $\theta = 1/r$ and $r$ is an integer.
Without loss of generality we can take these polynomials to be monic.
Thus $p_n$ is a polynomial of degree $n$ that is characterized by the 
property that 
\begin{align} 
	\label{biortho2}
	\int_0^{\infty} p_n(x) x^{k \theta} w(x) dx = 0, \qquad k = 0,1, \ldots, n-1. 
	\end{align}
	
\begin{lemma}
Suppose $\theta = 1/r$ with $r$ an integer.
Then $p_n$ is the unique monic polynomial of degree $n$ that satisfies 
\begin{align} 
\label{biortho3} 
	\int_0^{\infty} p_n(x) x^k  w_i(x) dx = 0,
	\qquad i = 1, \ldots, r, \quad k = 0, \ldots, \lfloor \tfrac{n-i}{r} \rfloor, 
\end{align}
where $w_i(x) = x^{(i-1) \theta} w(x)$ for $i =1, \ldots, r$.
\end{lemma}
\begin{proof}
If $i, k$  are as in \eqref{biortho3} then $rk + i-1 \leq n-1$, and
thus by \eqref{biortho2} we have
\[ 	\int_0^{\infty} p_n(x) x^{(rk + i-1) \theta} w(x) dx = 0 \]
which reduces to \eqref{biortho3} in view of the definition of $w_i(x)$.

Thus \eqref{biortho2} implies \eqref{biortho3} and it is easy to
see that the converse holds as well.
\end{proof}

The relations \eqref{biortho3} are multiple orthogonality conditions
with respect to  weight functions $w_1$, $w_2$, \ldots, $w_{r}$, and
with the multi-index $(n_1, \ldots, n_r)$ where 
$n_i = \lfloor \tfrac{n-i}{r} \rfloor + 1$
for $i=1, \ldots, r$. While \eqref{biortho3} is a simple reformulation
of \eqref{biortho} it has the advantage of leading to a RH
problem.

\subsection{Riemann-Hilbert problem}
We state the RH problem  for the case $r=2$. In this
case the RH problem has size $3 \times 3$. For general $r$ the size is
 $(r+1) \times (r+1)$, see \cite{VAGeKu}.

\begin{rhproblem} \label{RHPforY} \
\begin{description}
\item[RH-Y1] $Y : \mathbb{C}\setminus [0,\infty)\to \mathbb{C}^{3\times 3}$ is analytic.
\item[RH-Y2] $Y$ has boundary values for $x\in (0,\infty)$, denoted by $Y_{+}(x)$ (from the upper half plane) and $Y_{-}(x)$ (from the lower half plane), and
\begin{equation} \label{RHY2}
Y_{+}(x) = Y_{-}(x) \begin{pmatrix} 1 & w(x) & x^{\frac{1}{2}} w(x) \\ 
0 & 1 & 0\\ 0 & 0 & 1 \end{pmatrix}, \qquad  x > 0.
\end{equation}
\item[RH-Y3] As $|z|\to\infty$
\begin{align} \label{RHY3}
Y(z) = \left(\mathbb{I}+\mathcal{O}\left(\frac{1}{z}\right)\right) 
\begin{pmatrix} z^{n} & 0 & 0\\ 0 & z^{-\lceil \frac{n}{2} \rceil} & 0\\ 0 & 0 & z^{-\lfloor\frac{n}{2}\rfloor}\end{pmatrix}.
\end{align}
\item[RH-Y4] As $z\to 0$
\begin{align}
\label{RHY4}
Y(z) = \mathcal{O}\begin{pmatrix} 1 & h_{\alpha}(z) & h_{\alpha+\frac{1}{2}}(z)\\ 
1 & h_{\alpha}(z) & h_{\alpha+\frac{1}{2}}(z)\\ 
1 & h_{\alpha}(z) & h_{\alpha+\frac{1}{2}}(z) \end{pmatrix}
\text{ with } \, h_{\alpha}(z) = 
\begin{cases} |z|^{\alpha}, & \text{if } \alpha < 0, \\ 
	\log{|z|}, & \text{if } \alpha = 0,\\ 
	1, & \text{if } \alpha > 0.
	\end{cases}
\end{align}
\end{description}
The $\mathcal{O}$ condition in \eqref{RHY3} and \eqref{RHY4} is to be taken entry-wise.
\end{rhproblem}

The RH problem is the analogue of the Fokas-Its-Kitaev RH problem \cite{FoItKi}  for  orthogonal polynomials. 
The condition \textbf{RH-Y4} is an endpoint condition that is
analogous to the  endpoint  condition in \cite{KuMLVAVa}.

There is a unique solution of the RH problem (since the polynomial
$p_n$ uniquely exists \cite{ClRo})  and 
\[ Y_{11}(z) = p_n(z),
	\quad 
	Y_{12}(z) = \frac{1}{2\pi i} \int_0^{\infty}	
	\frac{p_n(x) w(x)}{x-z} dx,
	\quad
	Y_{13}(z) = \frac{1}{2\pi i} \int_0^{\infty}	
	\frac{p_n(x) x^{\frac{1}{2}} w(x)}{x-z} dx,
	\]
are the entries in the first row of the solution $Y$. The other two rows
are built in the same way out of suitable polynomials of lower degree.
Indeed, $Y_{21}$ and $Y_{31}$ are two polynomials of degree $\leq n-1$,
that are expressible in terms of the biorthogonal polynomials $p_{n-1}$
and $p_{n-2}$ of degrees $n-1$ and $n-2$, and for $j=2,3$,
\[ Y_{j2}(z) = \frac{1}{2\pi i} \int_0^{\infty}	
	\frac{Y_{j1}(x) w(x)}{x-z} dx,
	\quad
	Y_{j3}(z) = \frac{1}{2\pi i} \int_0^{\infty}	
	\frac{Y_{j1}(x) x^{\frac{1}{2}} w(x)}{x-z} dx. \]

\begin{lemma} \label{lemCorKerY}
The correlation kernel \eqref{KVn} with $\theta = \frac{1}{2}$ 
is expressed in terms of the solution of the RH problem
for $Y$ as follows
\begin{align} \nonumber
K^{\alpha,\frac{1}{2}}_{V,n}(x,y) & = \frac{1}{2\pi i(x-y)} 
\begin{pmatrix} 0 & w(y) & y^{\frac{1}{2}} w(y) \end{pmatrix} Y_{+}^{-1}(y) Y_{+}(x) \begin{pmatrix} 1 \\ 0 \\ 0\end{pmatrix} \\
& = \frac{y^{\beta} e^{-nV(y)}}{2\pi i(x-y)} 
\begin{pmatrix} 0 & y^{-\frac{1}{4}} & y^{\frac{1}{4}} \end{pmatrix} Y_{+}^{-1}(y) Y_{+}(x) \begin{pmatrix} 1 \\ 0 \\ 0\end{pmatrix},
	\qquad \beta = \alpha + \frac{1}{4},
\label{CorKerY}
\end{align}
for $x,y > 0$.
\end{lemma}
\begin{proof}
The first identity is due to  \cite{BlKu1}. See also \cite{DaKu} 
for the extension  to $r \geq 3$. The second identity in \eqref{CorKerY}
follows from the form \eqref{wx} of the weight function $w(y)$.
\end{proof}

The first identity in \eqref{CorKerY} is a Christoffel-Darboux formula that expresses
the kernel in terms of a finite number of biorthogonal polynomials. 
See also \cite{ClRo} for a Christoffel-Darboux formula for general rational
values of $\theta$. For asymptotic analysis the expression \eqref{CorKerY}
in terms of the RH problem is very convenient. 

\subsection{Vector equilibrium problem}

 The Deift-Zhou method of steepest descent consists of a 
 number of explicit transformations of the RH problem.
One of the transformations depends on the $\theta$-equilibrium measure
$\mu_{V,\theta}^*$ that minimizes the energy functional \eqref{IVtheta}. 

When $\theta = 1/r$ with $r$ an integer, the $\theta$-equilibrium measure
can also be characterized as the first component of the minimizer of a vector equilibrium problem \cite{Ku} and this will be important for us. 
We state it here for the case $r=2$.
In case $r = 2$, the vector equilibrium problem asks to minimize the
energy functional
\begin{equation} \label{Energy2} 
	I(\mu) - I(\mu,\nu) + I(\nu) + \int V(x) d \mu(x) 
	\end{equation}
among all pairs of measures $(\mu,\nu)$ such that
\begin{itemize}
\item $\mu$ is a probability measure on $[0,\infty)$,
\item $\nu$ is a measure on $(-\infty,0]$ of total mass $1/2$.
\end{itemize}
For general $\theta = 1/r$ the vector equilibrium problem
has $r$ measures.

It is shown in \cite{Ku} that there is a unique minimizer $(\mu^*,\nu^*)$
and the first component $\mu^*$ coincides with the minimizer $\mu^*_{V,
\frac{1}{2}}$
of the $\theta$-energy functional \eqref{IVtheta} with 
$\theta = \frac{1}{2}$.
In addition, the measure $\nu^*$ has full support
$\supp(\nu^*) = (-\infty,0]$, and the variational conditions
\begin{align} \label{ELvec1}
	2 \int \log |x-s| d\mu^*(s) - \int \log|x-s| d\nu^*(s)
		& 
		\begin{cases} = V(x) + \ell,  & x \in \supp(\mu^*), \\
		\leq V(x)  +  \ell, &  x \in [0,\infty), \end{cases} \\
		\label{ELvec2}
	2 \int \log |x-s| d\nu^*(s) - \int \log|x-s| d\mu^*(s)
		& = 0, \qquad \qquad x \in (-\infty,0],
\end{align}
are satisfied, with a constant $\ell$ that is possibly different 
from the one appearing in \eqref{ELcon1}. The identity \eqref{ELvec2} expresses
that $2\nu^*$ is the balayage of $\mu^*$ onto the negative real axis.

We use the two measures from the vector equilibrium problem 
in the steepest descent analysis that will follow. We also
use the vector equilibrium problem in the proof of the
following sufficient condition for an external field
to be one-cut regular.

\begin{proposition} \label{prop:ocr}
Suppose $V$ is twice differentiable on $[0,\infty)$ such
that $x V'(x)$ is increasing for $x > 0$.
Then $V$ is one-cut $\frac{1}{2}$-regular in the sense
of Definition \ref{ocr}.
\end{proposition}
\begin{proof}
Because of \eqref{Energy2} we have that $\mu^* = \mu^*_{V, \frac{1}{2}}$
is the minimizer of the energy functional
$ I(\mu) + \int \tilde{V} d\mu$ with
the modified potential
\[ \tilde{V}(x) = V(x) + \int \log(x-s) d\nu^*(s). \]
For every $s < 0$, we have that $x \mapsto \frac{x}{x-s}$
is increasing on $[0,\infty)$. Since $\nu^*$ is supported
on $(-\infty,0]$, it then easily follows from this and
the assumption of the proposition that 
$ x \mapsto x \tilde{V}'(x) = xV'(x) + \int \frac{x}{x-s} d\nu^*(s)$
is strictly increasing for $x \in [0,\infty)$. 
It then follows by well-known results on equilibrium
measures in external fields that $\supp(\mu^*) = [0,q]$
is an interval containing $0$, 
see e.g.\ \cite[Theorem IV.1.10 (c)]{SaTo}.

Knowing that the support is an interval we can use
\cite[Theorem 1.11 and formula (1.28)]{ClRo} to conclude that the density of
$\mu^*$ is positive on $(0,q)$ and has the required endpoint
behavior \eqref{EPreg}. Actually, the result in \cite{ClRo}
is stated for $\theta > 1$, but it also applies to $0 < \theta < 1$.

It remains to show that the inequality \eqref{ELcon1} is strict 
for $x > q$. From \eqref{ELvec1} we have
\[ \tilde{V}(x) + \ell - 2 \int \log(x-s) d\mu^*(s) \geq 0,
	\qquad x \geq q \]
with equality for $x = q$. Then the derivative is non-negative
at $x=q$, which implies
\begin{equation} \label{propocr} 
	x \tilde{V}'(x) - 2 \int \frac{x}{x-s} d\mu^*(s) \geq 0 
	\end{equation}
for $x =q$. The function in the left-hand side of \eqref{propocr} is strictly increasing 
for $x \geq q$.  To see this we use that $x \tilde{V}'(x)$ is strictly
increasing and the fact that $x \mapsto \frac{x}{x-s}$, $x > q$,
is decreasing for every $s \in [0,q]$.
It follows that strict inequality holds in \eqref{propocr}
for $x > q$, which in turn leads to the property 
that \eqref{ELcon1} is strict for $x > q$. 
This concludes the proof of the proposition. 
\end{proof}

For $V(x) = x$ we have the following explicit result.
\begin{proposition} \label{density}
When $V(x)=x$ and $\theta = \frac{1}{2}$, 
the equilibrium measure of the Muttalib-Borodin ensemble is
supported on $[0, \frac{27}{8}]$ with  density
\begin{multline} \label{Visxdensity}
	\frac{d\mu_{V,\frac{1}{2}}^*(s)}{ds}
	= \frac{\sqrt{3}}{4 \pi s^{2/3}}
	\left[ \left( 1 - \frac{4s}{3} + \frac{8s^2}{27}
		+ \sqrt{1 - \frac{8s}{27}} \right)^{1/3}
		+
	\left( -1 + \frac{4s}{3} - \frac{8s^2}{27}
		+ \sqrt{1 - \frac{8s}{27}} \right)^{1/3}
		 \right] \\
	\qquad 0 < s < \frac{27}{8}.
\end{multline}
In particular, \eqref{EPreg} holds with constants
$q = \frac{27}{8}$, $c_{0,V} = \frac{\sqrt{3}}{2^{5/3}\pi}$
and $c_{1,V} = \frac{16\sqrt{2}}{81 \pi}$.
\end{proposition} 
\begin{proof} 
By Proposition \ref{prop:ocr} we have $\supp(\mu^*) = [0,q]$
for some $q > 0$. By differentiating the identities in  
\eqref{ELvec1}-\eqref{ELvec2} we see
that the functions $F_1(z) = \ds \int \frac{d\mu^*(s)}{z-s}$
and $F_2(z) = \ds \int \frac{d\nu^*(s)}{z-s}$ satisfy
\begin{equation} \label{muVx1}
\begin{aligned} 
	F_{1,+}(x) + F_{1,-}(x) - F_2(x) & = V'(x) = 1, &&  x \in (0,q), \\
	F_{2,+}(x) + F_{2,-}(x) - F_1(x) & = 0, && x \in (-\infty,0).
	\end{aligned}
	\end{equation}
Let $\mathfrak R$ be the compact three sheeted Riemann surface
with sheets $\mathfrak R_0 = \overline{\mathbb C} \setminus [0,q]$,
$\mathfrak R_1 = \mathbb C \setminus (-\infty, q]$ and
$\mathfrak R_2 = \mathbb C \setminus (-\infty, 0]$ with sheet
structure as in Figure \ref{FigRS}. Then it follows from
\eqref{muVx1} that $\Psi$ defined by
\begin{equation} \label{muVx2} 
	\Psi(z) = \begin{cases}
	\Psi_0(z) = 1 - F_1(z), & \text{ if } z \in \mathfrak{R}_0, \\
	\Psi_1(z) = F_1(z) - F_2(z), & \text{ if } z \in \mathfrak{R}_1, \\
	\Psi_2(z) = F_2(z), & \text{ if } z \in \mathfrak{R}_2, \end{cases} 
	\end{equation}
is meromorphic on $\mathfrak{R}$. We note the asymptotic
behaviors
\begin{equation} \label{muVx3}
	\Psi_0(z)  = 1 - \frac{1}{z} + \mathcal{O}(z^{-2}), \qquad
	\Psi_1(z)  = \frac{1}{2z} + \mathcal{O}(z^{-3/2}), \qquad
	\Psi_2(z)  = \frac{1}{2z} + \mathcal{O}(z^{-3/2}), 
	\end{equation}
	as $z \to \infty$.
Thus $\Psi$ has a double zero at the point at infinity that 
is common to sheets $\mathfrak{R}_1$ and $\mathfrak{R}_2$. 
From Proposition \ref{prop:ocr} we know that $\mu^*$ has
a density that behaves like $c_{0,V} s^{-2/3}$ as $s \to 0+$,
see \eqref{EPreg}. Then it can be shown that  
$F_1(z) \sim -c z^{-2/3}$ as $z \to 0$ with $c= \frac{2\pi}{\sqrt{3}}c_{0,V}$, 
	which means that $\Psi$
has a double pole at $z=0$, as $z=0$ is a double branch point of
the Riemann surface. There are no other poles and zeros of $\Psi$.

\begin{figure}
	\centering
	\begin{picture}(200,120)(-20,160)
	\unitlength=2pt
	\put(0,100){\line(1,0){100}}
	\put(0,100){\line(2,1){20}}
	\put(20,110){\line(1,0){100}}
	\put(100,100){\line(2,1){20}}
	\put(0,80){\line(1,0){100}}
	\put(0,80){\line(2,1){20}}
	\put(20,90){\line(1,0){100}}
	\put(100,80){\line(2,1){20}}
	\put(0,120){\line(1,0){100}}
	\put(0,120){\line(2,1){20}}
	\put(20,130){\line(1,0){100}}
	\put(100,120){\line(2,1){20}}
	\multiput(60,105)(0,2){10}{\line(0,1){1}}
	\multiput(90,105)(0,2){10}{\line(0,1){1}}
	\multiput(60,105)(0,-2){10}{\line(0,-1){1}}
	\thicklines
	\put(60,105){\line(1,0){30}}
	\put(60,125){\line(1,0){30}}
	\put(10,105){\line(1,0){50}}
	\put(10,85){\line(1,0){50}}
	\put(60,105){\circle*{2}}
	\put(90,105){\circle*{2}}
	\put(60,125){\circle*{2}}
	\put(90,125){\circle*{2}}
	\put(60,85){\circle*{2}}
	\put(100,125){$\frak{R}_0$}
	\put(100,105){$\frak{R}_1$}
	\put(100,85){$\frak{R}_2$}
	\put(61,101){$0$}
	\put(90,101.5){$q$}
	\end{picture}
	\caption{The Riemann surface $\frak{R}$ used in the proof of Proposition \ref{density}. 
		The three sheets are defined by $\mathfrak{R}_{0} = \overline{\mathbb{C}}\setminus [0,q], \mathfrak{R}_{1} = \mathbb{C}\setminus (-\infty,q]$ and $\mathfrak{R}_{2} = \mathbb{C}\setminus (-\infty,0]$. The sheets are connected in
		the usual crosswise manner. \label{FigRS}}
\end{figure}
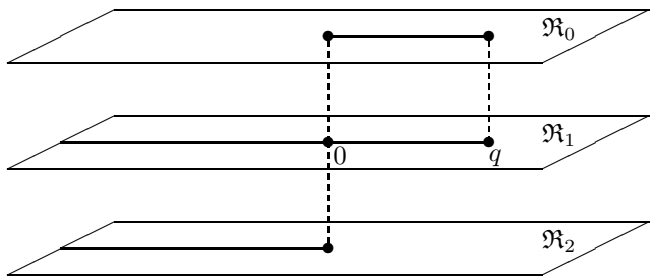

From \eqref{muVx3} we obtain
\begin{align*}
	\Psi_0(z) + \Psi_1(z) + \Psi_2(z) & = 1 + \mathcal{O}(z^{-3/2}), \\
	\Psi_0(z) \Psi_1(z) + \Psi_0(z) \Psi_2(z) + \Psi_1(z) \Psi_2(z) & = 
	\frac{1}{z} + \mathcal{O} (z^{-2}) \\
	\Psi_0(z) \Psi_1(z) \Psi_2(z) & = \frac{1}{4z^2} + \mathcal O(z^{-3})
\end{align*}
as $z \to \infty$.
The three above functions are meromorphic in the full complex plane
with a possible pole at $z=0$ only. Because $\Psi_j(z) = O(z^{-2/3})$
as $z \to 0$ for every $j=1,2,3$, it follows that all
$\mathcal O$ terms vanish identically. Thus
\[ 	\prod_{j=1}^3 (\zeta - \Psi_j(z))
	= \zeta^3 - \zeta^2 + \frac{1}{z} \zeta - \frac{1}{4z^2}
	\]
and
\begin{equation} \label{muVx4} 
	z^2 \zeta^3 - z^2 \zeta^2 + z \zeta - \frac{1}{4} = 0
\end{equation}
is the algebraic equation satisfied by $\Psi$.

The discriminant of \eqref{muVx4} is $z^4\frac{8z-27}{16}$
and it follows that there is a branch point at 
$z=\frac{27}{8}$. For real $z > \frac{27}{8}$ there are three 
real solutions of \eqref{muVx4}, and for real $z \in (0, \frac{27}{8})$ 
there is one real solution (namely $\zeta = \Psi_2(z)$) and 
two complex conjugate non-real solutions $\zeta = \Psi_j(z)$, $j=0,1$.

The density of $\mu^*$ satisfies
\begin{equation} \label{muVx5} 
	\frac{d\mu^*}{ds} = - \frac{1}{\pi} \Im F_{1,+}(s)
	= \frac{1}{\pi} \Im \Psi_{0,+}(s) \end{equation}
because of \eqref{muVx2} and the Stieltjes inversion formula.
Thus $\mu^*$ is supported on $[0, \frac{27}{8}]$.
With Cardano's method for solving a cubic equation
we find the solutions of \eqref{muVx4} for $z = s \in (0, \tfrac{27}{8})$.
$\Psi_0(s)$ is the unique solution with positive imaginary part
 and the calculations result in the formula \eqref{Visxdensity}.
\end{proof}

\subsection{Meijer G-functions}

The final main ingredient concerns the Meijer G-functions, that will
be used in the construction of a local parametrix at $0$.
In \cite{KuZh} the limiting kernel  
$\mathbb K^{(\alpha,\frac{1}{2})}(x,y)$ 
from our main result Theorem \ref{mainThm} was expressed
in terms of Meijer G-functions, see also 
Theorem \ref{thmKformula} below. 

The appearance of Meijer G-functions, or equivalently, generalized
hypergeometric functions, can also be expected
from the recent paper \cite{VA} on multiple Laguerre polynomials.
In case $w(x) = x^{\alpha} e^{-n x}$ the multiple orthogonal polynomials
from \eqref{biortho3} are such multiple Laguerre polynomials, 
and it was shown in \cite[Theorem 3]{VA} 
that in a suitable scaling regime these polynomials 
tend to a generalized hypergeometric function
\[ 
{{}_{0}F_{r}}\left( {- \atop 1+ \alpha_1,\dots, 1 +\alpha_r} ; -z\right)
= \sum_{k =0}^{\infty} \frac{1}{\left(1+\alpha_1\right)_{k}\cdots 
\left(1+\alpha_r\right)_{k}} \frac{(-z)^k}{k!}\]
where $\alpha_j = \alpha + \frac{j-1}{r}$ for $j=1,\ldots, r$.
This is the so-called Mehler-Heine asymptotics and it captures the
behavior of the polynomials near the origin. The case $r  = 2$ 
is relevant for the present paper and then we obtain
\begin{align} \label{phi0} 
	\phi_0(z) & := {{}_{0}F_{2}}\left( {- \atop 1+ \alpha, \frac{3}{2} +\alpha} ; -z\right)
	\end{align}
which is a solution of the third order differential equation
(we use $\vartheta = z \frac{d}{dz}$ as in \cite{OLBC}, for example),
\begin{align} \label{HGDE}
	\vartheta(\vartheta+ \alpha)(\vartheta+ \alpha +\tfrac{1}{2}) \phi
	+z\phi=0, \qquad \vartheta = z \frac{d}{dz}.
\end{align}

In case $2\alpha \not\in \mathbb Z$, other solutions of \eqref{HGDE} are $\ds z^{-\alpha} {{}_{0}F_{2}}\left( {- \atop 1 - \alpha, \frac{3}{2}} ; -z\right)$, and $\ds 
z^{-\alpha-\frac{1}{2}} {{}_{0}F_{2}}\left( {- \atop \frac{1}{2} - \alpha, \frac{1}{2}} ; -z\right) $.
Together with $\phi_0$ they are a basis of all solutions 
as they come from applying the Frobenius method
to \eqref{HGDE}. In case $2\alpha \in \mathbb Z$ (the resonant case), this method produces series solutions with additional log terms.
The Frobenius method provides a basis that is suitable for
the study of solutions near $z=0$.

We could have worked with the above ${{}_{0}F_{2}}$ 
hypergeometric functions. However, we prefer to use instead
a representation of solutions of \eqref{HGDE} in terms of 
Meijer G-functions, since  this representation is 
more convenient to describe the behavior  of solutions 
as $z \to \infty$.

The general Meijer G-function is defined by the following 
contour integral:
\begin{align}
\label{MeijerG}
\MeijerG{m}{n}{p}{q}{a_1, \ldots, a_p}{b_1, \ldots, b_q}{z} 
= \frac{1}{2\pi i} \int_{L} \frac{\prod_{j=1}^{m} \Gamma(b_j +s) \prod_{j=1}^{n} \Gamma(1-a_j -s)}{\prod_{j=m+1}^{q} \Gamma(1-b_j -s) \prod_{j=n+1}^{p} \Gamma(a_j +s)} z^{-s} ds,
\end{align}
where $\Gamma$ denotes the gamma function 
and empty products in \eqref{MeijerG} should be interpreted as $1$, as usual. 
The numbers $m, n, p$, and $q$ are integers with $0\leq m\leq q$ and 
$0\leq n \leq p$. See \cite[section 5.2]{Lu} for conditions on the parameters
$a_1, \ldots, a_p$ and $b_1, \ldots, b_q$ and the contour $L$.

By \cite[formula 16.18.1]{OLBC} and \eqref{phi0} we have
\begin{align} \nonumber
	\phi_0(z) & = \Gamma(1+\alpha) \Gamma(\tfrac{3}{2} + \alpha) 
	\MeijerG{1}{0}{0}{3}{ -}{0,-\alpha, - \alpha - \frac{1}{2}}{z} \\
	& =  \frac{\Gamma(1+\alpha) \Gamma(\frac{3}{2} + \alpha) }{2\pi i}
    	\int_L \frac{\Gamma(s)}{\Gamma(1+\alpha-s) \Gamma(\frac{3}{2} + \alpha-s)}
    	z^{-s} ds \label{MG1003} 
 \end{align}
with a Hankel contour $L$ that encircles the negative real axis.
Another solution of \eqref{HGDE}, see \cite[formula 16.21.1]{OLBC},
is
\begin{align}
\label{MG3003}
\MeijerG{3}{0}{0}{3}{ -}{0,-\alpha, - \alpha - \frac{1}{2}}{z} 
= \frac{1}{2\pi i} \int_{L} \Gamma(s) \Gamma(s-\alpha) \Gamma(s-\alpha - \tfrac{1}{2}) z^{-s} ds,
\end{align}
where now $L$ encircles the interval $(-\infty, \max(0,\alpha+\tfrac{1}{2})]$
in the complex $s$-plane (or alternatively it could be a vertical line  
$\Re s = c > \max(0,\alpha + \frac{1}{2})$ in the complex $s$-plane).

This solution behaves like $e^{-3 z^{\frac{1}{3}}}$ as $z \to \infty$
with $-4\pi < \arg z < 4 \pi$. In particular, it is
the recessive solution of \eqref{HGDE} in $-\pi < \arg z < \pi$.
More precise asymptotics are given in
\eqref{MeijerGexpansion}  below. Analytic continuations of 
\eqref{MG3003} to other Riemann sheets  provide other solutions.

Besides $\phi_0$ we will use  the following four solutions of \eqref{HGDE},
each defined in the sector $-\pi < \arg z < \pi$.

\begin{align} \label{phi1}
\phi_{1}(z) &= i e^{2\pi i\alpha} 
\MeijerG{3}{0}{0}{3}{-}{0,-\alpha, -\alpha-\frac{1}{2}}{z e^{2\pi i}}, \\
\label{phi2}
\phi_{2}(z) &= -i e^{-2\pi i\alpha} 
\MeijerG{3}{0}{0}{3}{-}{0,-\alpha, -\alpha-\frac{1}{2}}{z e^{-2\pi i}}, \\
\label{phi3}
\phi_{3}(z) &= 
\MeijerG{3}{0}{0}{3}{-}{0,-\alpha, -\alpha-\frac{1}{2}}{z}, \\
\label{phi4}
\phi_{4}(z) &= \phi_{1}(z)+\phi_{2}(z).
\end{align}

Here the notation $z e^{2\pi i}$ means that we have analytically continued $\phi_3$ along a counterclockwise loop around the origin. Likewise, the notation $z e^{-2\pi i}$ means that we have analytically continued $\phi_3$ along a clockwise loop around the origin. Indeed, these expressions are also given by \eqref{MG3003} by picking the argument of $z$ appropriately.

\begin{lemma} We have
\begin{equation} \label{phi4phi0}
 \phi_4(z) = - \frac{4 \pi^2}{\Gamma(1+\alpha)\Gamma(\frac{3}{2} + \alpha)}
	\phi_0(z) \end{equation}
\end{lemma}

\begin{proof}
We use the reflection formula 
$\Gamma(z) \Gamma(1-z) = \frac{\pi}{\sin \pi z}$ to 
write the combination of gamma functions appearing in \eqref{MG1003} as
\begin{align*}
	\frac{\Gamma(s)}{\Gamma(1+\alpha-s) \Gamma(\frac{3}{2}+\alpha -s)}
		& = 
	\Gamma(s) \Gamma(s-\alpha) \Gamma(s-\alpha-\tfrac{1}{2})
		\frac{\sin \pi (s-\alpha) \sin \pi(s-\alpha-\frac{1}{2})}{\pi^2} \\
	& =
	\Gamma(s) \Gamma(s-\alpha) \Gamma(s-\alpha-\tfrac{1}{2})
		\frac{1}{(2\pi i)^2} \left(e^{2\pi i(s-\beta)} + e^{-2\pi i(s-\beta)} \right).
	\end{align*}
Using this in \eqref{MG1003} and comparing with \eqref{MG3003},
	\eqref{phi1}, \eqref{phi2}, we see that \eqref{MG1003} is equal to 
$-\frac{1}{4\pi^2} (\phi_1(z) + \phi_2(z))$.
Then \eqref{phi4phi0} follows by \eqref{phi0}, \eqref{phi4}, and \eqref{MG1003}.	
\end{proof}

We collect the solutions into a matrix valued function
that is piecewise analytic in the complex plane.

\begin{definition}
We define (where $\vartheta = z \frac{d}{dz}$ as in \eqref{HGDE})
\begin{align} \label{defPhialpha}
\Phi_{\alpha}(z) = 
\begin{cases}
\begin{pmatrix} \phi_{1}(z) & \phi_{2}(z) & \phi_{3}(z)\\ 
	\vartheta \phi_{1}(z) & \vartheta \phi_{2}(z) & \vartheta \phi_{3}(z) \\
	\vartheta^2 \phi_{1}(z) & \vartheta^2 \phi_{2}(z) & \vartheta^2 \phi_{3}(z)\end{pmatrix}, & 0 < \arg (z)< \frac{\pi}{2}, \\
\begin{pmatrix} \phi_{4}(z) & \phi_{2}(z) & \phi_{3}(z)\\ 
	\vartheta \phi_{4}(z) & \vartheta \phi_{2}(z) & \vartheta \phi_{3}(z)\\
	\vartheta^2 \phi_{4}(z) & \vartheta^2 \phi_{2}(z) & \vartheta^2 \phi_{3}(z)\end{pmatrix}, & \frac{\pi}{2} < \arg (z)<\pi, \\
\begin{pmatrix}  \phi_{2}(z) & -\phi_{1}(z) & \phi_{3}(z)\\ 
\vartheta \phi_{2}(z) & -\vartheta \phi_{1}(z) & \vartheta \phi_{3}(z)\\  
\vartheta^2 \phi_{2}(z) & -\vartheta^{2}\phi_{1}(z) & \vartheta^2 \phi_{3}(z)\end{pmatrix}, & -\frac{\pi}{2}<\arg(z)<0,\\
\begin{pmatrix} \phi_{4}(z) & -\phi_{1}(z) & \phi_{3}(z)\\ 
	\vartheta \phi_{4}(z) & -\vartheta \phi_{1}(z) & \vartheta \phi_{3}(z)\\
	\vartheta^2 \phi_{4}(z) & -\vartheta^2 \phi_{1}(z) & \vartheta^2 \phi_{3}(z)\end{pmatrix}, & -\pi<\arg (z)<-\frac{\pi}{2}.
\end{cases}
\end{align}
\end{definition}
The definition is such that $\Phi_{\alpha}$ satisfies
a convenient RH problem that will be discussed in the next section.

We conclude this preliminary section by announcing the following
main result. It expresses the limiting kernel 
$\mathbb K^{(\alpha, \frac{1}{2})}$, see \eqref{KVnlimit2} and \eqref{Kalth}, 
in terms of $\Phi_{\alpha}$ and its inverse matrix 
$\Phi_{\alpha}^{-1}$.

\begin{theorem} \label{thmKformula}
We have for $x, y > 0$,
\begin{equation} \label{Kformula} 
	\mathbb K^{(\alpha, \frac{1}{2})}(x,y)
	= 
\frac{1}{2\pi i (x-y)} \begin{pmatrix} -1 & 1 & 0\end{pmatrix}
\Phi_{\alpha,+}^{-1}(y)  \Phi_{\alpha,+}(x) \begin{pmatrix} 1 \\ 1 \\ 0\end{pmatrix},
\end{equation}
\end{theorem}
We prove Theorem \ref{thmKformula} in the next section.

Theorem \ref{mainThm} follows from
a steepest descent analysis of the RH problem for $Y$,
where $\Phi_{\alpha}$ is used to construct a local parametrix
near $0$. In the appropriate scaling the expression \eqref{CorKerY} for
the correlation kernel tends to the kernel \eqref{Kformula}
as given in terms of $\Phi_{\alpha}$.

\section{RH problem for $\Phi_{\alpha}$ and proof of
Theorem \ref{thmKformula}}
\label{sec:RHproblem}

We discuss the RH problem satisfied by $\Phi_{\alpha}$.

\subsection{Jump conditions}
\begin{figure}
\centering
\begin{picture}(300,100)(10,0)
\thicklines
\put(20,50){\line(1,0){200}}
\put(120,0){\line(0,1){100}}
\put(120,50){\circle*{3}}
\put(85,50){\vector(-1,0){3}}
\put(155,50){\vector(1,0){3}}
\put(120,75){\vector(0,1){3}}
\put(120,25){\vector(0,-1){3}}
\put(113,40){$0$}
\end{picture}
\caption{The directions of the contour $\Sigma_{\Phi} = \mathbb R \cup i \mathbb R$ 
used in Lemma \ref{MGPhi} and in the model RH problem.
\label{RandiR}}
\end{figure}
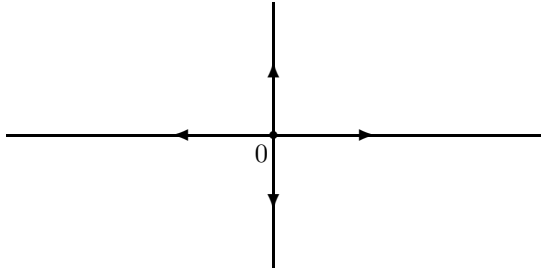

The definition \eqref{defPhialpha} is such that $\Phi_{\alpha}$ has the following
jumps along each of the contours in 
$\Sigma_{\Phi} = \mathbb R \cup i \mathbb R$,
\begin{lemma}
\label{MGPhi}
$\Phi_{\alpha}$ is analytic in each of the four quadrants with constant jumps
(all contours are oriented away from the origin, see Figure \ref{RandiR})
\begin{align} \label{Phialphajump}
	\Phi_{\alpha,+}(z) = \Phi_{\alpha,-}(z)
	\times \begin{cases} 
		\begin{pmatrix} 0 & 1 & 0 \\ -1 & 0 & 0 \\ 0 & 0 & 1 \end{pmatrix},
		&  z \in \mathbb R^+, \\
		\begin{pmatrix} 1 & 0 & 0 \\ 1 & 1 & 0 \\ 0 & 0 & 1 \end{pmatrix},
		& z \in i \mathbb R^{\pm}, \\
		\begin{pmatrix} 1 & 0 & 0 \\
			0 & 0 & i e^{2\pi i \alpha}  \\ 
			0 & -i e^{2\pi i \alpha} & 0 \end{pmatrix},
		& z \in \mathbb R^-
		\end{cases}
		\end{align}
\end{lemma}
\begin{proof} 
We only have
to verify the jump conditions for the first row,  since the jump matrices 
are constant along each of the four contours. Indeed, once we have the jump 
conditions for the first row, the conditions for the other rows will follow
by simple differentiation and multiplying by $z$. 

The jumps on $\mathbb R^+$ and on $i \mathbb R^{\pm}$ follow 
immediately from the definition \eqref{defPhialpha} combined with \eqref{phi4}.
To check the jump on $\mathbb R^-$, we note that $\phi_4$ is
entire because by \eqref{phi4phi0} and \eqref{phi0} it is a multiple
of the $\phantom{}_0 F_2$ hypergeometric function, which is entire.
Thus $\phi_{4,+} = \phi_{4,-}$ on $\mathbb R^-$.
The identities $\phi_{1,+} = e^{2\pi i \beta} \phi_{3,-}$
and $\phi_{3,+} = e^{2\pi i \beta} \phi_{2,-}$ are almost immediate from
\eqref{phi1}--\eqref{phi4}. For example if $z = -x$ with $x > 0$,
then by \eqref{phi1} and \eqref{phi3}
\[ \phi_{1,+}(z) = \phi_{1,+}(x e^{-\pi i}) = 
	e^{2\pi i \beta} G_{\alpha}(xe^{\pi i}) = e^{2\pi i \beta} \phi_{3,-}(z). \]
This shows that the jump conditions \eqref{Phialphajump} on $\mathbb R^-$
are indeed satisfied. 
\end{proof}

\subsection{Behavior as $z \to 0$}

\begin{lemma}
 As $z \to 0$ we have the behavior 
\begin{align} \label{Phiat0} 
	\Phi_{\alpha}(z)
	= \begin{cases} 
	\mathcal O \begin{pmatrix} h_{-\alpha - \frac{1}{2}}(z) & 
	h_{-\alpha - \frac{1}{2}}(z) & h_{-\alpha - \frac{1}{2}}(z) \\
	h_{-\alpha - \frac{1}{2}}(z) & 
	h_{-\alpha - \frac{1}{2}}(z) & h_{-\alpha - \frac{1}{2}}(z) \\
	h_{-\alpha - \frac{1}{2}}(z) & 
	h_{-\alpha - \frac{1}{2}}(z) & h_{-\alpha - \frac{1}{2}}(z) 
	\end{pmatrix}, & \text{ for } \Re z > 0, \\
	 \mathcal O \begin{pmatrix} 1 & 
	h_{-\alpha - \frac{1}{2}}(z) & h_{-\alpha - \frac{1}{2}}(z) \\
	1 & 
	h_{-\alpha - \frac{1}{2}}(z) & h_{-\alpha - \frac{1}{2}}(z) \\
	1  & 
	h_{-\alpha - \frac{1}{2}}(z) & h_{-\alpha - \frac{1}{2}}(z) 
	\end{pmatrix}, & \text{ for } \Re z < 0, 
	\end{cases}
	\end{align}
where $h_{-\alpha-\frac{1}{2}}$ is as in \eqref{RHY4}.
\end{lemma}
\begin{proof}
Because of \eqref{phi4phi0} and \eqref{phi0} the solution
$\phi_4$ is entire and thus bounded at $0$.
By the definition \eqref{defPhialpha} the function $\phi_4$
appears in the first column of $\Phi_{\alpha}$ in the left half-plane.
This accounts for the $\mathcal O(1)$ terms in  \eqref{Phiat0}
as $z \to 0$ with $\Re z < 0$.

The remaining terms come from the behavior of the functions 
$\phi_j$, $j=1,2,3$ as $z \to 0$. These  are solutions of the linear
differential equation \eqref{HGDE}, which can be solved by the
method of Frobenius. 
The indicial equation has roots $0$, $-\alpha$ and $-\alpha-\tfrac{1}{2}$.
If $\alpha > -\tfrac{1}{2}$ then all solutions behave
as $\mathcal O(z^{-\alpha - \frac{1}{2}})$ as $z \to 0$.
Also $\vartheta \phi_j$ and $\vartheta^2 \phi_j$ 
have the same behavior and \eqref{Phiat0} follows.

Similar considerations give \eqref{Phiat0} in case $\alpha = -\tfrac{1}{2}$,
in which case there is a generic $\mathcal O(\log z)$ behavior as $z \to 0$, 
or $-1 < \alpha < -\tfrac{1}{2}$, in which case solutions remain bounded 
as $z\to 0$.
\end{proof}

\begin{remark}
The behavior \eqref{Phiat0} suffices for the purposes of this paper.
However it is possible to obtain more precise information on
the behavior near $0$ by means of connection matrices $C_j$, $j=1,2,3,4$,
that are  such that
\[ \Phi_{\alpha}(z) C_j = 
\mathcal O \begin{pmatrix} 1 & z^{-\alpha} &  z^{-\alpha-\tfrac{1}{2}} \\
 1 & z^{-\alpha} &  z^{-\alpha-\tfrac{1}{2}} \\ 1 & z^{-\alpha} &  z^{-\alpha-\tfrac{1}{2}} \end{pmatrix} \]
 as $z \to 0$ in the $j$th quadrant.  For $\alpha = 0$ or $\alpha= -\frac{1}{2}$ there are 
additional log terms in either the second or the third columns.
\end{remark}

\subsection{Behavior as $z \to \infty$}

The behavior of $\Phi_{\alpha}(z)$ as $z \to \infty$ will be 
deduced from the known asymptotic behavior of Meijer G-functions
given in \cite[Theorem 5, page 179]{Lu}. It gives us 
\begin{align} \label{MeijerGexpansion}
\MeijerG{3}{0}{0}{3}{-}{0,-\alpha,-\alpha-\frac{1}{2}}{z}
\sim 
\frac{2\pi}{\sqrt{3}} z^{-\gamma} e^{-3 z^{\frac{1}{3}}} 
	 \left(1 + \sum_{k=1}^{\infty} M_{k} z^{-\frac{k}{3}} \right),
	\qquad -4 \pi < \arg z < 4 \pi
\end{align}
as $z\to\infty$ with explicit constants 
	(with $\beta = \alpha + \tfrac{1}{4}$ as in \eqref{CorKerY})
\begin{align} \label{defgamma} 
	\gamma & = \frac{2}{3}{\alpha}+\frac{1}{2} = 
	\frac{2}{3} \beta + \frac{1}{3}, \\
	M_1 & = \label{defM1}
	 \frac{1}{3} \alpha^2 + \frac{1}{6} \alpha - \frac{1}{36}, \\
	M_2  & = \label{defM2} 
	\frac{1}{18} \alpha^4 + \frac{1}{54}\alpha^3 - \frac{17}{216} \alpha^2 - \frac{1}{54} \alpha + \frac{25}{2592},
	\end{align}
	and other $M_k$ can be calculated, if necessary.

\begin{lemma} Let $\beta = \alpha + \tfrac{1}{4}$
and $\omega = e^{2\pi i/3}$. Then we have as $z \to \infty$,
\begin{multline}
\label{asympPhi}
{T}_{\alpha}
\Phi_{\alpha}(z) 
 = \frac{2\pi}{\sqrt{3}} z^{-\frac{2\beta}{3}} 
 \left( \mathbb{I} + A_{\alpha} z^{-1} +
 	\mathcal O(z^{-2}) \right)
 	L_{\alpha}(z) \times
\begin{cases}
\begin{pmatrix}  e^{-3\omega z^{1/3}} & 0 & 0\\ 0 & e^{-3\omega^2 z^{1/3}} & 0\\ 0 & 0 & e^{-3 z^{1/3}} \end{pmatrix}, & \Im(z)>0, \\
\begin{pmatrix} e^{-3\omega^2 z^{1/3}} & 0 & 0\\ 0 & e^{-3\omega z^{1/3}} & 0\\ 0 & 0 & e^{-3 z^{1/3}} \end{pmatrix}, & \Im(z)<0,
\end{cases}
\end{multline}
where 
\begin{equation} \label{defLalpha}
	L_{\alpha}(z) = \begin{pmatrix} z^{-\frac{1}{3}} & 0 & 0 \\
		0 & 1 & 0 \\ 0 & 0 & z^{\frac{1}{3}} \end{pmatrix}
		\times
	\begin{cases} \begin{pmatrix} 	\omega^2 & \omega & 1\\ 
	1 & 1 & 1\\ \omega & \omega^2 & 1 \end{pmatrix}  
	\begin{pmatrix} e^{\frac{2\pi i \beta}{3}} & 0 & 0 \\
	0 & e^{-\frac{2\pi i \beta}{3}} & 0 \\
	0 & 0 & 1 \end{pmatrix}, &\Im(z)>0,\\
	\begin{pmatrix}  \omega & -\omega^2 & 1\\ 1 & -1 & 1\\ 
	\omega^2 & -\omega & 1\end{pmatrix}
	\begin{pmatrix} e^{-\frac{2\pi i \beta}{3}} & 0 & 0 \\
	0 & e^{\frac{2\pi i \beta}{3}} & 0 \\
	0 & 0 & 1 \end{pmatrix},
	 &\Im(z)<0
\end{cases}
\end{equation}
$A_{\alpha}$ is a certain constant matrix and
\begin{equation} \label{defTalpha}
T_{\alpha} = \begin{pmatrix} 1 & 0 & 0 \\ t_1 & -1 & 0 \\
t_2 & t_3 & 1 \end{pmatrix} 
\end{equation}
is a lower triangular matrix with entries
\begin{equation} \label{defts}
	t_1 = -\gamma-M_1, \quad
	t_2 = \gamma(\gamma- \tfrac{1}{3}) + M_1(M_1 + \gamma- \tfrac{2}{3}) - M_2, \quad
	t_3  = 2 \gamma - \tfrac{1}{3} + M_1
	\end{equation}
where $\gamma$, $M_1$ and $M_2$ are given in \eqref{defgamma},
\eqref{defM1}, \eqref{defM2}.
\end{lemma}

\begin{proof}
Let us focus on the last column of $\Phi_{\alpha}$, see
\eqref{phi3} and \eqref{defPhialpha}, with $-\pi < \arg z < \pi$. 
Then by \eqref{phi3} and \eqref{MeijerGexpansion} we have
\begin{equation} \label{phi3expansion} 
	\phi_3(z) = \frac{2\pi}{\sqrt{3}}
	z^{-\gamma} e^{-3 z^{\frac{1}{3}}}  
		\left(1 + M_1 z^{-\frac{1}{3}} + M_2 z^{-\frac{2}{3}}
		+ \mathcal{O}(z^{-1}) \right) \end{equation}	
The expansion \eqref{MeijerGexpansion} can be differentiated 
termwise and we find
\begin{align} \label{dphi3expansion}
-\vartheta \phi_3(z)
	& = \frac{2\pi}{\sqrt{3}} z^{-\gamma}
	e^{-3 z^{\frac{1}{3}}} 
		\left(z^{\frac{1}{3}} + (M_1+\gamma) + \left(M_2 + 
		(\gamma + \tfrac{1}{3})M_1 \right) z^{-\frac{1}{3}} 
		+ \mathcal{O}(z^{-\frac{2}{3}}) \right) \\
	\vartheta^2 \phi_3(z)
	& = 	 \label{ddphi3expansion}
 	\frac{2\pi}{\sqrt{3}} z^{-\gamma} e^{-3 z^{\frac{1}{3}}} 
 	 \left(z^{\frac{2}{3}} + 
 \left(M_{1}+2\gamma-\tfrac{1}{3}\right) z^{\frac{1}{3}} 
 + \left(M_{2}+(2 \gamma + \tfrac{1}{3})M_{1} + \gamma^2\right) 
 +\mathcal{O}(z^{-\frac{1}{3}})\right).
\end{align}
If we choose constants  $t_1$, $t_2$, $t_3$ as in \eqref{defts}
then we get from \eqref{phi3expansion}, \eqref{dphi3expansion},
\eqref{ddphi3expansion},
\begin{align} \label{tdphiexpansion}
 t_1 \phi_3(z) -\vartheta \phi_3(z)
	& = \frac{2\pi}{\sqrt{3}}
	z^{-\gamma} e^{-3 z^{\frac{1}{3}}}  
		\left(z^{\frac{1}{3}} + \mathcal{O}(z^{-\frac{1}{3}}) \right), \\
	 t_2 \phi_3(z) + 
	t_3 \vartheta \phi_3(z) + \vartheta^2 \phi_3(z) 	& = 	 
 	\frac{2\pi}{\sqrt{3}} z^{-\gamma} e^{-3 z^{\frac{1}{3}}} 
 	 \left(z^{\frac{2}{3}} +  
  \mathcal{O}(z^{-\frac{1}{3}})\right). \label{tddphiexpansion}
\end{align}

With the matrix $T_{\alpha}$ from \eqref{defTalpha} we write 
\eqref{phi3expansion}, \eqref{tdphiexpansion} and 
\eqref{tddphiexpansion} as 
\[ T_{\alpha} \begin{pmatrix} \phi_3(z) \\ \vartheta \phi_3(z) \\ 
	\vartheta^2 \phi_3(z) \end{pmatrix} 	
	= \frac{2\pi}{\sqrt{3}} z^{-\gamma + \frac{1}{3}}
		e^{-3 z^{\frac{1}{3}}}
		\begin{pmatrix} z^{-\frac{1}{3}} + \mathcal{O}(z^{-\frac{2}{3}}) \\	
	  1 + \mathcal{O}(z^{-\frac{2}{3}}) \\
		z^{\frac{1}{3}} + \mathcal{O}(z^{-\frac{2}{3}})
		\end{pmatrix},  \]
which may also be written as
\begin{equation}
\label{asympPhiColumn}
 T_{\alpha} \Phi_{\alpha}(z) \begin{pmatrix} 0 \\ 0 \\ 1 \end{pmatrix}	
 	= \frac{2\pi}{\sqrt{3}} z^{-\gamma + \frac{1}{3}}
	\left( \mathbb{I} + \mathcal O(z^{-1})\right)
		L_{\alpha}(z) \begin{pmatrix} 0 \\ 0 \\ 1 \end{pmatrix}
		e^{-3 z^{\frac{1}{3}}}
\end{equation}
by the definitions \eqref{defLalpha} of $L_{\alpha}$ and \eqref{defPhialpha}
of $\Phi_{\alpha}$.
This is the equality  \eqref{asympPhi} for the third column with
$\mathbb I + A_{\alpha} z^{-1} + \mathcal{O}(z^{-2})$ 
replaced by $\mathbb{I} + \mathcal{O}(z^{-1})$. 
The  equality of the other two columns follows in a similar manner, and
with the same matrix $T_{\alpha}$, although the details are a bit more involved.
We  apply \eqref{MeijerGexpansion} with $z e^{2\pi i}$ and $z e^{-2\pi i}$ instead of $z$, 
and this leads to  the various factors of $\omega$ and $\omega^2$ in
\eqref{asympPhi} and \eqref{defLalpha}. The result is that \eqref{asympPhi} holds
with $\mathbb{I} + \mathcal{O}(z^{-1})$ instead of $\mathbb I + A_{\alpha} z^{-1} + \mathcal{O}(z^{-2})$.

To obtain the sharper order estimate, we consider 
\begin{align}
\label{defLhat}
\widehat L_{\alpha}(z) 
 = \frac{2\pi}{\sqrt{3}} z^{-\frac{2\beta}{3}} 
  	L_{\alpha}(z) \times
\begin{cases}
\begin{pmatrix}  e^{-3\omega z^{1/3}} & 0 & 0\\ 0 & e^{-3\omega^2 z^{1/3}} & 0\\ 0 & 0 & e^{-3 z^{1/3}} \end{pmatrix}, & \Im(z)>0, \\
\begin{pmatrix} e^{-3\omega^2 z^{1/3}} & 0 & 0\\ 0 & e^{-3\omega z^{1/3}} & 0\\ 0 & 0 & e^{-3 z^{1/3}} \end{pmatrix}, & \Im(z)<0,
\end{cases}
\end{align}
and we verify that $\widehat{L}_{\alpha}$ has the same jumps as $\Phi_\alpha$ 
has on $\mathbb{R}^\pm$.
The details of this calculation are similar to what we will do in the
the proof of Lemma \ref{lemEninv}(a) below. This lemma is about a slightly different function, $F_n(z)$, but the arguments are essentially the same. 
It follows that 
\begin{equation} \label{defRandO} 	
 R_{\alpha} = T_\alpha \Phi_\alpha \widehat{L}_\alpha^{-1} = \mathbb{I} + \mathcal{O}(z^{-1}) 
 \end{equation}
 has no jumps on $\mathbb{R}^\pm$, and therefore is analytic in $\mathbb C \setminus i \mathbb R$.
 
 Note that $\widehat{L}_{\alpha}$ is analytic across the imaginary axis, while
 $\Phi_{\alpha}$ has the jump \eqref{Phialphajump} there. For $R_{\alpha}$ we then have the
 jump matrix
 \[ R_{\alpha,-}^{-1}(z) R_{\alpha,+}(z) =
 	  \widehat{L}_{\alpha}(z) \begin{pmatrix} 1 & 0 & 0 \\ 1 & 1 & 0 \\ 0 & 0 & 1 \end{pmatrix}	
 	  \left(\widehat{L}_{\alpha}(z)\right)^{-1}, \quad z \in i \mathbb R,
 \]
 which by \eqref{defLhat} is equal to
 \begin{align*}  R_{\alpha,-}^{-1}(z) R_{\alpha,+}(z) & =
 		L_{\alpha}(z) \begin{pmatrix} 1 & 0 & 0 \\ e^{\pm 3(\omega - \omega^2) z^{1/3}}  & 1 & 0 \\ 0 & 0 & 1 \end{pmatrix}	
 	L_{\alpha}^{-1}(z),   	\quad z \in i \mathbb R^{\pm}, \\
 	&  	 = 
 	L_{\alpha}(z) \begin{pmatrix} 1 & 0 & 0 \\ e^{- \frac{3 \sqrt{3}}{2} |z|^{1/3}}  & 1 & 0 \\ 0 & 0 & 1 \end{pmatrix}	
 	L_{\alpha}^{-1}(z).
 \end{align*}
 The entries of $L_{\alpha}(z)$ and $L_{\alpha}^{-1}(z)$ are $\mathcal{O}(z^{1/3})$ as $z \to \infty$,
 as can be seen from \eqref{defLalpha}.
 Hence there is a constant $c > 0$ such that 
 $R_{\alpha,+}(z) = R_{\alpha,-}(z) \left(\mathbb{I} + \mathcal O(e^{-c |z|^{1/3}}) \right)$
 as $z \to \infty$ along the imaginary axis. 
 
 It follows from this and \eqref{defRandO} that  the $\mathcal{O}$ term in
  \eqref{defRandO} can be written as an asymptotic series in powers of $z^{-1}$
  as $z \to \infty$.
  In particular we find that there is a constant matrix $A_{\alpha}$ such that
  $R_{\alpha}(z) = \mathbb{I} + A_{\alpha} z^{-1} + \mathcal{O}(z^{-2})$ and the lemma follows.
\end{proof}

\subsection{RH problem for $\Phi_{\alpha}$}

We combine the above to conclude that $\Phi_{\alpha}$
satisfies the following RH problem.

\begin{rhproblem} \label{RHPforPhi} \
\begin{description}
\item[RH-$\Phi$1] $\Phi_{\alpha}$ is analytic in each of the four quadrants. 
\item[RH-$\Phi$2] $\Phi_{\alpha}$ has the constant jumps 
\eqref{Phialphajump} on $\Sigma_{\Phi} = \mathbb R \cup i \mathbb R$.
\item[RH-$\Phi$3] $\Phi_{\alpha}$ has the asymptotic behavior \eqref{asympPhi} as $z \to \infty$.
\item[RH-$\Phi$4] $\Phi_{\alpha}$ has the behavior \eqref{Phiat0} as $z \to 0$.
\end{description}
\end{rhproblem}

With standard RH arguments it can be shown that the solution
to the RH problem is unique (and thus is given by \eqref{defPhialpha}),
and that
\[ \det \Phi_{\alpha}(z) =  
	- \left(\frac{2\pi}{\sqrt{3}} \right)^3 z^{-2\beta}, \]
where $\beta = \alpha + \frac{1}{4}$ as before.

\begin{remark}
	In their analysis of coupled random matrices in a $p$-chain,
	Bertola and Bothner \cite{BeBo} also construct a model RH problem for Meijer G-functions.
 	 It is  called Bare Meijer G-parametrix in \cite[section 4.2.1]{BeBo}. For $p=2$ the 
	RH problem as well as its solution in terms of Mellin-Barnes type integrals shows
	great similarity with RH problem \ref{RHPforPhi}, although there does not seem to be a direct way
	to connect the two RH problems. 
\end{remark}

\subsection{The inverse of $\Phi_{\alpha}$}

The inverse of $\Phi_{\alpha}$ also appears in the
formula for  the limiting kernel 
$\mathbb K^{(\alpha,\frac{1}{2})}(x,y)$ in Theorem \ref{thmKformula}.
The inverse contains the Meijer G-functions 
$\MeijerG{3}{0}{0}{3}{-}{0, \alpha, \alpha+ \frac{1}{2}}{-z}$.
We define, for $-\pi < \arg(z) <\pi$,
\begin{align} \label{psi1}
\psi_{1}(z) &=
\MeijerG{3}{0}{0}{3}{-}{0,\alpha,\alpha+\frac{1}{2}}{z e^{-\pi i}}, \\
\psi_{2}(z) &= \label{psi2}
\MeijerG{3}{0}{0}{3}{-}{0,\alpha,\alpha+\frac{1}{2}}{z e^{\pi i}}, \\
\psi_{3}(z) &=  \label{psi3}
-i e^{2\pi i\alpha}
\MeijerG{3}{0}{0}{3}{-}{0,\alpha,\alpha+\frac{1}{2}}{z e^{-\pi i}} +
i e^{2\pi i\alpha}
\MeijerG{3}{0}{0}{3}{-}{0,\alpha,\alpha+\frac{1}{2}}{z e^{-3\pi i}}, \\
&=  \label{psi3bis}
	i e^{-2\pi i\alpha}
\MeijerG{3}{0}{0}{3}{-}{0,\alpha,\alpha+\frac{1}{2}}{z e^{\pi i}} 
- i e^{-2\pi i\alpha}
\MeijerG{3}{0}{0}{3}{-}{0,\alpha,\alpha+\frac{1}{2}}{z e^{3\pi i}}, \\
\psi_4(z) &= \psi_2(z) - \psi_1(z). \label{psi4}
\end{align}
The functions $\psi_j$ are solutions of the third order
differential equation \eqref{HGDE}, but with $\alpha$ replaced
by $-\alpha - \frac{1}{2}$, and $z$ replaced by $-z$, i.e.,
\begin{equation} \label{adjHGDE} 
	\vartheta(\vartheta-\alpha)(\vartheta-\alpha-\tfrac{1}{2}) \psi
	- z \psi = 0. \end{equation}
This differential equation has Frobenius indices $0,\alpha$ and $\alpha+\frac{1}{2}$. When $2\alpha\not\in \mathbb{Z}$, \eqref{adjHGDE} has a basis of solutions around $z=0$ formed by the functions
\begin{align} \label{psiFrob}
\ds {{}_{0}F_{2}}\left( {- \atop 1 - \alpha, \frac{3}{2}-\alpha} ; z\right),
\ds z^{\alpha} {{}_{0}F_{2}}\left( {- \atop 1 + \alpha, \frac{3}{2}} ; z\right)
\text{ and }
\ds z^{\alpha+\frac{1}{2}} {{}_{0}F_{2}}\left( {- \atop \frac{1}{2}+\alpha, \frac{1}{2}} ; z\right).
\end{align}
In particular, there are (real) constants $c_1,c_2,c_3$ such that for $-\pi<\arg(z)<\pi$
\begin{multline} \nonumber
\psi_1(z) =
c_1\hspace{0.1cm}\ds {{}_{0}F_{2}}\left( {- \atop 1 - \alpha, \frac{3}{2}-\alpha} ; z\right)+
c_2 e^{-\pi i\alpha}\ds z^{\alpha} {{}_{0}F_{2}}\left( {- \atop 1 + \alpha, \frac{3}{2}} ; z\right)
- c_3 i e^{-\pi i\alpha} \ds z^{\alpha+\frac{1}{2}} {{}_{0}F_{2}}\left( {- \atop \frac{1}{2}+\alpha, \frac{1}{2}} ; z\right).
\end{multline}
By analytically continuing this expression along circles we can write the other $\psi_j$ in terms of the Frobenius basis. Note that, regardless of whether we use \eqref{psi3} or \eqref{psi3bis}, the Frobenius basis \eqref{psiFrob} yields that
\begin{equation} \label{psi3Frob}
\psi_3(z) = 2 c_2 \sin(\pi\alpha) \ds z^{\alpha} {{}_{0}F_{2}}\left( {- \atop 1 + \alpha, \frac{3}{2}} ; z\right)
- 2 c_3 \cos(\pi\alpha) \ds z^{\alpha+\frac{1}{2}} {{}_{0}F_{2}}\left( {- \atop \frac{1}{2}+\alpha, \frac{1}{2}} ; z\right),
\end{equation}
showing that the two definitions \eqref{psi3} and \eqref{psi3bis} 
are indeed the same. By continuity, the equality also holds
true if $2\alpha \in \mathbb Z$.  
	
\begin{definition}
We define, with $\vartheta = z \frac{d}{dz}$ as before,
\begin{align} \label{defPsialpha}
\Psi_{\alpha}(z) =  
\begin{cases}
\begin{pmatrix} \vartheta^2 \psi_{1}(z) & \vartheta^2 \psi_{2}(z)
	& 	\vartheta^2 \psi_{3}(z) \\ 
	\vartheta \psi_{1}(z) & \vartheta \psi_{2}(z) & \vartheta \psi_{3}(z) \\
	\psi_{1}(z) & \psi_{2}(z) &  \psi_{3}(z) \end{pmatrix}, 
	& 0 < \arg (z)< \frac{\pi}{2}, \\
\begin{pmatrix} \vartheta^2 \psi_{1}(z) & \vartheta^2 \psi_{4}(z)
	& 	\vartheta^2 \psi_{3}(z) \\ 
	\vartheta \psi_{1}(z) & \vartheta \psi_{4}(z) & \vartheta \psi_{3}(z) \\
	\psi_{1}(z) & \psi_{4}(z) &  \psi_{3}(z) \end{pmatrix}, 
	& \frac{\pi}{2} < \arg(z) < \pi, \\
\begin{pmatrix} \vartheta^2 \psi_{2}(z) & -\vartheta^2 \psi_{1}(z)
	& 	\vartheta^2 \psi_{3}(z) \\ 
	\vartheta \psi_{2}(z) & -\vartheta \psi_{1}(z) & \vartheta \psi_{3}(z) \\
	\psi_{2}(z) & -\psi_{1}(z) &  \psi_{3}(z) \end{pmatrix}, 
	& -\frac{\pi}{2} < \arg(z) < 0, \\
\begin{pmatrix} \vartheta^2 \psi_{2}(z) & \vartheta^2 \psi_{4}(z)
	& 	\vartheta^2 \psi_{3}(z) \\ 
	\vartheta \psi_{2}(z) & \vartheta \psi_{4}(z) & \vartheta \psi_{3}(z) \\
	\psi_{2}(z) & \psi_{4}(z) &  \psi_{3}(z) \end{pmatrix}, 
	& -\pi<\arg (z)<-\frac{\pi}{2}.
\end{cases}
\end{align}
\end{definition}

The main property of $\Psi_{\alpha}$ is contained in
\begin{lemma}
We have
\begin{equation} \label{PhiInv}
\Phi_{\alpha} \Psi_{\alpha}^T = - 4 \pi^2
	\begin{pmatrix} 1 & 0 & 0 \\ -2\alpha-\tfrac{1}{2} & -1 & 0
	\\ \alpha(\alpha+\tfrac{1}{2}) & 2\alpha + \tfrac{1}{2} & 1
	\end{pmatrix}^{-1} 
\end{equation}
\end{lemma}

\begin{proof}
It is clear that $\Psi_{\alpha}$ is analytic in each
of the four quadrants with jumps 
\begin{equation} \label{Psijumps}
	\Psi_{\alpha,+}(z) = 
	\Psi_{\alpha,-}(z) \times
	\begin{cases} \begin{pmatrix} 0 & 1 & 0 \\ -1 & 0 & 0 \\ 0 & 0 & 1
	\end{pmatrix}, & z \in  \mathbb R^+, \\
	\begin{pmatrix} 1 & -1 & 0 \\ 0 & 1 & 0 \\ 0 & 0 & 1
	\end{pmatrix}, & z \in i \mathbb R^{\pm}, \\
	\begin{pmatrix} 1 & 0 & 0 \\ 0 & 0 & -i e^{-2\pi i \alpha} \\
	0 & i e^{-2\pi i \alpha} & 0 \end{pmatrix},
	& z \in \mathbb R^-. \end{cases} \end{equation}
Indeed, the jumps on $\mathbb R^+$ and $i\mathbb R^{\pm}$
are immediate from \eqref{defPsialpha}.
To verify the jump on $\mathbb R^-$ 
we have to check that on $\mathbb R^-$
\[ 	\psi_{2,+}  = \psi_{1,-}, \quad
	\psi_{3,+}  = i e^{-2\pi i \alpha} \psi_{4,-}, \quad
	\psi_{4,+}  = -i e^{-2\pi i \alpha} \psi_{3,-}.
	\]
These identities are indeed satisfied because of the definitions
	\eqref{psi1}--\eqref{psi4}. We  use that
	$z_- = z_+ e^{2\pi i}$ for $z \in \mathbb R^-$,
	since the negative real axis is oriented from right to left.

The jump matrices in \eqref{Psijumps} are the inverse transposes
of those in \eqref{Phialphajump}.
It follows that $\Phi_{\alpha} \Psi_{\alpha}^{T}$
is analytic across all contours, and so is analytic in 
$\mathbb C \setminus \{0\}$. The isolated singularity
at $0$ is at most a pole, since all functions involved
in $\Phi_{\alpha}$ and $\Psi_{\alpha}$ have a power law behavior
near $0$.

The functions $\psi_j$ are solutions of \eqref{adjHGDE}
which has Frobenius indices $0$, $\alpha$ and 
$\alpha + \frac{1}{2}$ at $0$. As such they are linear combinations of the functions \eqref{psiFrob} when $2\alpha\not\in\mathbb{Z}$.
If $\alpha < 0$ then all solutions are $\mathcal O(z^{\alpha})$
as $z \to 0$, and it follows from \eqref{defPsialpha} that
$\Psi_{\alpha}(z) = \mathcal O(z^{\alpha})$.
By \eqref{Phiat0} one has 
$\Phi_{\alpha}(z) = \mathcal O(z^{-\alpha- \frac{1}{2}})$
for $- \frac{1}{2} < \alpha < 0$ or $\Phi_{\alpha} = \mathcal O(1)$
for $\alpha < -\frac{1}{2}$. Thus 
\begin{equation} \label{PhiPsiat0} 
	\Phi_{\alpha}(z) \Psi_{\alpha}^T(z) 
	= \mathcal O\left(z^{ \min(\alpha, - \frac{1}{2})}\right)
		\qquad \text{ if } \alpha \in (-1, 0) \setminus
		\{ -\tfrac{1}{2}\}. \end{equation}
For $\alpha = -\frac{1}{2}$ or $\alpha = 0$ there is an
additional log term in the $\mathcal O$ term in \eqref{PhiPsiat0}.
Since $\alpha >-1$, we see from \eqref{PhiPsiat0}
that the isolated singularity is removable
in case $\alpha \leq 0$.

For $\alpha > 0$ we have to be a little more careful.
Now all solutions are bounded at $0$, but  we 
use that $\psi_3$ and $\psi_4$ are solutions of \eqref{adjHGDE}
which are $\mathcal O(z^{\alpha})$ as $z \to 0$. This fact can be checked by writing the solutions out in the Frobenius basis around $z=0$. 
We already did this for $\psi_3$ in \eqref{psi3Frob} and the situation for $\psi_4$ is analogous. 
Then it follows from \eqref{defPsialpha} that
\begin{equation} \label{Psiat0} \Psi_{\alpha}(z) = 
	\mathcal O \begin{pmatrix} 1 & z^{\alpha} & z^{\alpha} \\
	1 & z^{\alpha} & z^{\alpha} \\
	1 & z^{\alpha} & z^{\alpha} \end{pmatrix} 
	\end{equation}
as $z \to 0$ with $\Re z < 0$.
Thus
\begin{align*} 
	\Phi_{\alpha}(z) \Psi_{\alpha}^T(z) & = \mathcal O\begin{pmatrix}
	1 & z^{-\alpha-\frac{1}{2}} & z^{-\alpha-\frac{1}{2}} \\
	1 & z^{-\alpha-\frac{1}{2}} & z^{-\alpha-\frac{1}{2}} \\
	1 & z^{-\alpha-\frac{1}{2}} & z^{-\alpha-\frac{1}{2}} \end{pmatrix}
	\mathcal O\begin{pmatrix}
	1 & 1 & 1 \\
	z^{\alpha} & z^{\alpha} & z^{\alpha} \\
	z^{\alpha} & z^{\alpha} & z^{\alpha} 
	 \end{pmatrix} = \mathcal O\begin{pmatrix} 
	 z^{-\frac{1}{2}} & z^{-\frac{1}{2}} &z^{-\frac{1}{2}} \\
	 z^{-\frac{1}{2}} & z^{-\frac{1}{2}} &z^{-\frac{1}{2}} \\
	 z^{-\frac{1}{2}} & z^{-\frac{1}{2}} &z^{-\frac{1}{2}} 
	\end{pmatrix}
\end{align*}
as $z \to 0$ with $\Re z < 0$. The singularity of $\Phi_{\alpha} \Psi_{\alpha}^T$ at $0$ is at most
a pole, and thus it must be removable, since if it were
a pole the absolute value of some entry would be 
$\geq \frac{C}{|z|}$ for some $C > 0$ and all $z$ close
enough to $0$.
We conclude that
$\Phi_{\alpha} \Psi_{\alpha}^T$ is entire. 
	
Next we investigate the behavior as $z \to \infty$. The asymptotic behavior
of $\Phi_{\alpha}$ is given in  \eqref{asympPhi}. With similar arguments
we find the behavior of $\Psi_{\alpha}$.
We use \eqref{MeijerGexpansion}, \eqref{defgamma}, \eqref{defM1}, \eqref{defM2} 
with $\alpha$ replaced by $-\alpha-\tfrac{1}{2}$ and find
\begin{align} \label{MeijerGexpansionInverse}
\MeijerG{3}{0}{0}{3}{-}{0,\alpha,\alpha+\frac{1}{2}}{z}
\sim 
\frac{2\pi}{\sqrt{3}} e^{-3 z^{\frac{1}{3}}} 
z^{-\tilde\gamma} \left(1+ \sum_{k=1}^{\infty} \tilde M_{k} z^{-\frac{k}{3}}
\right),
	\qquad -4 \pi < \arg z < 4 \pi
\end{align}
with
\begin{align} \label{gammatilde}
	\tilde{\gamma} & = -\frac{2}{3}\alpha + \frac{1}{6} = 
	- \frac{2}{3} \beta + \frac{1}{3}, \\ \label{M1tilde}
	\tilde{M}_1  & = \frac{1}{3} \alpha^2 + \frac{1}{6} \alpha - \frac{1}{3}, \\ \label{M2tilde}
	\tilde{M}_2 & = \frac{1}{18} \alpha^4 + \frac{5}{54} \alpha^3 - \frac{5}{216} \alpha^2 - \frac{5}{108} \alpha + \frac{1}{2592}.
	\end{align}
	
We then have for $\psi_3$ 
\begin{align}
\psi_{3}(z) &= -\frac{2\pi}{\sqrt{3}} z^{-\tilde{\gamma}} 
e^{3 z^{\frac{1}{3}}} 
\left(1 - \tilde{M}_{1} z^{-\frac{1}{3}} + \tilde{M}_{2} z^{-\frac{2}{3}}+\mathcal{O}(z^{-1})\right), \\
\vartheta \psi_{3}(z) &= - \frac{2\pi}{\sqrt{3}} z^{-\tilde{\gamma}} 
e^{3 z^{\frac{1}{3}}} 
\left(z^{\frac{1}{3}} - (\tilde{M}_{1} + \tilde{\gamma}) 
+ (\tilde{M}_{2}+(\tfrac{1}{3}+\tilde\gamma) \tilde{M}_{1}) 
z^{-\frac{1}{3}}+\mathcal{O}(z^{-\frac{2}{3}}) \right), \\ 
\vartheta^2 \psi_{3}(z) &= - \frac{2\pi}{\sqrt{3}} z^{-\tilde\gamma}
e^{3 z^{\frac{1}{3}}}  \left(z^{\frac{2}{3}} + 
(-2\tilde\gamma+\tfrac{1}{3}-\tilde{M}_{1}) z^{\frac{1}{3}} 
+ (\tilde{M}_{2}+(\tfrac{1}{3} +2\tilde\gamma) M_{1} + \tilde\gamma^2) +\mathcal{O}(z^{-\frac{1}{3}}) 
\right).
\end{align}
These calculations come from applying \eqref{MeijerGexpansionInverse}
to the second term in either \eqref{psi3} or \eqref{psi3bis},
since this term is dominant in the full range $-\pi < \arg z < \pi$.

Then with numbers analogous to \eqref{defts}
\begin{align} \label{deftstilde}
\tilde{t}_{1} = \tilde{\gamma} + \tilde M_{1}, \quad 
\tilde{t}_{2} = \tilde\gamma (\tilde\gamma -\tfrac{1}{3}) + 
\tilde M_{1} (\tilde M_{1}+\tilde\gamma-\tfrac{2}{3}) - \tilde{M}_{2},
\quad
\tilde{t}_{3} = 2 \tilde{\gamma} - \tfrac{1}{3} + \tilde M_{1},
\end{align}
and
\begin{align} \label{Talphatilde}
	\tilde{T}_{\alpha} = \begin{pmatrix} 1 & \tilde{t}_3 & \tilde{t}_2 \\
	 0 & 1 & \tilde{t}_1 \\ 0 & 0 & 1 \end{pmatrix}
\end{align}
we have
\[ 
\tilde{T}_{\alpha} \begin{pmatrix} \vartheta^2 \psi_3(z) \\ 
\vartheta \psi_3(z) \\ \psi_3(z) \end{pmatrix}
	= -\frac{2 \pi}{\sqrt{3}} z^{\frac{2\beta}{3}} e^{3z^{\frac{1}{3}}}
	\begin{pmatrix} z^{\frac{1}{3}} + \mathcal O(z^{-\frac{2}{3}}) \\
	1 + \mathcal O(z^{-\frac{2}{3}}) \\
	 z^{-\frac{1}{3}} + \mathcal O(z^{-\frac{2}{3}})
	 \end{pmatrix}. \]
After similar calculations for $\psi_1, \psi_2, \psi_4$ we find
\begin{multline}
\label{asympPsi}
\tilde{T}_{\alpha} 
\Psi_{\alpha}(z) 
 = - \frac{2\pi}{\sqrt{3}} z^{\frac{2\beta}{3}} 
 	\left( \mathbb{I} + \mathcal{O}(z^{-1}) \right)
 	\tilde{L}_{\alpha}(z) \times
\begin{cases}
\begin{pmatrix}  e^{3\omega z^{1/3}} & 0 & 0\\ 0 & e^{3\omega^2 z^{1/3}} & 0\\ 0 & 0 & e^{3 z^{1/3}} \end{pmatrix}, & \Im(z)>0, \\
\begin{pmatrix} e^{3\omega^2 z^{1/3}} & 0 & 0\\ 0 & e^{3\omega z^{1/3}} & 0\\ 0 & 0 & e^{3 z^{1/3}} \end{pmatrix}, & \Im(z)<0,
\end{cases}
\end{multline}
with 
\begin{equation} \label{defLalphatilde}
	\tilde{L}_{\alpha}(z) = 
	\begin{pmatrix} z^{\frac{1}{3}} & 0 & 0 \\
		0 & 1 & 0 \\ 0 & 0 & z^{-\frac{1}{3}} \end{pmatrix}
		\times
	\begin{cases} \begin{pmatrix} 	\omega & \omega^2 & 1\\ 
	1 & 1 & 1\\ \omega^2 & \omega & 1 \end{pmatrix}  
	\begin{pmatrix} e^{-\frac{2\pi i \beta}{3}} & 0 & 0 \\
	0 & e^{\frac{2\pi i \beta}{3}} & 0 \\
	0 & 0 & 1 \end{pmatrix}, &\Im(z)>0,\\
	\begin{pmatrix}  \omega^2 & -\omega & 1\\ 1 & -1 & 1\\ 
	\omega & -\omega^2 & 1\end{pmatrix}
	\begin{pmatrix} e^{\frac{2\pi i \beta}{3}} & 0 & 0 \\
	0 & e^{-\frac{2\pi i \beta}{3}} & 0 \\
	0 & 0 & 1 \end{pmatrix},
	 &\Im(z)<0.
\end{cases}
\end{equation}

For $z \to \infty$ we then have from \eqref{asympPhi} and \eqref{asympPsi},
\begin{align} \nonumber
	T_{\alpha} \Phi_{\alpha}(z) \Psi_{\alpha}^T(z) 
	\tilde{T}_{\alpha}^T 
	& = T_{\alpha} \Phi_{\alpha}(z)
		\left( 	\tilde{T}_{\alpha}  \Psi_{\alpha}(z) \right)^T \\ 
	& = - \left(\frac{2\pi}{\sqrt{3}} \right)^2	
	\left(\mathbb{I} +   \mathcal O(z^{-\frac{1}{3}}) \right)
	L_{\alpha}(z) \tilde{L}_{\alpha}^T(z)
	\left(\mathbb{I}  + \mathcal O(z^{-\frac{1}{3}}) \right)^T  \label{PhiPsiasymp}
\end{align}
By multiplying \eqref{defLalpha} and \eqref{defLalphatilde} we get
\begin{equation} \label{LLinv} 
	L_{\alpha}(z) \tilde{L}_{\alpha}^T(z) = 3 \mathbb I. \end{equation}
Thus \eqref{PhiPsiasymp} remains bounded as $z \to \infty$.
Since $\Phi_{\alpha}(z) \Psi_{\alpha}^T(z)$ is entire, \eqref{PhiPsiasymp} 
is a  constant matrix by Liouville's theorem, and using \eqref{LLinv} 
we conclude
\begin{equation} \label{PhiPsi} 
	T_{\alpha} \Phi_{\alpha}(z) \Psi_{\alpha}^T(z) \tilde{T}_{\alpha}^T 
	= - 4 \pi^2 \mathbb {I} \end{equation}
for all $z \in \mathbb C$. A final computation (that we checked
with Maple), based on \eqref{defTalpha},
\eqref{Talphatilde}, and \eqref{defts}, \eqref{deftstilde} shows that
\[ \tilde{T}_{\alpha}^T T_{\alpha}
	= \begin{pmatrix} 1 & 0 & 0 \\ \tilde{t}_3 & 1 & 0 \\
	\tilde{t}_2 & \tilde{t}_1 & 1 \end{pmatrix}
	\begin{pmatrix} 1 & 0 & 0 & \\ t_1 & -1 & 0 \\
	t_2 & t_3 & 1 \end{pmatrix}
	= \begin{pmatrix} 1 & 0 & 0 \\ -2\alpha-\tfrac{1}{2} & -1 & 0 \\
	\alpha(\alpha+\tfrac{1}{2}) & 2 \alpha + \tfrac{1}{2} & 1 \end{pmatrix}
	\]
which by \eqref{PhiPsi} leads to the identity in \eqref{PhiInv}.
\end{proof}

\subsection{Proof of Theorem \ref{thmKformula}}

\begin{proof} We are going to use \cite[Proposition 5.4]{KuZh}
with parameters $M=2$ and $\nu_0 = 0, \nu_1 = \alpha$, and
$\nu_2 = \alpha+\frac{1}{2}$.
In the paper \cite{KuZh} the parameters $\nu_j$ are supposed
to be integers, but this fact does not play a role in the proof of
the formula (5.13) of \cite{KuZh}, and we can use it in our situation. The result is that
\begin{equation} \label{Kformula1} 
	\mathbb K^{(\alpha,\frac{1}{2})}(x,y) 
	= \frac{-1}{x-y} 
		\sum_{j=0}^2 \sum_{i=0}^{2-j} (-1)^j a_{i+j}
			\left( \vartheta_x^j f(x) \right)
			\left( \vartheta_y^i g(y) \right) 
			 \end{equation}
with $a_0 = \alpha(\alpha+ \frac{1}{2})$, $a_1 = - 2\alpha - \frac{1}{2}$, $a_2 = 1$, and
\[ f(x) = \MeijerG{1}{0}{0}{3}{-}{0,-\alpha, -\alpha-\frac{1}{2}}{x},
	\qquad 
	g(y) = \MeijerG{2}{0}{0}{3}{-}{\alpha,\alpha+\frac{1}{2},0}{y}. \]
Thus
\begin{equation} \label{Kformula2}
	\mathbb K^{(\alpha,\frac{1}{2})}(x,y) 
	= \frac{-1}{x-y}
		\begin{pmatrix} \vartheta_y^2 g(y) & \vartheta_y g(y) & g(y) 
		\end{pmatrix}
		\begin{pmatrix} 1 & 0 & 0 \\
		-2\alpha-\frac{1}{2} & 1 & 0 \\
		\alpha(\alpha+\frac{1}{2}) & - 2\alpha-\frac{1}{2} & 1
		\end{pmatrix}
		\begin{pmatrix} f(x) \\ - \vartheta_x f(x) \\
		\vartheta_x^2 f(x) \end{pmatrix}
		\end{equation}.

By \eqref{MG1003}, \eqref{phi4phi0}, and \eqref{phi4}
we have $f(x) = - \frac{1}{4\pi^2} (\phi_1(x) + \phi_2(x))$.
Also, by  \cite[page 151]{Lu} we have
\begin{align}
2\pi i 
\MeijerG{2}{0}{0}{3}{-}{\alpha,\alpha+\frac{1}{2},0}{z}  = 
\MeijerG{3}{0}{0}{3}{-}{0, \alpha,\alpha+\frac{1}{2}}{z e^{-\pi i}} -
\MeijerG{3}{0}{0}{3}{-}{0, \alpha,\alpha+\frac{1}{2}}{z e^{\pi i}}.
\end{align}
so that $g(y) = \frac{1}{2\pi i} (\psi_1(y) - \psi_2(y))$
by \eqref{psi1} and \eqref{psi2}.
Thus in view of \eqref{defPhialpha} and \eqref{defPsialpha},
\[  \begin{pmatrix} f(x) \\ \vartheta_x f(x) \\
		\vartheta_x^2 f(x) \end{pmatrix}
		= - \frac{1}{4\pi^2} \Phi_{\alpha,+}(x) \begin{pmatrix} 1 \\ 1 \\ 0 \end{pmatrix} \quad \text{and } \quad
 	\begin{pmatrix} \vartheta_y^2 g(y) \\ \vartheta_y g(y)\\ g(y) \end{pmatrix}
		= \frac{1}{2\pi i} \Psi_{\alpha_+}(y) 
		\begin{pmatrix} 1 \\ -1 \\ 0 \end{pmatrix} \]
with $x, y > 0$.
Thus \eqref{Kformula2} gives us
\begin{equation}
	\mathbb K^{(\alpha, \frac{1}{2})}(x,y)
	= \frac{1}{2\pi i(x-y)} 
		\begin{pmatrix} 1 & -1 & 0 \end{pmatrix}
		\frac{1}{4\pi^2} \Psi_{\alpha,+}^T(y)
		\begin{pmatrix} 1 & 0 & 0 \\
		-2\alpha-\frac{1}{2} & -1 & 0 \\
		\alpha(\alpha+\frac{1}{2}) & 2\alpha+\frac{1}{2} & 1
		\end{pmatrix}
		\Phi_{\alpha,+}(x) \begin{pmatrix} 1 \\ 1 \\ 0 \end{pmatrix}
\end{equation}
In view of \eqref{PhiInv} we arrive at \eqref{Kformula}.
\end{proof}

\section{Transformations of the Riemann-Hilbert problem}
\label{sec:first}

In sections \ref{sec:first}-\ref{sec:main} we apply 
the Deift-Zhou steepest descent method \cite{DeZh} to the RH problem
\ref{RHPforY} for $Y$. The method consists of a sequence of explict
transformations $Y \mapsto X \mapsto T \mapsto S \mapsto R$.

\subsection{First transformation $Y \mapsto X$}
We start with a preliminary transformation that simplifies the
jump condition in the RH problem on the positive real axis. 
It also introduces a jump on the negative real axis.

\begin{definition} \label{defX} 
If $n$  is even then  
\begin{align} \label{Xeven}
X(z) = \begin{pmatrix} 1 & 0 & 0\\ 0 & 1 & 0\\ 0 & 0 & i\end{pmatrix} Y(z)
\begin{pmatrix} 1 & 0 & 0\\ 0 & z^{\frac{1}{4}} & 0\\ 0 & 0 & z^{-\frac{1}{4}} \end{pmatrix}
\begin{pmatrix} 1 & 0 & 0\\ 0 & \frac{1}{\sqrt{2}} & \frac{1}{\sqrt{2}}\\ 0 & \frac{1}{\sqrt{2}} & -\frac{1}{\sqrt{2}}\end{pmatrix}
\begin{pmatrix} 1 & 0 & 0\\ 0 & 1 & 0\\ 0 & 0 & i\end{pmatrix},
\end{align}
with $z \in \mathbb C \setminus \mathbb R$, while if $n$ is odd then 
\begin{align} \label{Xodd}
X(z) = \begin{pmatrix} 1 & 0 & 0\\ 0 & 1 & 0\\ 0 & 0 & i\end{pmatrix}
	\begin{pmatrix} 1 & 0 & 0 \\ 0 & 0 & 1 \\ 0 & 1 & 0 \end{pmatrix}
	 Y(z) \begin{pmatrix} 1 & 0 & 0 \\ 0 & 0 & 1 \\ 0 & 1 & 0 \end{pmatrix}
\begin{pmatrix} 1 & 0 & 0\\ 0 & z^{-\frac{1}{4}} & 0\\ 0 & 0 & z^{\frac{1}{4}} \end{pmatrix}
\begin{pmatrix} 1 & 0 & 0\\ 0 & \frac{1}{\sqrt{2}} & \frac{1}{\sqrt{2}}\\ 0 & \frac{1}{\sqrt{2}} & -\frac{1}{\sqrt{2}}\end{pmatrix}
\begin{pmatrix} 1 & 0 & 0\\ 0 & 1 & 0\\ 0 & 0 & i\end{pmatrix}.
\end{align}
\end{definition}

The RH problem for $X$ follows from the definitions 
\eqref{Xeven}-\eqref{Xodd} and the RH problem \ref{RHPforY} for $Y$. 
While there is a different definition for $X$ in the cases 
$n$ even and $n$ odd, the RH problem for $X$ sees the two cases only 
in the jump \eqref{RHX2} on the negative real axis. We recall
\begin{align} \label{beta}
\beta=\alpha+\frac{1}{4}.
\end{align}

\begin{rhproblem} \label{RHPforX}  \
\begin{description}
\item[RH-X1] $X : \mathbb{C}\setminus \mathbb{R}\to \mathbb{C}^{3\times 3}$ is analytic.
\item[RH-X2] $X$ has boundary values for $x\in (-\infty,0) \cup 
(0,\infty)$, denoted by $X_{+}(x)$ and $X_{-}(x)$, and
\begin{align} \label{RHX1}
X_{+}(x) &= X_{-}(x) \begin{pmatrix} 1 & \sqrt{2} x^{\beta} e^{-nV(x)} & 0
\\ 0 & 1 & 0\\ 0 & 0 & 1\end{pmatrix}, && x > 0, \\
X_{+}(x) &= X_{-}(x) \begin{pmatrix} 1 & 0 & 0\\ 0 & 0 & (-1)^{n+1} 
\\ 0 & (-1)^n & 0
\end{pmatrix}, && x < 0, \label{RHX2}
\end{align}
\item[RH-X3] As $|z|\to\infty$
\begin{align}
X(z) = \left(\mathbb{I}+\mathcal{O}\left(\frac{1}{z}\right)\right) 
\begin{pmatrix} 1 & 0 & 0\\ 0 & z^{\frac{1}{4}} & 0\\ 0 & 0 & z^{-\frac{1}{4}}\end{pmatrix}
\begin{pmatrix} 1 & 0 & 0\\ 0 & \frac{1}{\sqrt{2}} & \frac{i}{\sqrt{2}}\\ 0 & \frac{i}{\sqrt{2}} & \frac{1}{\sqrt{2}}\end{pmatrix}
\begin{pmatrix} z^{n} & 0 & 0\\ 0 & z^{-\frac{n}{2}} & 0\\ 0 & 0 & z^{-\frac{n}{2}}\end{pmatrix}.
\end{align}
\item[RH-X4] As $z\to 0$
\begin{align}
\label{RHX4}
X(z) = \mathcal{O}\begin{pmatrix} 1 & z^{-\frac{1}{4}} h_{\alpha+\frac{1}{2}}(z) & z^{-\frac{1}{4}} h_{\alpha+\frac{1}{2}}(z)\\ 1 & z^{-\frac{1}{4}} h_{\alpha+\frac{1}{2}}(z) & z^{-\frac{1}{4}} h_{\alpha+\frac{1}{2}}(z)\\ 1 & z^{-\frac{1}{4}} h_{\alpha+\frac{1}{2}}(z) & z^{-\frac{1}{4}} h_{\alpha+\frac{1}{2}}(z)\end{pmatrix}.
\end{align}
\end{description}
\end{rhproblem}

While the RH problem \ref{RHPforX} follows in principle by
straightforward calculations, it may be good to point out
a few of the features leading to it.

First we may define $\tilde{Y} = Y$ in case $n$ is even
and 
\begin{equation} \label{defYtilde} \tilde{Y}(z) = 
  \begin{pmatrix} 1 & 0 & 0 \\ 0 & 0 & 1 \\ 0 & 1 & 0 \end{pmatrix}
  	Y(z) 
  \begin{pmatrix} 1 & 0 & 0 \\ 0 & 0 & 1 \\ 0 & 1 & 0 \end{pmatrix}
  \begin{pmatrix} 1 & 0 & 0 \\ 0 & z^{-1/2} & 0 \\ 0 & 0 & z^{1/2}
  \end{pmatrix} \end{equation}
in case $n$ is odd.
Then $\tilde{Y}$ satisfies of course the RH problem for $Y$ 
in case $n$ is even. For odd $n$, the RH problem is very close
to that of $Y$. Namely the jump on the positive real
axis remains the same with the jump matrix 
$\begin{pmatrix}
	1 & w(x) & x^{\frac{1}{2}} w(x) \\ 0 & 1 & 0 \\ 0 & 0 & 1
	\end{pmatrix}$ as before. For odd $n$, there is an additional 
jump on the negative real axis, but if we write the jump as
\[ \tilde{Y}_+ = \tilde{Y}_- \begin{pmatrix}
	1 & 0 & 0 \\ 0 & (-1)^n & 0 \\ 0 & 0 & (-1)^n \end{pmatrix}, \]
then it applies to both cases $n$ even and $n$ odd. 

In case $n$ is odd we obtain the asymptotic behavior of 
$\tilde{Y}$  from \eqref{RHY3} and \eqref{defYtilde}. It gives us
\begin{align} \nonumber
 \tilde{Y}(z) &=  \left( \mathbb{I} + \mathcal O(z^{-1}) \right)
	\begin{pmatrix} 1 & 0 & 0 \\ 0 & 0 & 1 \\ 0 & 1 & 0 \end{pmatrix}
	\begin{pmatrix} z^n & 0 & 0 \\  0 & z^{-\frac{n+1}{2}} & 0 \\
	0 & 0 & z^{-\frac{n-1}{2}} \end{pmatrix}
	\begin{pmatrix} 1 & 0 & 0 \\ 0 & 0 & 1 \\ 0 & 1 & 0 \end{pmatrix}
	\begin{pmatrix} 1 & 0 & 0 \\ 0 & z^{-\frac{1}{2}} & 0 \\ 0 & 0 & z^{\frac{1}{2}} \end{pmatrix} \\
	& = \left( \mathbb{I} + \mathcal O(z^{-1}) \right)
	\begin{pmatrix} z^n & 0 & 0 \\ 0 & z^{-\frac{n}{2}} & 0 \\ 0 & 0 & z^{-\frac{n}{2}}
	\end{pmatrix}
	\end{align}
and this is exactly the same as the asymptotic behavior \eqref{RHY3} for $n$ even.

The behavior near $0$ follows from \eqref{RHY4} and 
\eqref{defYtilde}. In case $n$ is odd it is
\begin{equation} \label{RHtY4} 
	\tilde{Y}(z) = \mathcal{O} 
\begin{pmatrix} 1 & z^{-1/2} h_{\alpha + \frac{1}{2}} & z^{1/2} h_{\alpha}(z) \\1 & z^{-1/2} h_{\alpha + \frac{1}{2}} & z^{1/2} h_{\alpha}(z) \\1 & z^{-1/2} h_{\alpha + \frac{1}{2}} & z^{1/2} h_{\alpha}(z) \end{pmatrix}. \end{equation}

From \eqref{Xeven}-\eqref{Xodd} and \eqref{defYtilde} we have
that in both case $n$ even and $n$ odd,
\begin{equation} \label{XinYt} 
	X(z) = \begin{pmatrix} 1 & 0 & 0\\ 0 & 1 & 0\\ 0 & 0 & i\end{pmatrix}\tilde{Y}(z)
\begin{pmatrix} 1 & 0 & 0\\ 0 & z^{\frac{1}{4}} & 0\\ 0 & 0 & z^{-\frac{1}{4}} \end{pmatrix}
\begin{pmatrix} 1 & 0 & 0\\ 0 & \frac{1}{\sqrt{2}} & \frac{1}{\sqrt{2}}\\ 0 & \frac{1}{\sqrt{2}} & -\frac{1}{\sqrt{2}}\end{pmatrix}
\begin{pmatrix} 1 & 0 & 0\\ 0 & 1 & 0\\ 0 & 0 & i\end{pmatrix}.
\end{equation}

To come from \eqref{XinYt} to the RH problem \ref{RHPforX}
requires some more calculations but these are more routine. 
We only note that for the asymptotic behavior near $0$,  
we use either \eqref{RHY4} or \eqref{RHtY4} in \eqref{XinYt}
and  \eqref{RHX4} follows in both cases since
\[ z^{\frac{1}{4}}h_{\alpha}(z) = 
	\mathcal{O}\left(z^{-\frac{1}{4}} h_{\alpha+\frac{1}{2}}(z)\right) \qquad \text{as } z \to 0, \]
as can be checked from the definition of $h_{\alpha}$
in \eqref{RHY4}. 
	
\subsection{Second transformation $X \mapsto T$}

The solution to the vector equilibrium problem \cite{Ku}, provides us with 
two measures $\mu^*$ and $\nu^*$ satisfying the variational conditions \eqref{ELvec1}
and \eqref{ELvec2}. Here $\mu^* = \mu_{V,\frac{1}{2}}^*$ is the 
$\theta = \frac{1}{2}$-equilibrium measure
in the external field $V$, which has a compact support 
\begin{align} \label{Delta1} 
	\supp(\mu^*) = \overline{\Delta_1}, \qquad \Delta_1 = (0,q) 
	\end{align}
by our assumption on the external field $V$ in Theorem \ref{mainThm}.
The measure $\nu^*$ has support
\begin{align} \label{Delta2} 
	\supp(\nu^*) = \overline{\Delta_2}, \qquad \Delta_2 = (-\infty,0). 
	\end{align}
Writing  
\[ U^{\mu^*}(z) = \int \log \frac{1}{|z-s|} d\mu^*(s),
	\qquad U^{\nu^*}(z) = \int \log \frac{1}{|z-s|} d\nu^*(s), \]
for the logarithmic potentials, we obtain from \eqref{ELvec1} and
\eqref{ELvec2},	
\begin{equation} 
\label{VarEqs}
\begin{aligned}
& 2 U^{\mu^*}-U^{\nu^*} + V  \begin{cases} = -\ell & \text{ on } [0,q], \\ 
	> -\ell  & \text{ on } (q,\infty), \end{cases}  \\ 
& 2 U^{\nu^*} = U^{\mu^*} \qquad \text{ on } (-\infty,0].
\end{aligned}
\end{equation}
The strict inequality on $(q,\infty)$ comes from the one cut
$\theta$-regular assumption in Theorem \ref{mainThm}, see Definition 
\ref{ocr}. Indeed the strict inequality that is assumed in
\eqref{ELvec2} for $x > q$, translates into the strict inequality 
in \eqref{VarEqs} for the minimizer of the vector equilibrium problem.

We define two $g$-functions by
\begin{align} \label{g1}
g_{1}(z) &= \int_{0}^{q} \log{(z-s)} \, d\mu^*(s), \qquad
	z \in \mathbb C \setminus (-\infty, q], \\ \label{g2}
g_{2}(z) & = \int_{-\infty}^{0} \log{(z-s)} \, d\nu^*(s), \qquad 
	z \in \mathbb C \setminus (-\infty, 0],
\end{align}
with principal branches of the logarithm. 
From the conditions \eqref{VarEqs} and definitions \eqref{g1}-\eqref{g2}
it follows that for $x>0$
\begin{equation} \label{VarGeq1}
\begin{aligned}
g_{1+}(x)-g_{1-}(x) &= 2\pi i \mu^*([x,\infty)), \\
g_{1-}(x)+g_{1+}(x)-g_{2}(x)-V(x) &
	\begin{cases} = \ell & \text{ on } [0,q], \\ 
	< \ell & \text{ on } (q,\infty), \end{cases} 
\end{aligned} 
\end{equation}
and for $x<0$
\begin{equation} \label{VarGeq2}
\begin{aligned}
g_{1+}(x)-g_{1-}(x) &= 2\pi i,\\
g_{2+}(x)-g_{2-}(x) &= 2\pi i\nu^*([x,0]),\\
-g_{1-}(x)+g_{2-}(x)+g_{2+}(x) &= \pi i.
\end{aligned}
\end{equation}
As before, the plus (minus) sign indicates that we use the limiting values
from the upper (lower) half-plane. 
		
\begin{proposition}
For the behavior as $z\to\infty$ we have
\begin{align} \label{g1expan}
g_{1}(z) &= \log{z}-\frac{m_{1}}{z}+\mathcal{O}\left(z^{-2}\right)\\
\label{g2expan}
g_{2}(z) &= \frac{1}{2}\log{z}+\frac{m_{\frac{1}{2}}}{\sqrt{z}}-\frac{m_{1}}{2 z}+\mathcal{O}\left(z^{-\frac{3}{2}}\right),
\end{align}
where
\begin{align} \label{m1}
m_{1} = \int_{0}^{q} s \, d\mu^*(s) \qquad \text{ and } \qquad 
m_{\frac{1}{2}} = \int_{0}^{q} \sqrt{s} d\mu^*(s).
\end{align}
\end{proposition}
\begin{proof}
The expansion \eqref{g1expan} is immediate from \eqref{g1}, since
$\mu^*$ is a probability measure with compact support. 
The first term of the expansion  \eqref{g2expan} for $g_2$ comes from
the fact that $\nu$ has total mass $1/2$. The full expansion 
can be obtained from the proofs of Propositions 3.1 and 3.2 in  \cite{Ku}. Namely, in the proof of Proposition 3.2 we find the relation
\begin{align}
\nu^* = \int_{0}^{q} \nu_t d\mu^*(t),
\end{align}
where $\nu_t$ for $t > 0$ is the minimizer of the energy functional
$I(\nu) + \int \log|s-t| d\nu(s)$ among measures on $(-\infty,0]$
with total mass $\frac{1}{2}$.
The proof of \cite[Proposition 3.1]{Ku}  yields
\begin{align*}
\int_{-\infty}^0 \log(z-s) d\nu_t(s) = 
	\log(z-t) - \log\left(z^{\frac{1}{2}} - \sqrt{t}\right) 
	= \log\left(z^{\frac{1}{2}} +\sqrt{t}\right).
\end{align*}
We therefore obtain from \eqref{g2} and the above
\begin{align} 
	g_{2}(z) & = 
	\int_0^q \left( \int_{-\infty}^0 \log(z-s) d\nu_t(s) \right) 
	d\mu^*(t)  \label{g2bis}
	 = \int_{0}^{q} \log\left(z^{\frac{1}{2}}+\sqrt{t}\right) d\mu^*(t). 
\end{align}
From the alternative expression \eqref{g2bis} for $g_2$ we 
find as $z \to \infty$
\begin{align*} 
	g_2(z) & = \log z^{\frac{1}{2}} + 
	\int_0^q \log\left(1+ \frac{\sqrt{t}}{z^{\frac{1}{2}}} \right) d\mu^*(t) 
	\\
	& = \frac{1}{2} \log z +
	\int_0^q \left(\frac{\sqrt{t}}{z^{\frac{1}{2}}}
		- \frac{t}{2 z} + \mathcal{O}\left(z^{-\frac{3}{2}}\right)
		\right) d\mu^*(t), 
	\end{align*}
which indeed yields the required expansion \eqref{g2expan} as $z\to\infty$. 
\end{proof}

The $g$ functions are used in the transformation $X \mapsto T$.
\begin{definition} \label{defT}  $T$ is defined by
\begin{align} \label{defT2}
T(z) = \begin{pmatrix} 1 & 0 & 0\\ 0 & 1 & i n m_{\frac{1}{2}}\\ 0 & 0 & 1\end{pmatrix}
L^{-1} X(z) \begin{pmatrix} e^{-n g_{1}(z)} & 0 & 0\\ 
	0 & e^{n(g_{1}(z)-g_{2}(z))} & 0\\ 0 & 0 & e^{n g_{2}(z)} 
		\end{pmatrix}  L,
\end{align}
where
\begin{align} \label{defL}
L = \begin{pmatrix} \sqrt{2} e^{\frac{2 n \ell}{3}} & 0 & 0\\ 
0 & e^{-\frac{n \ell}{3}} & 0\\ 0 & 0 & e^{-\frac{n \ell}{3}}\end{pmatrix}.
\end{align}
\end{definition}

In the RH problem for $T$ it will be convenient to use a function
that we call $\varphi$. To define it, we recall the assumption 
that $V$ is real analytic on $[0,\infty)$.
Thus $V$ has an analytic continuation to an open neighborhood $O_V$ 
of $[0,\infty)$ in the complex plane. Then 
\begin{align} \label{phi}
	\varphi(z) = -g_{1}(z)+\tfrac{1}{2} g_{2}(z) + \tfrac{1}{2} (V(z)+ \ell),
\end{align}
is defined and analytic for $z \in O_V \setminus (-\infty,q]$.
By \eqref{VarGeq1} and \eqref{phi}, 
\begin{align} \label{phipm}
	\varphi_{\pm}(x) & = \mp \pi i \mu^*([x, \infty)) = 
	\pm \pi i \mu^*([0,x])  \mp \pi i,
			\qquad x \in [0,q].
		\end{align}
From the variational inequality in  \eqref{VarGeq1} it follows that
\begin{equation} \label{phipos} 
	\varphi(x) >0 \qquad \text{ for real } x>q. 
	\end{equation}

Then $T$ satisfies the following RH problem.
\begin{rhproblem} \label{RHPforT} \
\begin{description}
\item[RH-T1] $T : \mathbb{C}\setminus \mathbb{R}\to \mathbb{C}^{3\times 3}$ is analytic.
\item[RH-T2] $T$ has boundary values for $x\in \Delta_{2} \cup \Delta_{1} \cup (q,\infty)$ and
\begin{align} \label{RHT2a}
T_{+}(x) &= T_{-}(x) \begin{pmatrix} e^{2 n \varphi_+(x)} & x^{\beta} & 0\\ 
0 & e^{2 n \varphi_-(x)} & 0\\ 0 & 0 & 1\end{pmatrix}, && x\in\Delta_{1} = (0,q), \\ \label{RHT2b}
T_{+}(x) &= T_{-}(x) \begin{pmatrix} 
1 & x^{\beta} e^{-2 n\varphi(x)} & 0\\ 0 & 1 & 0\\ 0 & 0 & 1
\end{pmatrix}, && x\in (q,\infty), \\ \label{RHT2c}
T_{+}(x) &= T_{-}(x) \begin{pmatrix} 1 & 0 & 0\\ 0 & 0 & -1\\ 0 & 1 & 0\end{pmatrix}, && x\in\Delta_{2} = (-\infty,0).   
\end{align}
\item[RH-T3] As $z\to\infty$
\begin{align}
\label{RHT3}
T(z) = \left(\mathbb{I}+\mathcal{O}\left(\frac{1}{z}\right)\right)
	\begin{pmatrix} 1 & 0 & 0\\ 0 & z^{\frac{1}{4}} & 0\\ 0 & 0 & z^{-\frac{1}{4}}\end{pmatrix}
\begin{pmatrix} 1 & 0 & 0\\ 0 & \frac{1}{\sqrt{2}} & \frac{i}{\sqrt{2}}\\ 0 & \frac{i}{\sqrt{2}} & \frac{1}{\sqrt{2}}\end{pmatrix}. 
\end{align}
\item[RH-T4] As $z\to 0$
\begin{align}
\label{RHT4}
T(z) = \mathcal{O}\begin{pmatrix} 1 & z^{-\frac{1}{4}} h_{\alpha+\frac{1}{2}}(z) & z^{-\frac{1}{4}} h_{\alpha+\frac{1}{2}}(z)\\ 1 & z^{-\frac{1}{4}} h_{\alpha+\frac{1}{2}}(z) & z^{-\frac{1}{4}} h_{\alpha+\frac{1}{2}}(z)\\ 1 & z^{-\frac{1}{4}} h_{\alpha+\frac{1}{2}}(z) & z^{-\frac{1}{4}} h_{\alpha+\frac{1}{2}}(z)\end{pmatrix}.
\end{align}
\end{description}
\end{rhproblem}

The behavior \eqref{RHT4} arises from \eqref{RHX4} 
since $g_1$ and $g_2$ remain bounded near $z$. The jump in
\eqref{RHT2a} follows from the variational conditions \eqref{VarGeq1}-\eqref{VarGeq2} and \eqref{phipm}.
For $x < 0$, we have from \eqref{defT2} and \eqref{RHX2}
\begin{multline}
	T_-^{-1}(x) T_+(x)  = 	
	L^{-1} \begin{pmatrix}
		e^{n g_{1,-}(x)} & 0 & 0 \\
		0 & e^{-n(g_{1,-}(x) - g_{2,-}(x))} & 0 \\
		0 & 0 & e^{-ng_{2,-}(x)} \end{pmatrix} \\
	\times
		\begin{pmatrix} 1 & 0 & 0 \\ 0 & 0 & -(-1)^n \\
		0 & (-1)^n & 0 \end{pmatrix}
		\begin{pmatrix}
		e^{-n g_{1,+}(x)} & 0 & 0 \\
		0 & e^{n(g_{1,+}(x) - g_{2,+}(x))} & 0 \\
		0 & 0 & e^{ng_{2,+}(x)} \end{pmatrix} L 
		\end{multline}
The matrices $L^{-1}$ and $L$ commute with all other
factors and thus cancel out because of
the special diagonal form \eqref{defL} of $L$.	
Then \eqref{RHT2c} follows since by \eqref{VarGeq2}
\[ e^{n(g_{1,+}(x) - g_{1,-}(x))} = 1 \quad \text{ and }
	\quad
	 e^{n(-g_{1,-}(x) + g_{2,+}(x) + g_{2,-}(x))} = (-1)^n \]
for $x < 0$.

Because of \eqref{phipos} the jump  matrix in \eqref{RHT2b}
tends to the identity matrix as $n \to \infty$.

The distinction between the cases $n$ even and
$n$ odd has disappeared in the RH problem for $T$.

\subsection{Third transformation $T \mapsto S$: opening of lenses}
\begin{figure}[t]
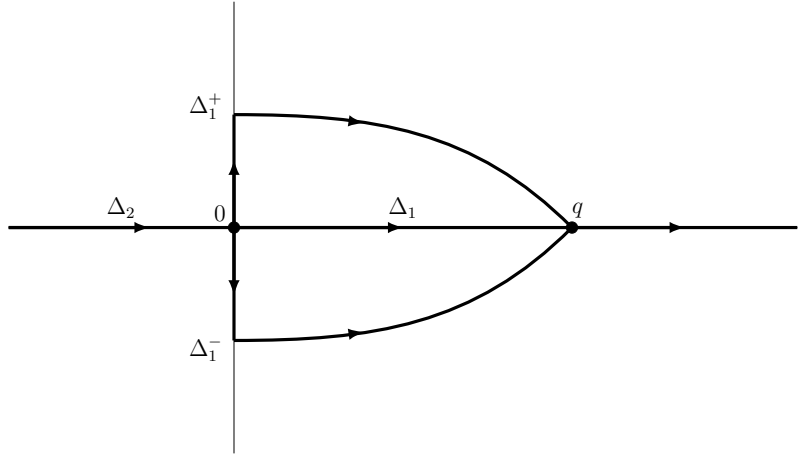

\begin{center}
\include{FiguurS}
\caption{Contour $\Sigma_{S} = \mathbb{R} \cup \Delta_{1}^{\pm}$ for the RH problem for $S$. The
	lense around $\Delta_1$ is contained in the domain $O_V$ where $V$ is analytic. \label{FigS}}
\end{center}
\end{figure}

In the next transformation $T \mapsto S$ we open a lense
around $\Delta_1 = (0,q)$ as in Figure \ref{FigS}. The lense is contained in
the neighborhood $O_V$ around $[0,\infty)$ where $V$ is analytic. 
We denote the upper and lower lips of the lense by $\Delta_1^+$
and $\Delta_1^-$, respectively, with orientation from $0$ to $q$.
The lips start from $0$ in a vertical direction. However, they do not necessarily
follow the imaginary axis exactly, contrary to how it is shown in the figure. Later we will construct
a conformal map $f$ in a neighborhood of $0$ with $f(0) = 0$. The  lips of
the lense near $0$ should be such that they are mapped by $f$ to the imaginary axis.

\begin{definition} \label{defS} $S$ is defined by 
\begin{align} \label{defS1}
S(z) & = T(z) \begin{pmatrix} 1 & 0 & 0\\ - z^{-\beta} e^{2 n\varphi(z)}  & 1 & 0\\ 0 & 0 & 1\end{pmatrix}, &&
\begin{array}{ll} \text{for $z$ in the upper} \\ \text{part of the lense}, \end{array} \\
\label{defS2}
S(z) & = T(z) \begin{pmatrix} 1 & 0 & 0\\ z^{-\beta} e^{2 n\varphi(z)}  & 1 & 0\\ 0 & 0 & 1\end{pmatrix}, &&
\begin{array}{ll} \text{for $z$ in the lower} \\ \text{part of the lense}, \end{array} \\
S(z) & = T(z), \label{defS3} && \text{ elsewhere}.
\end{align}
\end{definition}

The corresponding RH problem is
\begin{rhproblem} \label{RHPforS} \
\begin{description}
\item[RH-S1] $S : \mathbb{C}\setminus\Sigma_{S} \to \mathbb{C}^{3\times 3}$ is analytic (see Figure \ref{FigS} for the contour $\Sigma_S$). 
\item[RH-S2] On $\Sigma_S$ we have the jumps
\begin{align} \label{RHS2a}
S_{+}(x) &= S_{-}(x)\begin{pmatrix} 0 & x^{\beta} & 0\\ -x^{-\beta} & 0 & 0\\ 0 & 0 & 1\end{pmatrix}, && x\in\Delta_{1}, \\
\label{RHS2b}
S_{+}(z) &= S_{-}(z) \begin{pmatrix} 1 & 0 & 0\\ z^{-\beta} e^{2 n\varphi(z)} & 1 & 0\\ 0 & 0 & 1\end{pmatrix}, && z\in\Delta_{1}^{+} \cup \Delta_1^-,\\ \label{RHS2c}
S_{+}(x) &= S_{-}(x) \begin{pmatrix} 1 & x^{\beta} e^{-2 n\varphi(x)} & 0\\ 0 & 1 & 0\\ 0 & 0 & 1\end{pmatrix}, && x\in (q,\infty), \\
S_{+}(x) &= S_{-}(x) \begin{pmatrix} 1 & 0 & 0\\ 0 & 0 & -1\\ 0 & 1 & 0\end{pmatrix}, && x\in\Delta_{2}, \label{RHS2d}
\end{align}
\item[RH-S3] As $|z|\to\infty$
\begin{align}
\label{RHS3}
S(z) &= \left(\mathbb{I}+\mathcal{O}\left(\frac{1}{z}\right)\right)
\begin{pmatrix} 1 & 0 & 0\\ 0 & z^{\frac{1}{4}} & 0\\ 0 & 0 & z^{-\frac{1}{4}}\end{pmatrix}
\begin{pmatrix} 1 & 0 & 0\\ 0 & \frac{1}{\sqrt{2}} & \frac{i}{\sqrt{2}}\\ 0 & \frac{i}{\sqrt{2}} & \frac{1}{\sqrt{2}}\end{pmatrix}.
\end{align}
\item[RH-S4] As $z\to 0$
\begin{align}
\label{RHS4a}
S(z) &=  \mathcal{O}\begin{pmatrix} z^{-\alpha-\frac{1}{2}} h_{\alpha+\frac{1}{2}}(z) & z^{-\frac{1}{4}} h_{\alpha+\frac{1}{2}}(z) & z^{-\frac{1}{4}} h_{\alpha+\frac{1}{2}}(z)\\ z^{-\alpha-\frac{1}{2}} h_{\alpha+\frac{1}{2}}(z) & z^{-\frac{1}{4}} h_{\alpha+\frac{1}{2}}(z) & z^{-\frac{1}{4}} h_{\alpha+\frac{1}{2}}(z)\\ z^{-\alpha-\frac{1}{2}} h_{\alpha+\frac{1}{2}}(z) & z^{-\frac{1}{4}} h_{\alpha+\frac{1}{2}}(z) & z^{-\frac{1}{4}} h_{\alpha+\frac{1}{2}}(z)\end{pmatrix}
&& \text{ for } z \text{ inside the lense,}\\
\label{RHS4b}
S(z) &=  \mathcal{O}\begin{pmatrix} 1 & z^{-\frac{1}{4}} h_{\alpha+\frac{1}{2}}(z) & z^{-\frac{1}{4}} h_{\alpha+\frac{1}{2}}(z)\\ 1 & z^{-\frac{1}{4}} h_{\alpha+\frac{1}{2}}(z) & z^{-\frac{1}{4}} h_{\alpha+\frac{1}{2}}(z)\\ 1 & z^{-\frac{1}{4}} h_{\alpha+\frac{1}{2}}(z) & z^{-\frac{1}{4}} h_{\alpha+\frac{1}{2}}(z)\end{pmatrix}
&& \text{ for } z \text{ outside the lense.}
\end{align}
\end{description}
\end{rhproblem}

The asymptotic behavior \eqref{RHS4a} for $z\to 0$ inside the lens, comes from
\eqref{RHT4} combined with the  factor $z^{-\beta}$ that we have in 
the definitions \eqref{defS1}-\eqref{defS2}.
Note that $\varphi$ remains bounded as $z \to 0$. 

From the Cauchy-Riemann equations and \eqref{phipm} one can show 
that the real part of  $\varphi(z)$ is strictly negative for $z$ 
inside the lense with $\Re z  \in (0,q)$ provided we take the lense close enough to $[0,q]$. 
In section \ref{sec:local} we will do a closer analysis of $\varphi$ and
related functions near $0$. As a consequence of the assumed
one-cut $\theta$-regularity we have  \eqref{f1at0}
and together with \eqref{varphi1}, \eqref{varphi1f1f2} 
it implies that
\begin{align}
\varphi(z) = \begin{cases} -\pi i + 3\pi i c_{0,V} z^{\frac{1}{3}} + \mathcal{O}\left(z^\frac{2}{3}\right), & \text{for $z$ in the upper part of the lense},\\[5pt]
 \pi i - 3\pi i c_{0,V} z^{\frac{1}{3}} + \mathcal{O}\left(z^\frac{2}{3}\right), & \text{for $z$ in the lower part of the lense},\end{cases}
\end{align}
as $z\to 0$. Thus for $z$ close to the origin,
\begin{equation} \label{Revarphi} 
	\Re \varphi(z) \approx - 3 \pi c_{0,V} 
	\sin\left(\tfrac{1}{3}  |\arg (z)| \right)  |z|^{\frac{1}{3}} 
	\end{equation}
which is negative when $z$ is not on the positive real line.
In particular, $\Re \varphi(z) < 0$ for $z$ on the parts of $\Delta_1^{\pm}$ that are close
to $0$. 

By taking a smaller lense if necessary, we conclude that $\Re \varphi(z) < 0$
on the lips of the lense. Thus the jump matrices  on $\Delta_{1}^{\pm}$
and on  $(q,\infty)$ (this is due to \eqref{phipos}) in the RH problem for $S$ 
tend to the identity matrix  as $n\to\infty$.

\subsection{Global parametrix} 
We expect that $S$ is close to the solution of the 
following RH problem, at least away from the end points $0$ and $q$:

\begin{figure}[h]
\centering
\begin{picture}(300,50)(10,30)
\thicklines
\put(20,50){\line(1,0){180}}
\put(120,50){\circle*{3}}
\put(200,50){\circle*{3}}
\put(70,50){\vector(1,0){3}}
\put(160,50){\vector(1,0){3}}
\put(115,40){$0$}
\put(200,40){$q$}
\put(50,55){$\Delta_2$}
\put(165,55){$\Delta_1$}
\end{picture}
\caption{Contour $\Sigma_{N} = (-\infty,q]$ for the RH problem for $N$. \label{FigN}}
\end{figure}
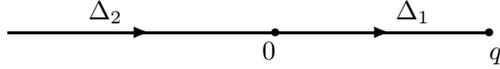

\begin{rhproblem} \label{RHPforN} \
\begin{description}
\item[RH-N1] $N : \mathbb{C}\setminus (-\infty,q] \to \mathbb{C}^{3\times 3}$ is analytic.
\item[RH-N2] On $\Sigma_N$ (see Figure \ref{FigN}) we have the jumps 
\begin{align}
N_{+}(x) &= N_{-}(x) \begin{pmatrix} 0 & x^\beta & 0\\ 
-x^{-\beta} & 0 & 0\\ 0 & 0 & 1\end{pmatrix}, && \text{for } x \in (0,q), \\
N_{+}(x) &= N_{-}(x) \begin{pmatrix} 1 & 0 & 0\\ 0 & 0 & -1\\ 0 & 1 & 0\end{pmatrix}, 
&& \text{for }x<0.
\end{align}
\item[RH-N3] As $z \to \infty$
\begin{align}
\label{RHN3}
N(z) = \left(\mathbb{I}+\mathcal{O}\left(\frac{1}{z}\right)\right)
\begin{pmatrix} 1 & 0 & 0\\ 0 & z^{\frac{1}{4}} & 0\\ 0 & 0 & z^{-\frac{1}{4}}\end{pmatrix}
\begin{pmatrix} 1 & 0 & 0\\ 0 & \frac{1}{\sqrt{2}} & \frac{i}{\sqrt{2}}\\ 0 & \frac{i}{\sqrt{2}} & \frac{1}{\sqrt{2}}\end{pmatrix}.
\end{align}
\end{description}
\end{rhproblem}

In \cite[section 4]{KuMFWi}  a global parametrix $N_{\alpha}$ is constructed that is
analytic in $\mathbb C \setminus (-\infty,q]$ with jumps
\begin{align*} 
N_{\alpha,+}(x) &= N_{\alpha,-}(x) \begin{pmatrix} 0 & x^\alpha & 0\\ 
-x^{-\alpha} & 0 & 0\\ 0 & 0 & 1\end{pmatrix}, \quad  x \in (0,q), \\
N_{\alpha,+}(x) &= N_{\alpha,-}(x) \begin{pmatrix} 1 & 0 & 0\\ 0 & 0 & -|x|^{-\alpha} \\ 0 & |x|^{\alpha} & 0\end{pmatrix}, \quad x<0,
\end{align*}
and asymptotic behavior
\begin{align*}
N_{\alpha}(z) = \left(\mathbb{I}+\mathcal{O}\left(\frac{1}{z}\right)\right)
\begin{pmatrix} 1 & 0 & 0\\ 0 & z^{\frac{1}{4}} & 0\\ 0 & 0 & z^{-\frac{1}{4}}\end{pmatrix}
\begin{pmatrix} 1 & 0 & 0\\ 0 & \frac{1}{\sqrt{2}} & \frac{i}{\sqrt{2}}\\ 0 & \frac{i}{\sqrt{2}} & \frac{1}{\sqrt{2}}\end{pmatrix} \begin{pmatrix} 1 & 0 & 0 \\ 0 & z^{\alpha/2} & 0 \\ 0 & 0 & z^{-\alpha/2}
\end{pmatrix} \quad \text{as } z \to \infty.
\end{align*}
It is straightforward to check that
\begin{equation} \label{defN}
N(z) = N_{2\beta}(z) \begin{pmatrix} 1 & 0 & 0 \\ 0 & z^{-\beta} & 0 \\ 0 & 0 & z^{\beta}
\end{pmatrix}
\end{equation}
satisfies the conditions in the RH problem \ref{RHPforN}. 

From \eqref{defN} and formulas (4.3) and (4.4) in \cite{KuMFWi} we also find the following
behaviors near $z = 0$ and $z =q$.
\begin{description}
\item[RH-N4] \
\begin{align}
\label{behavNin0}
N(z) \begin{pmatrix} z^{\frac{2\beta}{3}} & 0 & 0\\ 0 & z^{-\frac{\beta}{3}} & 0\\ 0 & 0 & z^{-\frac{\beta}{3}}\end{pmatrix} &= \mathcal{O}\begin{pmatrix} z^{-\frac{1}{3}} & z^{-\frac{1}{3}} & z^{-\frac{1}{3}}\\ z^{-\frac{1}{3}} & z^{-\frac{1}{3}} & z^{-\frac{1}{3}}\\ z^{-\frac{1}{3}} & z^{-\frac{1}{3}} & z^{-\frac{1}{3}}\end{pmatrix} &\text{ as }z\to 0, \\
N(z) &= \mathcal{O}\begin{pmatrix} (z-q)^{-\frac{1}{4}} & (z-q)^{-\frac{1}{4}} & 1\\ (z-q)^{-\frac{1}{4}} & (z-q)^{-\frac{1}{4}} & 1\\ (z-q)^{-\frac{1}{4}} & (z-q)^{-\frac{1}{4}} & 1\end{pmatrix} &\text{ as }z\to q.
\nonumber
\end{align}
\end{description}

\subsection{Local parametrix at the soft edge}

We consider an open disk $D(q,r_q)$ around $q$ with some small radius $r_q >0$, so that
$D(q,r_q)$ is contained in the region $O_V$ where $V$ is analytic. 
The local parametrix problem around $q$ is of the form:
\begin{rhproblem} \label{RHPforQ} \
\begin{description}
\item[RH-Q1] $Q$ is analytic on $D(q,r_q)\setminus (\Delta_{1}\cup \Delta_{1+} \cup \Delta_{1-}\cup [q,\infty))$. 
\item[RH-Q2] $Q$ has the same jumps as $S$ has on 
$(\Delta_{1}\cup \Delta_{1+} \cup \Delta_{1-}\cup [q,\infty)) \cap D(q,r_q)$. 
\item[RH-Q3] $Q$ has the same behavior as $z\to q$ as $S$ has on $D(q,r_q)$.
\item[RH-Q4] $Q$ matches with $N$ on the boundary:
\begin{align}
\label{RHQ4}
Q(z)N^{-1}(z) = \mathbb{I}+\mathcal{O}(n^{-1})\text{ as }n\to\infty,
\end{align}
uniformly for $z\in \partial D(q,r_q)$.
\end{description}
\end{rhproblem}

The construction of $Q$ can be done by a standard procedure using Airy functions \cite{De, DeKrML1, DeKrML2}. 
See \cite{Ku2} for a similar construction in a $3 \times 3$ matrix valued context.
We omit the details. 

\section{Local parametrix at the hard edge $0$}
\label{sec:local}

Our main task will be the construction of the local parametrix $P$ 
around the origin with the help of Meijer G-functions. Usually
such a local parametrix is constructed in a (possibly shrinking)
disk with a matching condition on the boundary of the disk. 
We deviate from this by constructing the local parametrix
on two concentric disks, which leads to a matching condition
on two circles. We introduce a new iterative technique to
improve the matching.  

\subsection{Statement of the RH problem for $P$}
\begin{figure}[t]
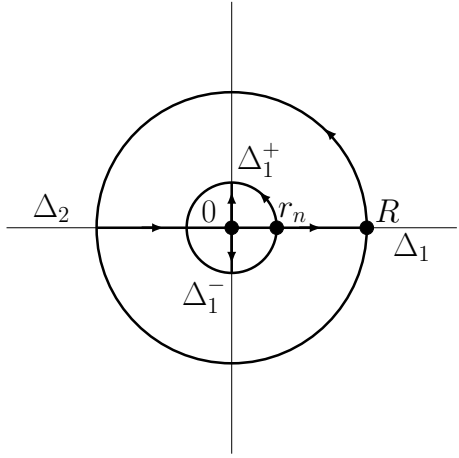

\centering
\include{FiguurP}
\caption{The contours for the local parametrix $P$. \label{FigP}}
\end{figure}

We consider two concentric disks with radii $r_n$ and $R$ 
which are centered at the origin, see Figure \ref{FigP}. For our purposes we will take 
\begin{align} \label{defrn}
	r_n & = n^{-\frac{3}{2}} 
	\end{align}
and $R$ is fixed, but small enough so that the boundary of the disk 
$D(0,R)$ is contained in the region $O_V$ where $V$ is analytic,
and such that the boundary of the disk meets the lips of the lens in
the (almost) vertical parts of $\Delta_{1}^{\pm}$, see Figure \ref{FigS}. 
We also want $R < q$. 
For large enough $n$ we of course
have $r_n < R$ and this is what we assume without further notice.
We use
\begin{equation} \label{annulus} 
	A(0; r_n, R) =
	\{ z \in \mathbb C \mid r_n < |z| < R \} \end{equation}
to denote the annular region bounded by the two circles
of radius $r_n$ and $R$. 
Our aim is to construct a local parametrix $P$ in the neighborhood
$D(0,R)$ of $0$ that also has a jump on the circle $|z| = r_n$.

\begin{rhproblem} \label{RHPforP} \
\begin{description}
\item[RH-P1] $P: D(0,R) \setminus (\partial D(0,r_n) \cup  
\Sigma_S) \to \mathbb C^{3 \times 3}$ is analytic. 
\item[RH-P2a] $P$ has the same jumps as $S$ has on 
$\Sigma_S \cap D(0,r_n)$, see \textbf{RH-S2} in the RH problem \ref{RHPforS}
for $S$.  That is, we have
\begin{align} \label{RHP2a}
P_{+}(x) &= P_{-}(x) \begin{pmatrix}  0 & x^\beta & 0\\ -x^{-\beta} & 0 & 0\\ 0 & 0 & 1\end{pmatrix}, 
&& x\in\Delta_{1}\cap D\left(0,r_n\right), \\
\label{RHP2b}
P_{+}(z) &= P_{-}(z) \begin{pmatrix} 1 & 0 & 0\\ z^{-\beta} e^{2 n\varphi(z)} & 1 & 0\\ 0 & 0 & 1\end{pmatrix}, 
&& 
z\in \left(\Delta_{1}^{+} \cup \Delta_2^+ \right) \cap D\left(0,r_n\right), \\
\label{RHP2c}
P_{+}(x) &= P_{-}(x) \begin{pmatrix} 1 & 0 & 0\\ 0 & 0 & -1\\ 0 & 1 & 0\end{pmatrix}, 
&& x\in\Delta_{2}\cap D\left(0,r_n\right),
\end{align}
where all contours are oriented as in Figure \ref{FigP}.
\item[RH-P2b]
$P$ has the same jumps as $N$ has on $\Sigma_N \cap A(0; r_n, R)$, 
see \textbf{RH-N2} in the RH problem \ref{RHPforN} for $N$, and no jump on 
$(\Delta_1^+ \cup  \Delta_1^-) \cap A(0; r_n, R)$. That is, we have
\begin{align} \label{RHPNa}
P_{+}(x) &= P_{-}(x) \begin{pmatrix}  0 & x^\beta & 0\\ -x^{-\beta} & 0 & 0\\ 0 & 0 & 1\end{pmatrix},
&& x\in\Delta_{1}\cap A(0; r_n, R), \\
\label{RHPNb}
P_{+}(x) &= P_{-}(x) \begin{pmatrix} 1 & 0 & 0\\ 0 & 0 & -1\\ 0 & 1 & 0\end{pmatrix}, 
&& x\in\Delta_{2}\cap A(0; r_n, R), \\
\label{RHPNc}
P_+(z) & = P_-(z), && z \in \Delta_1^{\pm} \cap A(0; r_n, R).
\end{align}
\item[RH-P3] As $z\to 0$, $P$ has the same behavior as $S$ has,
see \textbf{RH-S4} in the RH problem for $S$. 
That is, we have
\begin{align}
\label{RHP3a}
P(z) &=  \mathcal{O}\begin{pmatrix} z^{-\alpha-\frac{1}{2}} h_{\alpha+\frac{1}{2}}(z) & z^{-\frac{1}{4}} h_{\alpha+\frac{1}{2}}(z) & z^{-\frac{1}{4}} h_{\alpha+\frac{1}{2}}(z)\\ z^{-\alpha-\frac{1}{2}} h_{\alpha+\frac{1}{2}}(z) & z^{-\frac{1}{4}} h_{\alpha+\frac{1}{2}}(z) & z^{-\frac{1}{4}} h_{\alpha+\frac{1}{2}}(z)\\ z^{-\alpha-\frac{1}{2}} h_{\alpha+\frac{1}{2}}(z) & z^{-\frac{1}{4}} h_{\alpha+\frac{1}{2}}(z) & z^{-\frac{1}{4}} h_{\alpha+\frac{1}{2}}(z)\end{pmatrix}
&& \text{for } z \text{ to the right of $\Delta_1^{\pm} \cap D(0,R)$,}\\
\label{RHP3b}
P(z) &=  \mathcal{O}\begin{pmatrix} 1 & z^{-\frac{1}{4}} h_{\alpha+\frac{1}{2}}(z) & z^{-\frac{1}{4}} h_{\alpha+\frac{1}{2}}(z)\\ 1 & z^{-\frac{1}{4}} h_{\alpha+\frac{1}{2}}(z) & z^{-\frac{1}{4}} h_{\alpha+\frac{1}{2}}(z)\\ 1 & z^{-\frac{1}{4}} h_{\alpha+\frac{1}{2}}(z) & z^{-\frac{1}{4}} h_{\alpha+\frac{1}{2}}(z)\end{pmatrix}
&& \text{for } z \text{ to the left of $\Delta_1^{\pm} \cap D(0,R)$.}
\end{align}
\item[RH-P4] (Matching condition) We have
\begin{align}
\label{RHP4}
P_+(z)P_-^{-1}(z) &= \mathbb{I}+\mathcal{O}(n^{-1}) 
\quad \text{ as }n\to\infty, 
&& \text{uniformly for }z\in \partial D\left(0,r_n\right), \\
\label{RHPN4}
P(z)N^{-1}(z) &= \mathbb{I}+\mathcal{O}(n^{-1}) 
\quad \text{ as }n\to\infty, && \text{uniformly for }z\in \partial D\left(0,R\right).
\end{align}
\end{description}
\end{rhproblem}


\subsection{Ansatz for $P$}

We will look for $P$ in a form that contains the matrix valued function 
$\Phi_{\alpha}$ defined in \eqref{defPhialpha} and that 
is built out of Meijer G-functions. It has the piecewise constant jumps 
given in \eqref{Phialphajump}.
We will combine it with certain functions $\varphi_1$
and $\varphi_2$ to make sure that $P$ has the correct
jumps from \eqref{RHP2a}-\eqref{RHPNc}.
Of course, we also have to take care of the behavior 
\eqref{RHP3a}-\eqref{RHP3b} as $z \to 0$
and the matching condition \eqref{RHP4}-\eqref{RHPN4}.

Recall that $\varphi$ is given in \eqref{phi}. We are going to use 
the slightly modified version
\begin{align} \label{varphi1}
\varphi_{1}(z) =  \varphi(z) + 
	\begin{cases} \pi i, & \text{if }  z \in O_V, \, \Im(z)>0, \\
	-\pi i, & \text{if }  z \in O_V, \, \Im(z)<0. 
\end{cases}
\end{align}
We also introduce
\begin{align}
\label{varphi2}
\varphi_{2}(z) = -g_{2}(z)+\tfrac{1}{2} g_{1}(z) + 
\begin{cases}  -\frac{\pi i}{2}, & \text{if }  \Im(z)>0,\\ 
\frac{\pi i}{2}, & \text{if }  \Im(z)<0,
\end{cases}
\end{align}
where $g_1$ and $g_2$ are the functions from \eqref{g1}, \eqref{g2}.
The $\varphi_1$ and $\varphi_2$ functions appear in the following diagonal
matrix
\begin{align} \label{defD0}
D_0(z) & = \begin{pmatrix} 
	e^{\frac{2}{3}(2\varphi_{1}(z)+\varphi_{2}(z))} & 0 & 0\\ 
0 & e^{\frac{2}{3}(\varphi_{2}(z)-\varphi_{1}(z))} & 0\\ 
0 & 0 & e^{-\frac{2}{3}(\varphi_{1}(z)+2\varphi_{2}(z))}\end{pmatrix}.
\end{align}
 We will look for $P$ in the following form.
 
 \begin{definition}
Suppose $f$ is a conformal map and $E_{in}$ and $E_{out}$ are non-singular analytic matrix-valued functions on $D(0,r_n)$ and $A(0; r_n, R)$ respectively. 
Then we define
\begin{align} \label{Pansatz}
P(z) & = 
	E_{in}(z) \Phi_{\alpha}(n^3 f(z)) D_0^n(z)
	\begin{pmatrix} 1 & 0 & 0\\ 
	0 & z^{\beta} & 0\\ 
	0 & 0 & z^{\beta} \end{pmatrix}, && \text{ for } z \in D(0,r_n),
		\\  \label{Pansatz2}
	P(z) &= E_{out}(z) N(z), && \text{ for }  z \in A(0; r_n, R),
\end{align}
where $D_0$ is given by \eqref{defD0}.
\end{definition}

\begin{lemma} \label{Pansatzlemma}
Suppose that $f$ is a conformal map in $D(0,R)$
that maps the positive (negative) real line inside $D(0,R)$ to
part of the positive (negative) real line. 
Suppose also that the lense is opened such that
$\Delta_1^+ \cap D(0,R)$ is mapped to the positive imaginary axis
and $\Delta_1^- \cap D(0,R)$ is mapped to the negative imaginary axis.
Then for any analytic non-singular $3 \times 3$
matrix valued functions $E_{in}$ and $E_{out}$ 
the function $P(z)$ defined in \eqref{Pansatz} and \eqref{Pansatz2}
 satisfies the conditions \textbf{RH-P1}, \textbf{RH-P2a}, 
 \textbf{RH-P2b}, and \textbf{RH-P3} of the RH problem  for $P$.
\end{lemma}
\begin{proof}
The analyticity condition \textbf{RH-P1} is clearly satisfied.
The condition \textbf{RH-P2b} is also immediate since
the jumps in \eqref{RHPNa}-\eqref{RHPNc} are satisfied by $N$,
and they do not change if we multiply by an analytic prefactor
$E_{out}$ as in \eqref{Pansatz2}.
For \textbf{RH-P3} we note that by the definition of 
$h_{\alpha}$ in \eqref{RHY4} we check that 
$h_{-\alpha-\frac{1}{2}}(z) = z^{-\alpha-\frac{1}{2}} h_{\alpha+\frac{1}{2}}(z)$.
Then we use \eqref{Phiat0} and \eqref{Pansatz} where we note that
$\varphi_1(z)$ and $\varphi_2(z)$ remain bounded as $z \to 0$. Also $E_{out}(z)$
is bounded as $z \to 0$ and the conditions \eqref{RHP3a} and \eqref{RHP3b}  of \textbf{RH-P3} follow.
What remains is to prove \textbf{RH-P2a}.

We first check the jumps on the imaginary axis.
The functions $z \mapsto z^{\beta}$, $\varphi_1$ and $\varphi_2$ 
are analytic on the imaginary axis. Then by \eqref{Phialphajump}
and \eqref{Pansatz}, we have for $z \in \Delta_1{\pm} \cap D(0,r_n)$,
\begin{align} \nonumber
	P_{-}^{-1}(z) P_{+}(z) & = 
		\begin{pmatrix} 1 & 0 & 0\\ 
0 & z^{-\beta} & 0\\ 0 & 0 & z^{-\beta} 
\end{pmatrix} D_0^{-n}(z) 
\begin{pmatrix} 1 & 0 & 0 \\ 1 & 1 & 0 \\ 0 & 0 & 1 \end{pmatrix}
	D_0^n(z) \begin{pmatrix} 1 & 0 & 0\\ 
0 & z^{\beta}  & 0\\ 0 & 0 & z^{\beta} 
\end{pmatrix}  \\ \label{Pjump1calculation}
& = \begin{pmatrix} 1 & 0 & 0 \\
	z^{-\beta} e^{2n \varphi_1(z)} & 1 & 0 \\ 0 & 0 & 1 \end{pmatrix}.
\end{align}
Because of \eqref{varphi1} we may replace $\varphi_1$
by $\varphi$ in \eqref{Pjump1calculation}, since $n$  is an integer.
Thus \eqref{RHP2b} holds.

When checking the jumps on the real line, we have to take
into account the fact that $\varphi_1$ and $\varphi_2$
have jumps there. Also $z \mapsto z^{\beta}$ is defined
as the principal branch and so it has a jump on the negative
real axis (unless $\beta$ is an integer). 
For $x \in \Delta_1$, we find by \eqref{phipm} that
\begin{equation} \label{phi1Delta1}
	\varphi_{1,\pm}(x) = \pm \pi i \mu^*([0,x]) = 
	\pm \pi i(1- \mu^*([x,\infty))
\end{equation}
and therefore
\begin{equation} \label{phi1jumpDelta1} 
	\varphi_{1,-}(x) = - \varphi_{1,+}(x), \qquad x \in \Delta_1. 
	\end{equation}
From \eqref{varphi2} and \eqref{VarGeq1}, we find (since
$g_2$ is analytic across $\Delta_1$)
\[
	\varphi_{2,+}(x) - \varphi_{2,-}(x) = 
	\frac{1}{2} (g_{1,+}(x) - g_{1,-}(x) ) - \pi i
	= \pi i \mu^*([x,\infty)) - \pi i
	\]
which by \eqref{phi1Delta1} leads to
\begin{equation} \label{phi2jumpDelta1}
	 \varphi_{2,-}(x) = \varphi_{1,+}(x) + \varphi_{2,+}(x),
	\qquad x \in \Delta_1. \end{equation}
The two jump relations \eqref{phi1jumpDelta1} and \eqref{phi2jumpDelta1} 
are used in the calculation of the jump of $P$ across $\Delta_1$.

By  \eqref{Phialphajump} and \eqref{Pansatz}, we have 
for $x \in \Delta_1 \cap D(0,r_n)$,
\begin{multline} \label{Pjump2calculation}
	P_{-}^{-1}(x) P_{+}(x) = 
		\begin{pmatrix} 1 & 0 & 0\\ 
0 & x^{-\beta} & 0\\ 
0 & 0 & x^{-\beta} \end{pmatrix}  \left(D_{0}^{-n}\right)_-(z)  
\begin{pmatrix} 0 & 1 & 0 \\ -1 & 0 & 0 \\ 0 & 0 & 1 \end{pmatrix}
\left(D_{0}^n\right)_+(z) \begin{pmatrix}
	1  & 0 & 0\\  0 & x^{\beta}  & 0\\ 0 & 0 & x^{\beta} 
\end{pmatrix}.
\end{multline}
If we do the multiplication we find three non-zero entries with
an exponential factor containing certain linear combinations of
$\varphi_{1,\pm}(x)$ and $\varphi_{2,\pm}(x)$. Because of 
\eqref{phi1jumpDelta1} and \eqref{phi2jumpDelta1} all three linear 
combinations turn out to vanish identically
for $x \in \Delta_1$. The factors $x^{\beta}$ and $x^{-\beta}$
end up in the $12$- and $21$-entries, respectively, and 
the jump condition \eqref{RHP2a} follows.

Something similar happens on $\Delta_2$.
From \eqref{varphi2} and \eqref{g1}, \eqref{g2} we deduce
\begin{equation} \label{phi2delta2}
	\varphi_{2,\pm}(x) = \mp \pi i \nu^*([x,0]), \qquad x \in \Delta_2,
	\end{equation}
and therefore
\begin{equation} \label{phi2jumpDelta2}
	\varphi_{2,-}(x) = - \varphi_{2,+}(x), \qquad x \in \Delta_2,
\end{equation}
Also from \eqref{VarGeq2}, \eqref{phi}, \eqref{varphi1}, and the
fact that $V$ is analytic in $O_V$, we get for $x \in \Delta_2 \cap O_V$,
\[ \varphi_{1,+}(x) - \varphi_{1,-}(x) =
	-g_{1,+}(x) + g_{1,-}(x) - \tfrac{1}{2}(g_{2,+}(x) - g_{2,-}(x))
		+2\pi i = \pi i \nu^*([x,0]), \]
which by \eqref{phi2delta2} leads to
\begin{equation} \label{phi1jumpDelta2}
	\varphi_{1,-}(x) = \varphi_{1,+}(x) + \varphi_{2,+}(x), \qquad
	x \in \Delta_2 \cap O_V.
	\end{equation}
By  \eqref{Phialphajump} and \eqref{Pansatz}, we then have for $x \in \Delta_2 \cap D(0,r_n)$,
\begin{multline} \label{Pjump3calculation}
	P_{-}^{-1}(x) P_{+}(x) \\ = 
		\begin{pmatrix} 1 & 0 & 0\\ 
	0 & x_-^{-\beta} & 0\\ 
	0 & 0 & x_-^{-\beta} \end{pmatrix} \left( D_0^{-n}\right)_-(z) 
\begin{pmatrix} 1 & 0 & 0 \\ 0 & 0 & i e^{2\pi i \alpha} \\ 
0 & -i e^{2\pi i\alpha} & 0 \end{pmatrix}^{-1}
\left(D_0^n\right)_+(z) \begin{pmatrix}  
1  & 0 & 0\\ 0 & x_+^{\beta} & 0\\ 
0 & 0 & x_+^{\beta} 
\end{pmatrix}.
\end{multline}
Note that the inverse of the jump matrix in \eqref{Phialphajump}
appears in \eqref{Pjump3calculation}, since the orientation on the negative real 
axis is from left to right,  while in \eqref{Phialphajump} it is from right to left.

Multiplying the matrices in \eqref{Pjump3calculation} we get a jump matrix with
three non-zero entries.
Each entry has an exponential factor containing a linear combination of
$\varphi_{1,\pm}(x)$ and $\varphi_{2,\pm}(x)$. Each of these three linear
combinations vanishes because of the identities \eqref{phi1jumpDelta2}
and \eqref{phi2jumpDelta2}.
What we are left with from \eqref{Pjump3calculation} then is
\begin{equation} \label{Pjump3calculation2} 
	P_{-}^{-1}(x) P_{+}(x) = 
	\begin{pmatrix} 1 & 0 & 0 \\ 0 & 0 & i e^{-2\pi i \alpha} x_-^{-\beta} x_+^{\beta} \\
	0 & -i e^{2\pi i \alpha} x_-^{-\beta} x_+^{\beta} & 0 
	\end{pmatrix}, \qquad x \in \Delta_2 \cap O_V. \end{equation}
Finally, we also use that $x_-^{-\beta} = |x|^{-\beta} e^{\pi i \beta}$
and $x_+^{\beta}  = |x|^{\beta} e^{\pi i \beta}$, so that
$x_-^{-\beta} x_+^{\beta} = e^{2\pi i \beta} = i e^{2\pi i \alpha}$ for 
$x \in \Delta_2$,
where we recall $\beta = \alpha + \frac{1}{4}$, see \eqref{beta}. 
Then \eqref{Pjump3calculation2} reduces to \eqref{RHP2b}, as claimed.
\end{proof}

In view of Lemma \ref{Pansatzlemma} our remaining task
is to construct a specific conformal map $f$ and 
analytic prefactors $E_{in}$ and $E_{out}$ such that the 
matching condition \textbf{RH-P4} is satisfied as well. 

\subsection{The conformal map $f$}

Remember that $V$ can be analytically continued to an open neighborhood 
$O_V$ of $[0,\infty)$, and that $\overline{D(0,R)} \subset O_V$.

We let $O_V^0$  be such that $\overline{D(0,R)} \subset O_V^0 \subset O_V$
 and $O_V^0$ is slighly larger than $\overline{D(0,R)}$.
We make sure that $\Re z < q$ for $z \in O_V^0$.
On $O_V^0\setminus\mathbb{R}$ we define the auxiliary functions
\begin{align}
\label{deff1}
f_{1}(z) &= -z^{-\frac{1}{3}} \times \begin{cases} -\omega^2 \varphi_{1}(z)+\varphi_{2}(z), & \text{ for } z \in O_V^0,
 \ \Im(z)>0,\\
 -\omega \varphi_{1}(z)+\varphi_{2}(z), & \text{ for } z \in O_V^0, 
 \  \Im(z)<0, \end{cases} \\
\label{deff2}
f_{2}(z) &= -z^{-\frac{2}{3}} \times \begin{cases} -\omega \varphi_{1}(z)+\varphi_{2}(z), & \text{ for } z \in O_V^0, \  \Im (z)>0,\\
 -\omega^2 \varphi_{1}(z)+\varphi_{2}(z), & \text{ for } z \in O_V^0, \
   \Im(z)<0.\end{cases} 
\end{align}

\begin{proposition}\label{f1f2ana}
The functions $f_1$ and $f_2$ have analytic continuation to
$O_V^0$, and 
\begin{equation} \label{f1at0} 
	f_1(0) = 3 \sqrt{3} \pi c_{0,V} > 0. 
	\end{equation}
 \end{proposition}

\begin{proof}
It follows from \eqref{phi1jumpDelta1}, \eqref{phi2delta2}, \eqref{phi1jumpDelta2} and \eqref{phi2jumpDelta1} that $f_1$ and $f_2$ have no jumps on $O_V^0 \cap \mathbb R$. For example for $x\in \Delta_2 \cap O_V^0$ we have
\begin{align}
\nonumber
f_{1-}(x)  &= -e^{\frac{\pi i}{3}} |x|^{-\frac{1}{3}} (-\omega \varphi_{1-}(x) + \varphi_{2-}(x))\\
\nonumber
&= -e^{\frac{\pi i}{3}} |x|^{-\frac{1}{3}} (-\omega (\varphi_{1+}(x)+\varphi_{2+}(x)) + \varphi_{2-}(x))\\
\nonumber
&= -e^{\frac{\pi i}{3}} |x|^{-\frac{1}{3}} (-\omega \varphi_{1+}(x) - (1+\omega) \varphi_{2+}(x))\\
&=  -e^{-\frac{\pi i}{3}} |x|^{-\frac{1}{3}} (-\omega^2 \varphi_{1+}(x) + \varphi_{2+}(x)) = f_{1+}(x),
\end{align}
where we have used \eqref{phi1jumpDelta2} in the second line and \eqref{phi2delta2} in the third line. We conclude that $f_1$ and $f_2$ have a Laurent series around $z=0$.
In the upper half plane we have
\begin{align}
z^{1/3} f_{1}(z) + z^{2/3} f_{2}(z) &= -\varphi_{1}(z)-2\varphi_{2}(z)\\
\omega^2 z^{1/3} f_{1}(z) + \omega z^{2/3} f_{2}(z) &= \varphi_{2}(z)-\varphi_{1}(z).
\end{align}
It follows that (where we use
$1 + 2 \omega = i \sqrt{3}$ and $1+ 2 \omega^2 = - i \sqrt{3}$),
\begin{align} \label{varphi1f1f2}
  i \sqrt{3} z^{1/3} f_{1}(z) - i \sqrt{3} z^{2/3} f_{2}(z) &=  3\varphi_{1}(z), \qquad \Im z > 0.
\end{align}
Thus using the behavior \eqref{EPreg} of the density of the equilibrium measure $\mu^*$ we have
\begin{align}
\label{f10c0}
\sqrt{3}i x^{1/3} f_{1}(x) - \sqrt{3} i x^{2/3} f_{2}(x) = 3\varphi_{1+}(x)= 3\pi i\mu^{*}([0,x]) \sim 9\pi i c_{0,V} x^{1/3} \qquad
	\text{ as } x \downarrow 0.
\end{align}
This is only possible if $f_{1}$ and $f_{2}$ do not have negative powers in their respective Laurent series. 
We conclude that $f_{1}$ and $f_{2}$ can be analytically continued to $O_V^0$.

The value \eqref{f1at0} follows from \eqref{f10c0}. 
\end{proof}

Inverting the relations \eqref{deff1}, \eqref{deff2}, we see
that for $z \in O_V^0$, $\Im z > 0$, 
\begin{equation} 
\label{matchingn}
\begin{aligned}
2\varphi_{1}(z)+\varphi_{2}(z) &= \omega z^{1/3} f_{1}(z) + \omega^2 z^{2/3} f_{2}(z), \\
\varphi_{2}(z)-\varphi_{1}(z) &= \omega^2 z^{1/3} f_{1}(z) + \omega z^{2/3} f_{2}(z), \\
-\varphi_{1}(z)-2\varphi_{2}(z) & = z^{1/3} f_{1}(z) + z^{2/3} f_{2}(z),
\end{aligned}
\end{equation}
A similar relation, with $\omega$ and $\omega^2$ interchanged, holds for $\Im z < 0$. Observe that the combinations \eqref{matchingn}
appear in the diagonal matrix $D_0$, see \eqref{defD0}.

Our choice of conformal map is dictated by the desire that
$D_1(n^3 f(z))$ should cancel with $D_0^n(z)$ in \eqref{Pansatz}
as much as possible, where
\begin{equation} \label{defD1}
	D_1(z) = \begin{cases}
		\begin{pmatrix}
		e^{-3 \omega z^{1/3}} & 0 & 0 \\0 & e^{-3 \omega^2 z^{1/3}} & 0 \\
		0 & 0 & e^{-3z^{1/3}} \end{pmatrix},
		& \text{ for } \Im z > 0, \\
		\begin{pmatrix}
		e^{-3 \omega^2 z^{1/3}} & 0 & 0 \\0 & e^{-3 \omega z^{1/3}} & 0 \\
		0 & 0 & e^{-3z^{1/3}} \end{pmatrix},
		& \text{ for } \Im z < 0, 
		\end{cases} 
\end{equation}
is the exponential part of the asymptotic formula for $\Phi_{\alpha}$.
Ideally, we would have $D_1(n^3 f(z)) D_0^n(z) = 
\mathbb I$, but this is not
possible. 

\begin{definition} Let $f_1$ be given by \eqref{deff1}.
We define for $z \in O_V^0$,
\begin{align} \nonumber
f(z) & = \frac{8}{729} z f_1^3(z) \\
	& = \frac{8}{729} \times \begin{cases}
		\left(\omega^2 \varphi_1(z) - \varphi_2(z) \right)^3,
			& \Im z > 0, \\[5pt]
		\left(\omega \varphi_1(z) - \varphi_2(z) \right)^3,
			& \Im z < 0. \end{cases} 
\label{defConformalf}
			\end{align}
\end{definition}
The second equality in \eqref{defConformalf} follows from
\eqref{deff1}.

\begin{proposition}
\label{conformalMapReal}
$f$ is a conformal map near $0$ that maps positive and negative numbers to positive and negative numbers, respectively.
\end{proposition}

\begin{proof}
From \eqref{defConformalf} and \eqref{f1at0} we learn that
$f'(0) = \frac{8}{729} f_1(0)^3 = \frac{8\pi^3}{3 \sqrt{3}} c_{0,V}^3 > 0$.
Thus $f$ is conformal in some small enough neighbourhood of $0$. 
Now let $x$ be a real number in this open neighbourhood, then
\begin{align}
\overline{f(x)} &= \frac{8}{729} \left(\omega \overline{\varphi_{1+}(x)}-\overline{\varphi_{2+}(x)}\right)^3
= \frac{8}{729} \left(\omega \varphi_{1-}(x) - \varphi_{2-}(x)\right)^3 = f(x).
\end{align}
Here we have used that $\overline{\varphi_{1+}(x)}=\varphi_{1-}(x)$ and $\overline{\varphi_{2+}(x)}=\varphi_{2-}(x)$, which follows from the definitions \eqref{varphi1}, \eqref{varphi2} and the fact that $\overline{g_{1+}(x)}=g_{1-}(x)$ and $\overline{g_{2+}(x)}=g_{2-}(x)$, which is immediate from the definitions of the $g$-functions. It follows that real numbers are mapped to real numbers. 
Since $f(0) = 0$ and $f'(0) > 0$ it then follows that
positive numbers are mapped to positive numbers, and
 negative numbers are mapped to negative numbers. 
\end{proof}

\begin{remark}
From now on we adjust the construction so that the conditions
regarding the conformal map in Lemma \ref{Pansatzlemma} are satisfied.
That is, we let $f$ be as in \eqref{defConformalf}, and
if necessary we decrease $R$ so that $f$ is a conformal map
in $D(0,R)$. We open the lense so that $f$ maps
\[ f \left(\Delta_1^{\pm} \cap D(0,R) \right) \subset i \mathbb R^{\pm}. \]
\end{remark}

\begin{proposition} \label{propD2} 
The conformal map $f$ is such that $D_1(f(z)) D_0(z) = D_2(z)$,
$z \in O_V^0$, with
\begin{equation} \label{defD2}
	D_2(z) = \begin{cases}
		\begin{pmatrix}
		e^{\frac{2}{3} \omega^2 z^{2/3} f_2(z)} & 0 & 0 \\
		0 & e^{\frac{2}{3} \omega z^{2/3} f_2(z)} & 0 \\
		0 & 0 & e^{\frac{2}{3} z^{2/3} f_2(z)} \end{pmatrix},
		& \text{ for } \Im z > 0, \\
		\begin{pmatrix}
		e^{\frac{2}{3} \omega z^{2/3} f_2(z)} & 0 & 0 \\
		0 & e^{\frac{2}{3} \omega^2 z^{2/3} f_2(z)} & 0 \\
		0 & 0 & e^{\frac{2}{3} z^{2/3} f_2(z)} \end{pmatrix},
		& \text{ for } \Im z < 0, 
		\end{cases} 
\end{equation}
\end{proposition}
\begin{proof}
	We combine 	\eqref{matchingn} with the definitions
	\eqref{defD0}, \eqref{defD1}, and \eqref{defConformalf}.
\end{proof}
Observe that for $|z| = r_n = n^{-3/2}$ the diagonal entries of $D_2^n(z)$ 
remain bounded and bounded away from $0$ as $n \to \infty$, since $f_2$ is
analytic by Proposition \ref{f1f2ana}.

\subsection{First step towards the matching} \label{prefactorEn}
From the asymptotic behavior of $\Phi_{\alpha}$
in \eqref{asympPhi} and Proposition \ref{propD2} we obtain
\begin{align} \label{firstStepBehav}
	T_{\alpha} \Phi_{\alpha}(n^3 f(z)) D_0^n(z) 
	= \frac{2\pi}{\sqrt{3}} \left(n^3 f(z)\right)^{-\frac{2\beta}{3}}
	\left( \mathbb{I} + \frac{A_{\alpha}}{n^3 f(z)} + 
	O\left(\frac{1}{n^6 f(z)^2}\right)	\right)
	L_{\alpha} (n^3 f(z))  D_2^n(z)
	\end{align}
with a constant matrix $A_\alpha$.
In the first step towards the matching condition, we forget about
the terms $\frac{A_{\alpha}}{n^3 f(z)}$ and $O\left(\frac{1}{n^6 f(z)^2}\right)$ and in view of \eqref{Pansatz}-\eqref{Pansatz2}, 
we aim to find $E_{in}^{(1)}$ and $E_{out}^{(1)}$ such that
\begin{equation} \label{Eattempt1}
	E_{in}^{(1)}(z) T_{\alpha}^{-1}
	 \frac{2\pi}{\sqrt{3}} \left(n^3 f(z)\right)^{-\frac{2\beta}{3}}
	L_{\alpha} (n^3 f(z))  D_2^n(z)
		\begin{pmatrix} 1 & 0 & 0 \\ 0 & z^{\beta} & 0 \\
		0 & 0 & z^{\beta} \end{pmatrix} = E_{out}^{(1)}(z) N(z),
		\qquad |z| = r_n, 
		\end{equation}
To that end we introduce the following definitions.
\begin{definition}
We define
\begin{align} \label{defFn}
	F_n(z) & = L_{\alpha}(n^3 f(z)) D_2^n(z), \\ \label{defEn}
	E_n(z) & = 	N(z) \begin{pmatrix} z^{\frac{2\beta}{3}} & 0 & 0 \\
		0 & z^{-\frac{\beta}{3}} & 0 \\ 0 & 0 & z^{-\frac{\beta}{3}}
		\end{pmatrix} F_n^{-1}(z) T_{\alpha}, \\ \label{defEin1}
	E_{in}^{(1)}(z) & = \frac{\sqrt{3}}{2\pi} n^{2\beta} 
	\left( \frac{f(z)}{z} \right)^{\frac{2\beta}{3}} E_n(z), \\ \label{defEout1}
	E_{out}^{(1)}(z) & = \mathbb{I}.
\end{align}
\end{definition}
With these definitions \eqref{Eattempt1} holds. In \eqref{defEin1} we separated the scalar factors from the matrix valued
factors which are in $E_n$. We need that $E_{in}^{(1)}$ is analytic
in a neighborhood of $0$, which in view of \eqref{defEin1} means
that we have to show that $E_n$ is analytic. 
	
\begin{lemma} \label{lemEninv}
\begin{enumerate}
\item[(a)] $N(z) \begin{pmatrix} z^{\frac{2\beta}{3}} & 0 & 0 \\
0 & z^{-\frac{\beta}{3}} & 0 \\ 0 & 0 & z^{-\frac{\beta}{3}}
\end{pmatrix}$,
$L_{\alpha}(z)$ 
and $F_n(z)$ all satisfy the same jump conditions
near $z=0$. Namely, all three matrix valued functions 
are analytic in $O_V^0 \setminus \mathbb R$,  with jump matrix 
$\begin{pmatrix} 0 & 1 & 0 \\ - 1 & 0 & 0 \\ 0 & 0 & 1 \end{pmatrix}$
for $ x > 0$, and
$\begin{pmatrix} e^\frac{4\pi i\beta}{3} & 0 & 0\\ 0 & 0 & -e^{-\frac{2\pi i\beta}{3}}\\ 0 & e^{-\frac{2\pi i\beta}{3}} & 0\end{pmatrix}$
for $x < 0$.
\item[(b)] In addition, all three matrix valued functions 
are $\mathcal{O}(z^{-1/3})$ as $z \to 0$.
\item[(c)]
$E_n(z)$ is an invertible and analytic matrix valued function 
in a neighborhood of $z=0$.
\end{enumerate}
\end{lemma}
\begin{proof}
(a) Using \textbf{RH-N2} we obtain for $x\in(0,q)$,
\begin{multline}
 \begin{pmatrix} x_-^{\frac{2\beta}{3}} & 0 & 0 \\ 0 & x_-^{-\frac{\beta}{3}} & 0 \\ 0 & 0 & x_-^{-\frac{\beta}{3}} \end{pmatrix}^{-1}
 N_-^{-1}(x) N_+(x) 
 \begin{pmatrix} x_+^{\frac{2\beta}{3}} & 0 & 0 \\ 0 & x_+^{-\frac{\beta}{3}} & 0 \\ 0 & 0 & x_+^{-\frac{\beta}{3}} \end{pmatrix}\\
 = \begin{pmatrix} x^{-\frac{2\beta}{3}} & 0 & 0 \\ 0 & x^{\frac{\beta}{3}} & 0 \\ 0 & 0 & x^{\frac{\beta}{3}}\end{pmatrix}
 \begin{pmatrix} 0 & x^\beta & 0\\ -x^{-\beta} & 0 & 0\\ 0 & 0 & 1\end{pmatrix}
 \begin{pmatrix} x^{\frac{2\beta}{3}} & 0 & 0 \\ 0 & x^{-\frac{\beta}{3}} & 0 \\ 0 & 0 & x^{-\frac{\beta}{3}} \end{pmatrix}
 =  \begin{pmatrix} 0 & 1 & 0\\ -1 & 0 & 0\\ 0 & 0 & 1\end{pmatrix},
\end{multline}
and for $x<0$ we get
\begin{multline}
 \begin{pmatrix} x_-^{\frac{2\beta}{3}} & 0 & 0 \\ 0 & x_-^{-\frac{\beta}{3}} & 0 \\ 0 & 0 & x_-^{-\frac{\beta}{3}} \end{pmatrix}^{-1}
 N_-^{-1}(x) N_+(x) 
 \begin{pmatrix} x_+^{\frac{2\beta}{3}} & 0 & 0 \\ 0 & x_+^{-\frac{\beta}{3}} & 0 \\ 0 & 0 & x_+^{-\frac{\beta}{3}} \end{pmatrix}\\
 = \begin{pmatrix} |x|^{-\frac{2\beta}{3}} e^{\frac{2\pi i\beta}{3}} & 0 & 0 \\ 0 & |x|^{\frac{\beta}{3}} e^{-\frac{\pi i\beta}{3}} & 0 \\ 0 & 0 & |x|^{\frac{\beta}{3}} e^{\frac{\pi i\beta}{3}}\end{pmatrix}
 \begin{pmatrix} 1 & 0 & 0\\ 0 & 0 & -1\\ 0 & 1 & 0\end{pmatrix}
 \begin{pmatrix} |x|^{\frac{2\beta}{3}} e^{\frac{2\pi i\beta}{3}} & 0 & 0 \\ 0 & |x|^{-\frac{\beta}{3}} e^{-\frac{\pi i\beta}{3}} & 0 \\ 0 & 0 & |x|^{-\frac{\beta}{3}} e^{-\frac{\pi i\beta}{3}} \end{pmatrix}\\
 =  \begin{pmatrix} e^\frac{4\pi i\beta}{3} & 0 & 0\\ 0 & 0 & -e^{-\frac{2\pi i\beta}{3}}\\ 0 & e^{-\frac{2\pi i\beta}{3}} & 0\end{pmatrix}.
\end{multline}
Next we notice using \eqref{defLalpha} that for $x > 0$
\begin{multline} \label{Ljump1}
L_{\alpha,-}^{-1}(x) L_{\alpha,+}(x) = 
\begin{pmatrix} e^{-\frac{2\pi i\beta}{3}} & 0 & 0\\ 0 & e^{\frac{2\pi i\beta}{3}} & 0\\ 0 & 0 & 1\end{pmatrix}^{-1}
\begin{pmatrix} \omega & -\omega^2 & 1\\ 1 & -1 & 1\\ \omega^2 & -\omega & 1\end{pmatrix}^{-1}
\begin{pmatrix} \omega^2 & \omega & 1\\ 1 & 1 & 1\\ \omega & \omega^2 & 1\end{pmatrix}
\begin{pmatrix} e^{\frac{2\pi i\beta}{3}} & 0 & 0\\ 0 & e^{-\frac{2\pi i\beta}{3}} & 0\\ 0 & 0 & 1\end{pmatrix}\\
= \begin{pmatrix} e^{\frac{2\pi i\beta}{3}} & 0 & 0\\ 0 & e^{-\frac{2\pi i\beta}{3}} & 0\\ 0 & 0 & 1\end{pmatrix}
\begin{pmatrix} 0 & 1 & 0\\ -1 & 0 & 0\\ 0 & 0 & 1\end{pmatrix}
\begin{pmatrix} e^{\frac{2\pi i\beta}{3}} & 0 & 0\\ 0 & e^{-\frac{2\pi i\beta}{3}} & 0\\ 0 & 0 & 1\end{pmatrix}
= \begin{pmatrix} 0 & 1 & 0\\ -1 & 0 & 0\\ 0 & 0 & 1\end{pmatrix},
\end{multline}
and for $x<0$ we obtain 
\begin{multline} \label{Ljump2}
L_{\alpha,-}^{-1}(x) L_{\alpha,+}(x) = \\
\begin{pmatrix} e^{-\frac{2\pi i\beta}{3}} & 0 & 0\\ 0 & e^{\frac{2\pi i\beta}{3}} & 0\\ 0 & 0 & 1\end{pmatrix}^{-1}
\begin{pmatrix} \omega & -\omega^2 & 1\\ 1 & -1 & 1\\ \omega^2 & -\omega & 1\end{pmatrix}^{-1}
\begin{pmatrix} \omega^2 & 0 & 0\\ 0 & 1 & 0\\ 0 & 0 & \omega\end{pmatrix}
\begin{pmatrix} \omega^2 & \omega & 1\\ 1 & 1 & 1\\ \omega & \omega^2 & 1\end{pmatrix}
\begin{pmatrix} e^{\frac{2\pi i\beta}{3}} & 0 & 0\\ 0 & e^{-\frac{2\pi i\beta}{3}} & 0\\ 0 & 0 & 1\end{pmatrix}\\
= \begin{pmatrix} e^{\frac{2\pi i\beta}{3}} & 0 & 0\\ 0 & e^{-\frac{2\pi i\beta}{3}} & 0\\ 0 & 0 & 1\end{pmatrix}
\begin{pmatrix} 1 & 0 & 0\\ 0 & 0 & -1\\ 0 & 1 & 0\end{pmatrix}
\begin{pmatrix} e^{\frac{2\pi i\beta}{3}} & 0 & 0\\ 0 & e^{-\frac{2\pi i\beta}{3}} & 0\\ 0 & 0 & 1\end{pmatrix}
= \begin{pmatrix} e^\frac{4\pi i\beta}{3} & 0 & 0\\ 0 & 0 & -e^{-\frac{2\pi i\beta}{3}}\\ 0 & e^{-\frac{2\pi i\beta}{3}} & 0\end{pmatrix}.
\end{multline}

Since $f$ is a conformal map that maps the positive (negative)
real line near $0$ to the positive (negative) real line near $0$,
we find that $L_{\alpha}(n^3 f(z))$ has the same jump matrices
\eqref{Ljump1}-\eqref{Ljump2} as $L_{\alpha}(z)$. 
Then by the definition \eqref{defFn} of $F_n$, we obtain for $x \in 
	O_V^0 \cap \mathbb R$,
\begin{align} \label{Fnjump1}
	F_{n,-}^{-1}(x) F_{n,+}(x) & =
		\left(D_{2}^{-n}\right)_-(x) 
		\begin{pmatrix} 0 & 1 & 0\\ -1 & 0 & 0\\ 0 & 0 & 1\end{pmatrix}
		\left(D_{2}^n\right)_+(x),			\qquad x > 0, \\
	F_{n,-}^{-1}(x) F_{n,+}(x) & =  \label{Fnjump2}
		\left(D_{2}^{-n}\right)_-(x) 
		\begin{pmatrix} e^\frac{4\pi i\beta}{3} & 0 & 0
		\\ 0 & 0 & -e^{-\frac{2\pi i\beta}{3}}\\ 
		0 & e^{-\frac{2\pi i\beta}{3}} & 0\end{pmatrix}
		\left(D_{2}^n\right)_+(x),	\qquad x < 0. 
	\end{align}
The definition \eqref{defD2} of $D_2$  is such that
(where we recall $f_2$ is analytic in a neighborhood of $0$)
\begin{align*}
	D_{2,+}(x) & = \begin{pmatrix} 0 & 1 & 0 \\
		1 & 0 & 0 \\ 0 & 0 & 1 \end{pmatrix} D_{2,-}(x)
		\begin{pmatrix} 0 & 1 & 0 \\
		1 & 0 & 0 \\ 0 & 0 & 1 \end{pmatrix}, \qquad x > 0, \\
	D_{2,+}(x) & = \begin{pmatrix} 1 & 0 & 0 \\
		0 & 0 & 1 \\ 0 & 1 & 0 \end{pmatrix} D_{2,-}(x)
		\begin{pmatrix} 1 & 0 & 0 \\
		0 & 0 & 1 \\ 0 & 1 & 0 \end{pmatrix}, \qquad x < 0.
		\end{align*}
Combining this with \eqref{Fnjump1}-\eqref{Fnjump2} and using the fact that
$D_2$ is a diagonal matrix, we obtain that $F_n$ has the correct jumps as claimed in part (a).		

\medskip

(b) Part (b) is a direct consequence of  \textbf{RH-N4}, 
\eqref{defLalpha} and \eqref{defConformalf}.

\medskip
(c) It follows from (a) and the definition \eqref{defEn} 
that $E_n$ has the identity jump on $O_V^0 \cap (\mathbb R \setminus \{0\})$
and therefore $E_n$ is analytic in $O_V^0 \setminus \{0\}$. 
By (b) it is $\mathcal{O}(z^{-2/3})$ as $z\to 0$ and we conclude that 
the singularity at $z=0$ is removable, and thus 
$E_n$ is analytic in $O_V^0$. Furthermore, since 
all the factors in the definition on $E_n$ are invertible, $E_n$ itself is also invertible.  
\end{proof}

As in \eqref{Pansatz}-\eqref{Pansatz2} we now make the first attempt in defining $P$.
\begin{definition}
	We define
\begin{equation}
\label{defP1}
	P^{(1)}(z) = 
	\begin{cases} E_{in}^{(1)}(z) \Phi_{\alpha}(n^3 f(z))
		D_0^n(z) \begin{pmatrix} 1 & 0 & 0 \\ 0 & z^{\beta} & 0 \\
		0 & 0 & z^{\beta} \end{pmatrix},
		& \text{ for } z \in D(0,r_n) \\
		E_{out}^{(1)}(z) N(z) = N(z), & \text{ for } z \in A(0; r_n, R),
		\end{cases}
	\end{equation}
with $E_{in}^{(1)}(z)$ and $E_{out}^{(1)}(z)$ as in \eqref{defEin1}-\eqref{defEout1}.
\end{definition}
	
In what follows we will give $\partial D(0,r_n)$ a positive orientation. 

\begin{corollary} \label{cormatch1} For $|z| = r_n$,
\begin{equation} \label{match1} 
	P^{(1)}_+(z) P^{(1)}_-(z)^{-1}
	= \mathbb{I} + \frac{A_n^{(1)}(z)}{n^3 z}  + E_n(z) \mathcal{O}(n^{-3}) E_n^{-1}(z) 
\end{equation}
where
\begin{equation} \label{defAn1}
	 A_n^{(1)}(z) =  \frac{1}{f'(0)}
	 E_n(z) T_{\alpha}^{-1} A_{\alpha} T_{\alpha}
	E_n^{-1}(z). \end{equation}
\end{corollary}
\begin{proof}
Let $|z|=r_n$. Then we have $\frac{1}{n^6 f(z)^2} = \mathcal{O}(n^{-3})$ and
$ \frac{1}{n^3 f(z)} = \frac{1}{n^3 f'(0) z} + \mathcal{O}(n^{-3})$.
Plugging this into \eqref{firstStepBehav} and using \eqref{defFn} gives us
\begin{align} 
	T_{\alpha} \Phi_{\alpha}(n^3 f(z)) D_0^n(z) 
	&= \frac{2\pi}{\sqrt{3}} \left(n^3 f(z)\right)^{-\frac{2\beta}{3}}
	\left( \mathbb{I} + \frac{A_{\alpha}}{n^3 f'(0) z} + 
	O\left(n^{-3}\right)	\right)
	F_n(z), \label{firstStepBehav2}
	\end{align}
Now using \eqref{defEn} and \eqref{firstStepBehav2} we get
\begin{multline*}
T_{\alpha} \Phi_{\alpha}(n^3 f(z)) D_0^n(z) 
\begin{pmatrix} 1 & 0 & 0 \\ 0 & z^{\beta} & 0 \\ 0 & 0 & z^{\beta} \end{pmatrix} N^{-1}(z)\\
= \frac{2\pi}{\sqrt{3}} \left( \frac{n^3 f(z)}{z} \right)^{-\frac{2\beta}{3}}
	\left( \mathbb{I} + \frac{A_{\alpha}}{n^3 f'(0) z} + 
	O\left(n^{-3}\right)	\right) T_\alpha E_n^{-1}(z).
\end{multline*}
Finally, from this and \eqref{defP1}, \eqref{defEin1}, \eqref{defEout1}, and \eqref{defEn} we arrive at
\begin{align*}
P^{(1)}_+(z) P^{(1)}_-(z)^{-1}
&= E_{in}^{(1)}(z) T_\alpha^{-1} \frac{2\pi}{\sqrt{3}} \left(\frac{n^3 f(z)}{z}\right)^{-\frac{2\beta}{3}}
	\left( \mathbb{I} + \frac{A_{\alpha}}{n^3 f'(0) z} + 
	O\left(n^{-3}\right)	\right)  T_\alpha E_n^{-1}(z) \\
&= E_n(z) T_\alpha^{-1} \left( \mathbb{I} + \frac{A_{\alpha}}{n^3 f'(0) z} + 
	O\left(n^{-3}\right)	\right) T_\alpha E_n^{-1}(z) \\
&= \mathbb{I} + \frac{A_n^{(1)}(z)}{n^3 z} + E_n(z) O\left(n^{-3}\right) E_n^{-1}(z),
\end{align*}
where we have used \eqref{defAn1} in the last line. 
\end{proof}

Unfortunately, \eqref{match1} is not $\mathbb{I} + "small"$ since both $E_n(z)$
and $E_n^{-1}(z)$ are $\mathcal{O}(n)$ as $n \to \infty$ with $|z| = r_n$,
as we prove in the next subsection.
Then it follows from \eqref{defAn1} that $A_n^{(1)}(z) = \mathcal{O}(n^2)$
and thus $\ds \frac{A_n^{(1)}(z)}{n^3 z} = \mathcal{O}(n^{1/2})$ as $n \to \infty$
with $|z|= r_n = n^{-\frac{3}{2}}$.

\subsection{Estimates on $E_n$}

The following result will turn out to be key to obtain the matching.

\begin{lemma}
\label{En1-En2Lem}
We have $E_n(z) = \mathcal{O}(n)$ and $E_n^{-1}(z) = \mathcal{O}(n)$
as $n \to \infty$, uniformly for $z \in \overline{D(0,r_n)}$.
Furthermore, we have
\begin{align} \label{Enestimate} E_n^{-1}(z_1) E_n(z_2) = 
\mathbb{I} + \mathcal{O}\left(n^\frac{5}{2}(z_1-z_2)\right) 
\end{align}
as $n \to \infty$, uniformly for $z_1, z_2 \in \overline{D(0,r_n)}$.
\end{lemma}
\begin{proof}
We write \eqref{defEn} with $F_n$ given in \eqref{defFn} as
\begin{multline} \label{Enfactors}
	E_n(z) = M_{\alpha}(z)
	\begin{pmatrix} n^{\frac{1}{2}} & 0 & 0 \\
		0 & 1 & 0 \\ 0 & 0 & n^{-\frac{1}{2}} \end{pmatrix} 
		L_{\alpha}(n^{3/2} z) D_2^{-n}(z) 
L_{\alpha}^{-1}(n^{3/2} f(z))
	\begin{pmatrix} n^{\frac{1}{2}} & 0 & 0 \\ 0 & 1 & 0 \\ 0 & 0 & n^{-\frac{1}{2}}
	\end{pmatrix} T_{\alpha} \end{multline}
with
\[ M_{\alpha}(z) = 
	N(z) \begin{pmatrix} z^{\frac{2\beta}{3}} & 0 & 0 \\
		0 & z^{-\frac{\beta}{3}} & 0 \\ 0 & 0 & z^{-\frac{\beta}{3}}
		\end{pmatrix} L_{\alpha}^{-1}(z).
\]
To obtain \eqref{Enfactors} we also used the identities 
\begin{align*} 
	L_{\alpha}(n^{3/2} z) & = \begin{pmatrix} n^{-\frac{1}{2}} & 0  & 0 \\
0 & 1 & 0 \\ 0 & 0 & n^{\frac{1}{2}} \end{pmatrix} L_{\alpha}(z), \quad
	L_{\alpha}(n^{3/2} f(z))  = \begin{pmatrix} n^{\frac{1}{2}} & 0 & 0 \\
	0 & 1 & 0 \\ 0 & 0 & n^{-\frac{1}{2}} \end{pmatrix} L_{\alpha}(n f(z)) 
	\end{align*}
	that follow from the definition \eqref{defLalpha} of $L_{\alpha}$.

Then $M_{\alpha}(z)$ is analytic by Lemma \ref{lemEninv} and does not depend on $n$.
Therefore it is certainly uniformly bounded for $z \in \overline{D(0,r_n)}$. 

Each of the factors in 
$L_{\alpha}(n^{3/2} z) D_2^{-n}(z) L_{\alpha}^{-1}(n^{3/2} f(z))$
remains bounded for $|z| = r_n$, and since the product is
analytic (which follows from Lemma \ref{lemEninv} and \eqref{defFn}), 
the product also remains bounded for $D(0,r_n)$, 
uniformly in $n$, by the maximum modulus principle. 
So the only factors in \eqref{Enfactors} that grow are the two
diagonal matrices and they are $\mathcal{O}(n^{1/2})$, which in total
leads to $E_n(z) = \mathcal{O}(n)$. Similarly $E_n^{-1}(z) = \mathcal{O}(n)$.
 
 To prove \eqref{Enestimate}, let us denote
 \begin{align}
\mathcal L_n(z) := L_{\alpha}(n^{3/2} z) D_2^{-n}(z) L_{\alpha}^{-1}(n^{3/2} f(z)).
 \end{align}
 We already observed that $\mathcal L_n$ is an analytic function that is bounded by an 
 $n$-independent constant on $\overline{D(0,r_n)}$.  The uniform bound 
 also works on $\overline{D(0,2r_n)}$.  Cauchy's integral formula  then gives us
 \begin{align}
\mathcal L_n(z_2)- \mathcal L_n(z_1) = 
	\frac{z_1-z_2}{2\pi i} \oint_{|z|=2 r_n} \frac{\mathcal L_n(z)}{(z-z_1)(z-z_2)} dz = 
	\mathcal O\left(n^{\frac{3}{2}}(z_1-z_2)\right)
 \end{align}
 uniformly for $z_1,z_2 \in \overline{D(0,r_n)}$. Combining this with $\mathcal L_n^{-1}(z_1) = \mathcal{O}(1)$, we get
 \begin{align}
 \mathcal L_n^{-1}(z_1) \mathcal L_n(z_2) = \mathbb{I} + \mathcal L_n^{-1} (z_1) (\mathcal L_n(z_2)-\mathcal L_n(z_1)) = 
 \mathbb{I} + \mathcal{O}\left(n^\frac{3}{2}(z_1-z_2)\right).
 \end{align} 
 uniformly for $z_1,z_2\in\overline{D(0,r_n)}$.
 
We also have $M_{\alpha}^{-1}(z_1) M_{\alpha}(z_2) = \mathbb{I} + \mathcal{O}(z_1-z_2)$ if $z_1, z_2 \in \overline{D(0,r_n)}$. Applying these  estimates and using \eqref{Enfactors} 
yields for $z_1, z_2 \in \overline{D(0,r_n)}$,
\begin{align*} 
E_n^{-1}(z_1) E_n(z_2) &= T_\alpha^{-1} 
\begin{pmatrix} n^{-\frac{1}{2}} & 0 & 0 \\
		0 & 1 & 0 \\ 0 & 0 & n^{\frac{1}{2}} \end{pmatrix} 
		\mathcal L_n^{-1}(z_1) 
		\begin{pmatrix} n^{-\frac{1}{2}} & 0 & 0 \\
		0 & 1 & 0 \\ 0 & 0 & n^{\frac{1}{2}} \end{pmatrix} 
		\left(\mathbb{I}+\mathcal{O}(z_1-z_2)\right)\\
		& \qquad \qquad \times \begin{pmatrix} n^{\frac{1}{2}} & 0 & 0 \\
		0 & 1 & 0 \\ 0 & 0 & n^{- \frac{1}{2}} \end{pmatrix} 
		\mathcal L_n(z_2) 
		\begin{pmatrix} n^{\frac{1}{2}} & 0 & 0 \\
		0 & 1 & 0 \\ 0 & 0 & n^{-\frac{1}{2}} \end{pmatrix} T_\alpha\\
		&= T_\alpha^{-1} 
\begin{pmatrix} n^{-\frac{1}{2}} & 0 & 0 \\
		0 & 1 & 0 \\ 0 & 0 & n^{\frac{1}{2}} \end{pmatrix} 
		\mathcal L_n^{-1}(z_1) 
		\left(\mathbb{I}+\mathcal{O}(n(z_1-z_2))\right) \mathcal L_n(z_2) 
		\begin{pmatrix} n^{\frac{1}{2}} & 0 & 0 \\
		0 & 1 & 0 \\ 0 & 0 & n^{-\frac{1}{2}} \end{pmatrix} T_\alpha \\
		&= 		T_\alpha^{-1} \begin{pmatrix} n^{-\frac{1}{2}} & 0 & 0 \\
		0 & 1 & 0 \\ 0 & 0 & n^{\frac{1}{2}} \end{pmatrix} 
		\left(\mathbb{I} + \mathcal{O}\left(n^{\frac{3}{2}}(z_1-z_2)\right)\right)
		\begin{pmatrix} n^{\frac{1}{2}} & 0 & 0 \\
		0 & 1 & 0 \\ 0 & 0 & n^{-\frac{1}{2}} \end{pmatrix} T_\alpha\\
		&= \mathbb{I} + \mathcal{O}\left(n^\frac{5}{2} (z_1-z_2)\right).
\end{align*}
as claimed in the lemma.
\end{proof}

Note that \eqref{Enestimate} implies that for $z_1, z_2 \in \overline{D(0,r_n)}$ we have
\[ E_n^{-1}(z_1) E_n(z_2) = \mathcal{O}(n) \]
while for $z_1 = \frac{x}{(c_V n)^{3}}$, $z_2 = \frac{y}{(c_V n)^{3}}$ with fixed $x,y> 0$ we have
\[ E_n^{-1}(z_1) E_n(z_2) = \mathbb{I} + \mathcal{O}
\left(n^{-\frac{1}{2}}(x-y)\right).\]

We find from Lemma \ref{En1-En2Lem} the following  estimates for the terms in \eqref{match1}.
\begin{corollary}
Uniformly for $|z| = r_n$, we have
$\ds \frac{A_n^{(1)}(z)}{n^3 z} = \mathcal{O}(n^{\frac{1}{2}})$ and  
$E_n(z) \mathcal{O}(n^{-3}) E_n^{-1}(z) = \mathcal{O}(n^{-1})$
as $n \to \infty$.
\end{corollary}

\subsection{Estimates on $A_n^{(1)}(z)$}

After the first step towards the matching condition 
we obtained $E_{in}^{(1)}$ and $E_{out}^{(1)}$ and a corresponding $P^{(1)}$, as in \eqref{defEin1}, \eqref{defEout1} and $\eqref{defP1}$, such that \eqref{match1} holds for $|z| = r_n$.
To proceed we need estimates on $A_n^{(1)}(z)$.

\begin{proposition}
Let $k=1,2,3,\ldots$. 
\begin{itemize}
\item[\rm (a)] We have uniformly for $z\in\overline{D(0,r_n)}$ that
\begin{equation} \label{Anestimate0} 
	\left( A_n^{(1)}(z) \right)^k = \mathcal{O}(n^2). \end{equation}
\item[\rm (b)] We have uniformly for $z_1,z_2,\ldots,z_k \in \overline{D(0,r_n)}$ that
\begin{equation} \label{Anestimate1}
	 A_n^{(1)}(z_1) \cdot \cdots \cdot A_n^{(1)}(z_k) = \mathcal{O}(n^{k+1}). 
	 \end{equation}
When $\ell$ indices $j\in\{1,2,\ldots,k-1\}$ are such that $z_j=z_{j+1}$ then this is replaced by
\begin{equation} \label{Anestimate2}
	 A_n^{(1)}(z_1) \cdot \cdots \cdot A_n^{(1)}(z_k) = \mathcal{O}(n^{k+1-\ell}). 
	 \end{equation}
\item[\rm (c)] We have uniformly for $z_1,z_2\in\overline{D(0,r_n)}$ that
\begin{align} \label{AnEnestimate1}
A_n^{(1)}(z_1) E_n(z_2) & = \mathcal{O}(n^2) \\ \label{AnEnestimate2}
 E_n^{-1}(z_1) A_n^{(1)}(z_2) & = \mathcal{O}(n^2). 
\end{align}
\end{itemize}
\end{proposition}

\begin{proof}
To prove (a) we notice using definition \eqref{defAn1} that
\begin{align}
\left( A_n^{(1)}(z) \right)^k = \frac{1}{f'(0)^k} E_n(z) T_\alpha^{-1} A_\alpha^k T_\alpha E_n^{-1}(z).
\end{align}
The only $n$-dependent factors here are $E_n(z)$ and $E_n^{-1}(z)$ which are both $\mathcal{O}(n)$ for $z \in \overline{D(0,r_n)}$ 
by Lemma \ref{En1-En2Lem}, leading to an overall $\mathcal{O}(n^2)$ behavior.

For (b) we first remark that by Lemma \ref{En1-En2Lem} we have for $j=1,2,\ldots,k-1$ that
\begin{align}
E_n^{-1}(z_j) E_n(z_{j+1}) = \mathcal{O}(n)
\end{align}
uniformly on $\overline{D(0,r_n)}$. This implies by \eqref{defAn1} that
\begin{align}
\nonumber
A_n^{(1)}(z_1) \cdot \cdots \cdot A_n^{(1)}(z_k) &= 
\frac{1}{f'(0)^k} E_n(z_1) \left(\prod_{j=1}^{k-1} T_\alpha^{-1} A_\alpha  T_\alpha E_n^{-1}(z_j) E_n(z_{j+1})\right) T_\alpha^{-1} A_\alpha  T_\alpha E_n^{-1}(z_k)\\
&= E_n(z_1) \mathcal{O}(n^{k-1}) E_n^{-1}(z_k).
\end{align}
The factors $E_n(z_1)$ and $E_n^{-1}(z_k)$ are $\mathcal{O}(n)$ uniformly by Lemma \ref{En1-En2Lem} and this leads to the $\mathcal{O}(n^{k+1})$ behavior.

A combination of this argument and the proof of part (a) 
gives  the second statement of (b).

To prove (c) we notice that
\begin{align}
A_n^{(1)}(z_1) E_n(z_2) = \frac{1}{f'(0)} E_n(z_1) T_\alpha^{-1} A_\alpha T_\alpha \left(E_n^{-1}(z_1) E_n(z_2)\right)
\end{align}
and this is indeed $\mathcal{O}(n^2)$ by Lemma \ref{En1-En2Lem}. A similar argument gives us \eqref{AnEnestimate2}.  
\end{proof}

\subsection{Second step towards the matching}
To improve the matching we now define

\begin{definition}
\begin{align} E^{(2)}_{in}(z) & = \label{E2in}
\left( \mathbb{I} - \frac{A_n^{(1)}(z)-A_n^{(1)}(0)}{n^3 z} \right) E_{in}^{(1)}(z), && z \in D(0,r_n) \\ \label{E2out}
E^{(2)}_{out}(z) & = \left(\mathbb{I} -  \frac{A_n^{(1)}(0)}{n^3z} \right)^{-1}, && z \in A(0; r_n, R),
	 \end{align}
and
\begin{equation} \label{defP2} 
	P^{(2)}(z) = \begin{cases}	\ds
\left( \mathbb{I} - \frac{A_n^{(1)}(z)-A_n^{(1)}(0)}{n^3 z} \right) P^{(1)}(z), &
	z \in D(0,r_n), \\ \ds
\left( \mathbb{I} - \frac{A_n^{(1)}(0)}{n^3 z} \right)^{-1} P^{(1)}(z), &
	z \in A(0; r_n, R).
	\end{cases} 
	\end{equation}
	\end{definition}
Note that $E_{in}^{(2)}(z)$ is analytic at $z=0$.

For $n$ large enough, the inverse $\left(\mathbb{I} -  \frac{A_n^{(1)}(0)}{n^3z} \right)^{-1}$ in \eqref{E2out} does 
indeed exist. The reason for this is 
that $\left(\frac{A_n^{(1)}(0)}{n^3 z} \right)^2 = 
\mathcal{O}(n^{-1})$ for $|z| \geq r_n$ 
by \eqref{Anestimate0}, and this implies that  for sufficiently
large $n$, $\frac{A_n^{(1)}(0)}{n^3 z}$ does not have an eigenvalue equal to $1$. 

\begin{lemma}
On $\partial D(0,r_n)$ we have the jump
\begin{align} \label{match2} 
	P^{(2)}_+(z) \left(P^{(2)}_-(z)\right)^{-1} = \mathbb{I} +
		\frac{A_n^{(2)}(z)}{n^6 z} 
	- \frac{A_n^{(1)}(0)A_n^{(1)}(z) A_n^{(1)}(0)}{n^9z^3} + \mathcal{O}(n^{-1})
	\end{align}
with
\begin{equation} \label{An2}
	A_n^{(2)}(z) =  
	\frac{ A_n^{(1)}(0)\left(A_n^{(1)}(z) - A_n^{(1)}(0)\right)}{z}
	\end{equation}
The terms in \eqref{match2} satisfy, for  $|z| = r_n$,
\begin{align}  \label{Anestimate3}
	\frac{A_n^{(2)}(z)}{ n^6 z} 
	& = \mathcal{O}(1), \quad \text{ and }  \quad
		\frac{A_n^{(1)}(0)A_n^{(1)}(z) A_n^{(1)}(0)}{n^9z^3} 
	 = \mathcal{O}(n^{- \frac{1}{2}}).
\end{align} 
\end{lemma}
\begin{proof}
On $\overline{D(0,r_n)}$ we have by \eqref{match1} and  \eqref{defP2}, 
\begin{multline*} 
P^{(2)}_+(z) \left(P^{(2)}_-(z)\right)^{-1} \\
=
\left( \mathbb{I} - \frac{A_n^{(1)}(z)-A_n^{(1)}(0)}{n^3 z} \right)
\left( \mathbb{I} + \frac{A_n^{(1)}(z)}{n^3 z} + E_n(z) \mathcal{O}(n^{-3}) E_n^{-1}(z) \right)
	\left(\mathbb{I} - \frac{A_n^{(1)}(0)}{n^3 z} \right) 
	 \end{multline*}
We first expand the product of the second and third factors
on the right to obtain
\begin{multline*} P^{(2)}_+(z) \left(P^{(2)}_-(z)\right)^{-1} 
= 
\left( \mathbb{I} - \frac{A_n^{(1)}(z)-A_n^{(1)}(0)}{n^3 z} \right) \\
\times 
	\left(\mathbb{I}+\frac{A_n^{(1)}(z)-A_n^{(1)}(0)}{n^3 z}
	- \frac{A_n^{(1)}(z)A_n^{(1)}(0)}{n^6 z^2} + E_n(z) \mathcal{O}(n^{-3}) E_n^{-1}(z) + \mathcal{O}(n^{-\frac{3}{2}}) \right)
\end{multline*}
The $\mathcal{O}(n^{-\frac{3}{2}})$ term comes from 
\begin{align}
E_n(z) \mathcal{O}(n^{-3}) E_n^{-1}(z) \frac{A_n^{(1)}(0)}{n^3z} = \mathcal{O}(n^{-\frac{3}{2}}), \qquad |z| = r_n,
\end{align}
where we used Lemma \ref{En1-En2Lem} and \eqref{AnEnestimate2}.
 When we multiply this with the remaining factor on the left we get
\begin{multline} P^{(2)}_+(z) \left(P^{(2)}_-(z)\right)^{-1} 
= \\
\mathbb{I} - \frac{A_n^{(1)}(z) A_n^{(1)}(0)}{n^6 z^2} - \left(\frac{A_n^{(1)}(z)-A_n^{(1)}(0)}{n^3 z}\right)^2
+\frac{(A_n^{(1)}(z)-A_n^{(1)}(0))A_n^{(1)}(z) A_n^{(1)}(0)}{n^9 z^3} + \mathcal{O}\left(\frac{1}{n}\right).
\end{multline}
Here we again used \eqref{Enestimate} and \eqref{Anestimate0}
and \eqref{AnEnestimate1} to
obtain the $\mathcal{O}\left(\frac{1}{n}\right)$ behavior. 

By \eqref{Anestimate0} and \eqref{Anestimate2} respectively we get
$
\frac{A_n^{(1)}(z)^2}{n^6 z^2} = \mathcal{O}\left(\frac{1}{n}\right)
$
and
$\frac{A_n^{(1)}(z)^2 A_n^{(1)}(0)}{n^9  z^3} = \mathcal{O}\left(n^{-\frac{3}{2}}\right)$ 
and we obtain \eqref{match2}.
The estimates in \eqref{Anestimate3} follow directly from \eqref{Anestimate1}.
\end{proof}

\subsection{Third step towards the matching}

In the next step we define

\begin{definition}
\begin{align} \label{E3in}
	E_{in}^{(3)}(z) & = \left( \mathbb{I} - \frac{A_n^{(2)}(z) - A_n^{(2)}(0)}{n^6z} \right) E_{in}^{(2)}(z),
		&& z \in D(0,r_n) \\
	\label{E3out}
	E_{out}^{(3)}(z) & = 
	\left( \mathbb{I} - \frac{A_n^{(2)}(0)}{n^6 z} \right)^{-1} E_{out}^{(2)}(z), && z \in A(0; r_n, R),
\end{align}
and
\begin{equation}
	P^{(3)}(z) = \begin{cases} \label{defP3}
	\ds \left( \mathbb{I} - \frac{A_n^{(2)}(z) - A_n^{(2)}(0)}{n^6z} \right) P^{(2)}(z), & z \in D(0,r_n) \\
	\ds \left( \mathbb{I} - \frac{A_n^{(2)}(0)}{n^6 z} \right)^{-1}
	P^{(2)}(z), & z \in A(0; r_n, R),
\end{cases} \end{equation}
\end{definition}

By \eqref{Anestimate1} we have
\begin{align}
\left(\frac{A_n^{(2)}(0)}{n^6 z}\right)^2 = \frac{A_n^{(1)}(0)(A_n^{(1)}(z) - A_n^{(1)}(0))A_n^{(1)}(0)(A_n^{(1)}(z) - A_n^{(1)}(0))}{n^{12} z^4} = \mathcal{O}\left(\frac{1}{n}\right)
\end{align}
and we may use this to conclude that for $n$ large enough, 
$\ds \frac{A_n^{(2)}(0)}{n^6 z}$ does not have eigenvalue $1$ for $|z| \geq r_n$. Therefore the inverse in \eqref{E3out} and 
\eqref{defP3} is well-defined. 

\begin{lemma}
On $\partial D(0,r_n)$ we have the jump
\begin{equation} \label{match3} 
	P^{(3)}_+(z) \left( P^{(3)}_-(z) \right)^{-1} = 
	\mathbb{I} + \frac{A_n^{(3)}(z)}{n^9 z} + \mathcal{O}(n^{-1})
	\end{equation}
with 
\begin{equation} \label{An3}
 A_n^{(3)}(z) = - 
	\frac{\left( A_n^{(1)}(z) - A_n^{(1)}(0) \right) A_n^{(1)}(z)
	\left( A_n^{(1)}(z) - A_n^{(1)}(0) \right)}{z^2},
	\end{equation}
and for $|z| = r_n$,
\begin{equation} \label{Anestimate4}
 \frac{A_n^{(3)}(z)}{n^9 z} = \mathcal{O}\left(n^{-\frac{1}{2}} \right).
	\end{equation}
\end{lemma}

\begin{proof}
By \eqref{match2} and \eqref{defP3} we have
\begin{multline} \label{P3jump1} 
	P^{(3)}_+(z) \left( P^{(3)}_-(z) \right)^{-1} = 	
	\left( \mathbb{I} - \frac{A_n^{(2)}(z) - A_n^{(2)}(0)}{n^6z} \right) \\
	\times
	\left( 	\mathbb{I} + \frac{A_n^{(2)}(z)}{n^6 z}
	+ \frac{A_n^{(1)}(0)A_n^{(1)}(z) A_n^{(1)}(0)}{n^9z^3} 
	+ \mathcal{O}(n^{-1}) \right)
	\left(\mathbb{I} - \frac{A_n^{(2)}(0)}{n^6 z} \right).
\end{multline}
By counting the number of $A_n^{(1)}$ factors it follows from \eqref{Anestimate1} and \eqref{An2} that for $|z| = r_n$,
\begin{align*}
\frac{A_n^{(2)}(z)}{n^6 z} \frac{A_n^{(2)}(0)}{n^6 z} &= \mathcal{O}\left(n^{-1}\right)\\
\frac{A_n^{(1)}(0)A_n^{(1)}(z) A_n^{(1)}(0)}{n^9z^3} \frac{A_n^{(2)}(0)}{n^6 z} &= \mathcal{O}\left(n^{-\frac{3}{2}}\right)\\
\frac{A_n^{(2)}(0)}{n^6 z} &= \mathcal{O}(1).
\end{align*}
Hence the product of the middle and the right factors in \eqref{P3jump1} equals
\[\mathbb{I} + \frac{A_n^{(2)}(z)-A_n^{(2)}(0)}{n^6 z}
	+ \frac{A_n^{(1)}(0)A_n^{(1)}(z) A_n^{(1)}(0)}{n^9z^3} 
	+ \mathcal{O}(n^{-1}). \]
To also incorporate the factor on the left, we again count the number of $A_n^{(1)}$ factors and use \eqref{Anestimate1} and \eqref{An2} to get
\begin{align*}
\left(\frac{A_n^{(2)}(z) - A_n^{(2)}(0)}{n^6z} \right)^2  &= \mathcal{O}\left(n^{-1}\right)\\
\frac{A_n^{(2)}(z) - A_n^{(2)}(0)}{n^6z} \frac{A_n^{(1)}(0)A_n^{(1)}(z) A_n^{(1)}(0)}{n^9z^3}  &= \mathcal{O}\left(n^{-\frac{3}{2}}\right)
\end{align*}
We conclude that
\begin{equation} \label{match3a} 
	P^{(3)}_+(z) \left( P^{(3)}_-(z) \right)^{-1} = 
	\mathbb{I} -  
	 \frac{A_n^{(1)}(0)A_n^{(1)}(z) A_n^{(1)}(0)}{n^9z^3} 
	+ \mathcal{O}(n^{-1}). \end{equation}

Note that by \eqref{An3} and \eqref{Anestimate2}
\begin{equation} \label{Anestimate5} 
	 - \frac{A_n^{(1)}(0)A_n^{(1)}(z) A_n^{(1)}(0)}{n^9z^3} 
	= \frac{A_n^{(3)}(z)}{n^9 z^3}
		+ \mathcal{O}\left(n^{-\frac{3}{2}}\right) \end{equation}
for $|z|=r_n$. We use \eqref{Anestimate5} in \eqref{match3a} 
and the result \eqref{match3} follows.
The estimate \eqref{Anestimate4} follows from \eqref{Anestimate3} and \eqref{Anestimate5}.
\end{proof}
Note that $A_n^{(3)}(z)$ is analytic at $z=0$.

We now have by \eqref{match3} and \eqref{Anestimate4}
that  uniformly on $\partial D(0,r_n)$
\begin{equation} \nonumber 
	P^{(3)}_+(z) \left( P^{(3)}_-(z) \right)^{-1} = \mathbb{I} 
			+ \mathcal{O}(n^{-\frac{1}{2}}) 
			\end{equation}
as $n\to\infty$. This is a weaker form of the matching condition
\eqref{RHP4}, with $\mathcal{O}(n^{-1})$ replaced by 
$\mathcal O(n^{-\frac{1}{2}})$. The weaker form is not sufficient
to prove the scaling limit of the correlation kernel. Therefore
we make another step. 

\subsection{Final step towards the matching}
The final step takes a by now familiar form.
\begin{definition}
We define
\begin{align} \label{E4in}
	E_{in}^{(4)}(z) & = 
	\left( \mathbb{I} - \frac{A_n^{(3)}(z)-A_n^{(3)}(0)}{n^9 z} \right)
		E_{in}^{(3)}(z) \\ \label{E4out}
	E_{out}^{(4)}(z) & = 
	\left( \mathbb{I} - \frac{A_n^{(3)}(0)}{n^9 z} \right)^{-1}
		E_{out}^{(3)}(z)
	\end{align}
and
\begin{equation} \label{defP4} 
	P^{(4)}(z) = \begin{cases} \ds
	\left( \mathbb{I} - \frac{A_n^{(3)}(z) - A_n^{(3)}(0)}{n^9z} \right) P^{(3)}(z), & z \in D(0,r_n) \\ \ds
	\left( \mathbb{I} - \frac{A_n^{(3)}(0)}{n^9 z} \right)^{-1}
	P^{(3)}(z), & z \in A(0; r_n, R),
\end{cases} \end{equation}
\end{definition}
Again, we argue in similar fashion as before that the inverse is well-defined provided $n$ is large enough. 

\begin{lemma}
On $\partial D(0,r_n)$ we have the jump
\begin{equation} \label{match4} 
	P^{(4)}_+(z) \left( P^{(4)}_-(z) \right)^{-1} = \mathbb{I} + \mathcal{O}(n^{-1}) \end{equation}
\end{lemma}

\begin{proof}
By \eqref{match3} and \eqref{defP4} the jump is given by
\begin{align}
\nonumber
	P^{(4)}_+(z) \left( P^{(4)}_-(z) \right)^{-1} &= 
	\left( \mathbb{I} - \frac{A_n^{(3)}(z)-A_n^{(3)}(0)}{n^9 z} \right)
	\left(\mathbb{I} - \frac{A_n^{(1)}(0)A_n^{(1)}(z) A_n^{(1)}(0)}{n^9z^3} + \mathcal{O}(n^{-1})\right)
	\left( \mathbb{I} - \frac{A_n^{(3)}(0)}{n^9 z} \right)\\ \label{P4P4}
	&= \left( \mathbb{I} - \frac{A_n^{(3)}(z)-A_n^{(3)}(0)}{n^9 z} \right)
	\left(\mathbb{I} + \frac{A_n^{(3)}(z)}{n^9 z} + \mathcal{O}(n^{-1})\right)
	\left( \mathbb{I} - \frac{A_n^{(3)}(0)}{n^9 z} \right).
\end{align}
By counting the number of $A_n^{(1)}$ factors and applying \eqref{Anestimate1} we obtain
\begin{align*}
 \frac{A_n^{(3)}(z)}{n^9 z}\frac{A_n^{(3)}(0)}{n^9 z} &= \mathcal{O}\left(n^{-2}\right) \qquad \text{and} \qquad 
\frac{A_n^{(3)}(0)}{n^9 z} = \mathcal{O}\left(n^{-\frac{1}{2}}\right).
\end{align*}
Hence the product of the middle and right factor in \eqref{P4P4} equals
\[
\mathbb{I} + \frac{A_n^{(3)}(z)-A_n^{(3)}(0)}{n^9 z} + \mathcal{O}(n^{-1}).
\]
Again by counting the number of $A_n^{(1)}$ factors and applying \eqref{Anestimate1} we obtain
\begin{align*}
 \left(\frac{A_n^{(3)}(z)-A_n^{(3)}(0)}{n^9 z}\right)^2 &= \mathcal{O}\left(n^{-2}\right)
\end{align*}
and we arrive at \eqref{match4}.
\end{proof}

We are ready to fix the definitions of $E_{in}$ and $E_{out}$.

\begin{definition} \label{defEinEout}
We set $E_{in} = E_{in}^{(4)}$ and $E_{out} = E_{out}^{(4)}$, where $E_{in}^{(4)}$ and $E_{out}^{(4)}$ are as in \eqref{E4in} and $\eqref{E4out}$.
\end{definition}

\begin{theorem} \label{matchingConditionThm2}
$P$, as defined in \eqref{Pansatz} and \eqref{Pansatz2} with
$E_{in} = E_{in}^{(4)}$ and $E_{out} = E_{out}^{(4)}$, 
satisfies the matching condition. That is, we have
\begin{align} \label{Pjump1}
P_+(z)P_-^{-1}(z) &= \mathbb{I}+\mathcal{O}(n^{-1}) 
\quad \text{ as }n\to\infty, 
&& \text{uniformly for }z\in \partial D\left(0,r_n\right), \\
P(z)N^{-1}(z) &= \mathbb{I}+\mathcal{O}(n^{-1}) 
\quad \text{ as }n\to\infty, && \text{uniformly for }z\in \partial D\left(0,R\right). \label{Pjump2}
\end{align}
\end{theorem}

\begin{proof}
The matching condition \eqref{Pjump1} on $\partial D(0,r_n)$ was just proved
in \eqref{match4}. 
To prove the matching condition \eqref{Pjump2} on $\partial D(0,R)$,
we notice that $P = P^{(4)}$ satisfies, by
\eqref{defP4}, \eqref{defP3}, \eqref{defP2}, \eqref{defP1}
\begin{equation} \label{P4product} 
	P(z) = \left(\mathbb{I} - \frac{A_n^{(3)}(0)}{n^9 z}\right)^{-1}
\left(\mathbb{I} - \frac{A_n^{(2)}(0)}{n^6 z}\right)^{-1}
\left(\mathbb{I} - \frac{A_n^{(1)}(0)}{n^3 z}\right)^{-1} N(z),
	\qquad z \in \partial D(0,R). \end{equation}

In \eqref{AnEnestimate1} we already noted that $A_n^{(1)}(0) = \mathcal O(n^2)$ and therefore
\[ 
\left(\mathbb{I} - \frac{A_n^{(1)}(0)}{n^3 z}\right)^{-1}
	= \mathbb{I} + \mathcal O \left(n^{-1}\right) \qquad \text{ as } n \to \infty, \]
	uniformly for $|z| = R$.
From the definition \eqref{An2} we see that $A_n^{(2)}$ contains terms with
products of two $A_n^{(1)}$ factors which by \eqref{AnEnestimate2}
gives the estimate $A_n^{(2)}(z) = O(n^{9/2})$ for $|z| = r_n$.
Since $A_n^{(2)}$ is analytic, we may then use the maximum
modulus principle to conclude that also $A_n^{(2)}(0) = O(n^{9/2})$
as $n \to \infty$.  As a result
\[ 
\left(\mathbb{I} - \frac{A_n^{(2)}(0)}{n^6 z}\right)^{-1}
	= \mathbb{I} + \mathcal O \left( n^{-\frac{3}{2}} \right) \]
	uniformly for $|z| = R$.
From \eqref{An3} we get that $A_n^{(3)}$ has terms with three  $A_n^{(1)}$ factors. 
We again use the estimate \eqref{AnEnestimate2} (now with $k=3$). Then 
reasoning as before  we obtain 
$A_n^{(3)}(0) = O(n^{7})$ as $n \to \infty$.
This leads to 
\[ 
\left(\mathbb{I} - \frac{A_n^{(3)}(0)}{n^9 z}\right)^{-1}
	= \mathbb{I} + \mathcal O \left( n^{-2} \right) \]
	uniformly for $|z| = R$.
Inserting all the $\mathcal{O}$ estimates in \eqref{P4product}is, we obtain \eqref{Pjump2}
and we are done. 
\end{proof}

\begin{remark}
As already noted at the end of the introduction, the matching condition can be an issue
in larger size Riemann-Hilbert problems, and other constructions have been proposed
before in the literature. In \cite{DeKuZh,DeKu,KuMFWi} the matching  condition was
established by means of a modification of the global parametrix. 
It also relied on the nilpotency of a term that appears in the jump matrix which is
analogous to our matrix $A_n^{(1)}$.

An iterative method was designed by Bertola and Bothner \cite{BeBo} in a situation
that is similar to ours. At the initial step of the construction of the local
parametrix also a ``matching" of the form $\mathbb{I} + \mathcal{O}(n^{1/2})$ is obtained.
Bertola and Bothner also modify the analytic prefactor step by step to improve
the matching. They use very detailed information on the $\mathcal{O}(n^{1/2})$ term
which is also nilpotent.

It differs from our method since we do not use the nilpotency, but instead we have
control on the behavior of the powers of $A_n^{(1)}$ as $n \to \infty$, see  \eqref{Anestimate0}.
We also need technical estimates on the analytic prefactor $E_n$ 
as in Lemma \ref{En1-En2Lem}. Having these estimates our transformations proceed
in a systematic fashion where the estimates come from a rather straightforward
counting procedure, namely  we count the number of $A_n^{(1)}$ factors as we go along.
We also introduced the novel trick to construct a matching condition on two circles.

We expect that the ideas behind our method may be more wider applicable. 
The first test case will be
the generalization of this work to the case $\theta = \frac{1}{r}$ with $r \geq 3$.
\end{remark}


\section{Final transformation and proof of Theorem \ref{mainThm}} \label{sec:main}

\subsection{Final transformation $S \mapsto R$}

Having the global parametrix $N$, and the local parametrices
$P$ (in the neighborhood $D(0,R)$ of $0$) and $Q$ (in the 
neighborhood $D(q,r_q)$ of $q$) we define
the  final transformation as follows. 

\begin{definition} We define $R$  as
\begin{align} \label{Rdef}
R(z) = \begin{cases}  
S(z) N^{-1}(z), & z \in \mathbb{C}\setminus 
	\left(\Sigma_S \cup \overline{D(0,R)}
	\cup \overline{D(q,r_q)} \right),\\
S(z) P^{-1}(z), & z \in D(0,R) \setminus 
	\left(\Sigma_S \cup \partial D(0,r_n) \right), \\
S(z) Q^{-1}(z), & z\in D(q,r_q) \setminus \Sigma_S.
\end{cases}
\end{align}
\end{definition}

\begin{figure}[t]
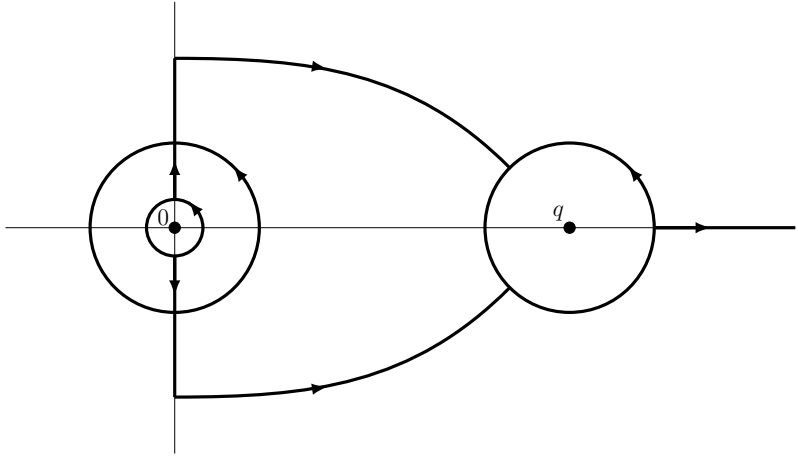

\centering
\include{FiguurR}
\caption{The contour $\Sigma_R$ for $R$. \label{FigR}}
\end{figure}

Since the jumps of $S$ and $N$ agree on $(-\infty,0) \cup (0,q)$
we find that $R$ is analytic across $(-\infty,-R)$ and
$(R,q-r_q)$. Also the jumps of $S$ and $Q$ agree inside
$D(q, r_q)$ and so $R$ is analytic in $D(q,r_q)$.
Finally, the jumps of $S$ and $P$ agree inside
$D(0,r_n)$ and on $(-R,0)$ and $(0,R)$. Therefore $R$ has analytic
continuation to $\mathbb C \setminus (\Sigma_R \cup \{0\})$
where
\begin{equation} \label{SigmaR} 
\Sigma_R = \partial D(0,r_n) \cup \partial D(0,R)
	\cup \partial D(q,r_q) \cup 
	\left( \Delta_1^{\pm} \setminus 
	\left( D(0,r_n) \cup D(q,r_q) \right) \right)
	\cup [q+r_q, \infty)
	\end{equation}
see Figure \ref{FigR}. We show that the isolated singularity
at $0$ is removable.

\begin{lemma} The singularity of $R$ at $0$ is removable,
and thus $R$ has analytic continuation to $\mathbb C \setminus
\Sigma_R$.
\end{lemma}

\begin{proof} We analyse the behavior of $R(z)$ as $z \to 0$ in the
left half plane. We have from \eqref{PhiInv} and \eqref{Pansatz}
\begin{align} \label{Rexpand}
	R(z) = S(z) \begin{pmatrix} 1 & 0 & 0 \\ 0 & z^{-\beta} & 0 
	\\ 0 & 0 & z^{-\beta} \end{pmatrix}
		D_0^{-n}(z) \Psi_{\alpha}^T(n^3f(z)) 
		\frac{-1}{4\pi^2} \begin{pmatrix} 1 & 0 & 0 \\
		-2\alpha - \frac{1}{2} & -1 & 0 \\
		\alpha(\alpha+\frac{1}{2}) & 2 \alpha + \frac{1}{2} & 1 
		\end{pmatrix} E_{in}^{-1}(z) \end{align}
We have from \eqref{RHS4b} that, as $z\to 0$ 
with $\Re z < 0$,
\begin{align*}
S(z) \begin{pmatrix} 1 & 0 & 0\\ 0 & z^{-\beta} & 0\\ 0 & 0 & z^{-\beta}\end{pmatrix} = 
\mathcal{O}\begin{pmatrix} 1 & h_{-\alpha-\frac{1}{2}}(z) & h_{-\alpha-\frac{1}{2}}(z)\\ 1 & h_{-\alpha-\frac{1}{2}}(z) & h_{-\alpha-\frac{1}{2}}(z)\\ 1 & h_{-\alpha-\frac{1}{2}}(z) & h_{-\alpha-\frac{1}{2}}(z)\end{pmatrix},
\end{align*}
since $z^{-\beta - \frac{1}{4}} h_{\alpha+\frac{1}{2}}(z)
	= h_{-\alpha-\frac{1}{2}}(z)$. Since $D_0^{-n}(z)$
	is a diagonal matrix which remains bounded as $z \to 0$,
	we also find
\begin{align} \label{Rexpand2}
S(z) \begin{pmatrix} 1 & 0 & 0\\ 0 & z^{-\beta} & 0\\ 0 & 0 & z^{-\beta}\end{pmatrix} D_0^{-n}(z) = 
\mathcal{O}\begin{pmatrix} 1 & h_{-\alpha-\frac{1}{2}}(z) & h_{-\alpha-\frac{1}{2}}(z)\\ 1 & h_{-\alpha-\frac{1}{2}}(z) & h_{-\alpha-\frac{1}{2}}(z)\\ 1 & h_{-\alpha-\frac{1}{2}}(z) & h_{-\alpha-\frac{1}{2}}(z)\end{pmatrix}, \end{align}
as $z \to 0$ with $\Re z < 0$. For $\Psi_{\alpha}$ we recall
the behavior \eqref{Psiat0} which for general $\alpha > -1$ is
\begin{align} \label{Rexpand3}
\Psi_\alpha^T(z) = 
\mathcal{O}\begin{pmatrix} h_{\alpha}(z) & h_{\alpha}(z) & h_{\alpha}(z)\\ z^\alpha & z^\alpha & z^\alpha \\ z^\alpha & z^\alpha & z^\alpha\end{pmatrix}
\end{align}
as $z \to 0$ with $\Re z < 0$. 

The remaining factors in \eqref{Rexpand} are either constant
or analytic at $z=0$ and therefore remain bounded as $z \to 0$. 
Then using  \eqref{Rexpand2} and  \eqref{Rexpand3} in
\eqref{Rexpand} we get
\[ R(z) =  \mathcal{O} \left( h_{\alpha}(z) \right) + 
	\mathcal {O} \left(z^{\alpha} h_{-\alpha-\frac{1}{2}}(z) \right)
	=   \begin{cases} z^\alpha, & -1 < \alpha < -\frac{1}{2}, \\ 
z^{-\frac{1}{2}}\log(z), & \alpha = -\frac{1}{2},\\ 
z^{-\frac{1}{2}}, & \alpha > -\frac{1}{2},\end{cases}
\]
as $z\to 0$ with $\Re z < 0$. 

In all cases we have $z R(z) \to 0$ as $z \to 0$ with 
$\Re z < 0$. This implies that $R$ cannot have a pole at $z=0$.
There cannot be an essential singularity either, and
therefore the singularity at $z=0$ is indeed removable.
\end{proof}

As a result of the steepest descent analysis,
and in particular the matching condition in Theorem \ref{matchingConditionThm2} we now have.

\begin{proposition}

\begin{enumerate}
\item[\rm (a)]
As $n \to \infty$ we have
\begin{align} \label{Rjump1}
R_{+}(z) & = R_{-}(z) 
	\left(\mathbb{I}+\mathcal{O}\left(\frac{1}{n}\right)\right),
	&& \text{ for }  z \in \partial D\left(0,R\right)
	\cup \partial D\left(0,r_n\right) \cup \partial D(q,r_q), \\
	\label{Rjump2}
R_+(z)	& = R_-(z) \left( \mathbb I + \mathcal O\left(e^{-c_1\sqrt{n}}\right) \right),
	&& \text{ for } z \in \Delta_1^{\pm} \cap A(0; r_n, R), \\
	\label{Rjump3}
R_+(z) & = R_{-}(z) \left(\mathbb{I}+\mathcal{O}\left(e^{-c_2 n} \right)\right), && \text{ on remaining parts of } \Sigma_R.
	 \end{align}
where $c_1, c_2$ are certain positive constants.
All $\mathcal O$ terms in \eqref{Rjump1}-\eqref{Rjump3}
are uniform for $z$ on the indicated contours.
\end{enumerate}
\item[\rm (b)] We have
\begin{align} \label{RisI}
R(z) = \mathbb{I}+\mathcal{O}\left(\frac{1}{n}\right) 	\qquad \text{as }n\to\infty
\end{align}
uniformly for $z\in \mathbb{C}\setminus \Sigma_R$.
\end{proposition}

\begin{proof}
For $z \in \partial D(0,r_n)$ we have by the definition \eqref{Rdef}
that
\[ R_-^{-1}(z) R_+(z) = \left(S(z) P_-(z)\right)^{-1} S(z) P_+(z)
	= P_-^{-1}(z) P_+(z) \]
since $S$ has no jump on $\partial D(0,r_n)$.
Then \eqref{Rjump1} for $z \in \partial D(0,r_n)$ 
follows from \eqref{Pjump1}.
Similarly, \eqref{Rjump1} for $z \in \partial D(0,R)$
follows from \eqref{Pjump2} and for $z \in \partial D(q,r_q)$
it follows from the matching condition \eqref{RHQ4} for
the local parametrix at $q$.

For $z \in \Sigma_R$ outside the disks $D(0,R)$ and $D(q,r_q)$ 
we have by \eqref{Rdef}
\[ R_-^{-1}(z) R_+(z) = N(z) S^{-1}_-(z)S_+(z) N^{-1}(z) \]
and $S^{-1}_-(z) S_+(z) = \mathbb{I} + \mathcal{O}(e^{-cn})$,
which follows from the jumps in the RH problem for $S$, since
$\Re \varphi(z) \leq -c <  0$ on $(\Delta_1^+ \cup \Delta_1^-) 
\setminus (D(0,R) \cup D(q, r_q))$, while
$\Re \varphi(x) \geq c > 0$ on $(q+r_q, \infty)$.
Then the estimate \eqref{Rjump3} also follows.

For \eqref{Rjump2} we have to be a little bit more careful.
For $z \in \Delta_1^{\pm} \cap A(0;r_n,R_n)$ we have $R(z) = S(z) P^{-1}(z)$ and $P_+(z) = P_-(z)$  by \eqref{RHPNc}. 
By the jump of $S$ in \eqref{RHS2b}, we get
\begin{align*}
	R_-^{-1}(z) R_+(z) & = P(z) S_-^{-1}(z) S_+(z) P^{-1}(z) \\
	& = \mathbb I + z^{-\beta} e^{2n \varphi(z)}
		P(z) \begin{pmatrix} 0 & 0 & 0 \\ 
	1 &  0 & 0 \\
	0 & 0 & 0 \end{pmatrix} P^{-1}(z).
	\end{align*}
On $\Delta_1^{\pm}$ we have by  \eqref{Revarphi}  that
$\Re \varphi(z) \leq - c |z|^{1/3}$   for some constant $c > 0$.
Since $|z| \geq n^{-\frac{3}{2}}$ on the annulus $A(0;r_n,R)$
we have $\left|e^{2n \varphi(z)} \right| \leq e^{-2 c \sqrt{n}}$
for $z \in \Delta_1^{\pm} \cap A(0;r_n,R)$.
Then \eqref{Rjump2} follows since the entries of $P(z)$ and $P^{-1}(z)$ 
have at most a power law singularity as $z \to 0$.

\medskip

Part (b) follows from part (a) by standard arguments from
Riemann-Hilbert theory. Compare for example with
the arguments in Appendix A from \cite{BlKu2} that deal explicitly
with a RH problem containing a varying shrinking contour.
\end{proof}

\subsection{Rewriting of the correlation kernel}

Before we take the scaling limit we follow the transformations
$Y \mapsto X \mapsto T \mapsto S \mapsto R$ to see
their effect on the correlation kernel for $x,y$ near $0$.

\begin{lemma} \label{lemKVninPhi}
For $x,y \in (0,r_n)$, we have
\begin{multline} \label{KVninPhi}
	K_{V,n}^{\alpha,\frac{1}{2}}(x,y) = e^{\frac{2}{3}n (V(x) - V(y))} \\
	\frac{1}{2\pi i(x-y)} \begin{pmatrix} -1 & 1 & 0 \end{pmatrix}
	\Phi_{\alpha,+}^{-1}(n^3 f(y))
	E_{in}^{-1}(y) R^{-1}(y) R(x) E_{in}(x)
	\Phi_{\alpha,+}(n^3f(x)) \begin{pmatrix} 1 \\ 1 \\ 0 \end{pmatrix}.
	\end{multline}
\end{lemma}
\begin{proof}
We start with the formula \eqref{CorKerY} for the correlation kernel
in terms of $Y$. The first two transformations
$Y \mapsto X \mapsto T$ from Definitions \ref{defX} and \ref{defT} lead to
\begin{align*} K_{V,n}^{\alpha,\frac{1}{2}}(x,y) & =
	\frac{\sqrt{2} y^{\beta} e^{-nV(y)}}{2\pi i(x-y)}
		\begin{pmatrix} 0 & 1 & 0 \end{pmatrix}
		X_+^{-1}(y) X_+(x) \begin{pmatrix} 1 \\ 0 \\ 0 \end{pmatrix}\\
	& = \frac{y^{\beta} e^{-n (V(y)-g_{1,+}(y)+g_2(y) + \ell)}
	e^{n g_{1,+}(x)}}{2\pi i (x-y)}
			\begin{pmatrix} 0 & 1 & 0 \end{pmatrix}
		T_+^{-1}(y) T_+(x) \begin{pmatrix} 1 \\ 0 \\ 0 \end{pmatrix}.
		\end{align*}
The next transformation $T \mapsto S$ is the opening of lenses.
For $x,y \in (0,q)$ we use \eqref{defS1} to obtain
\begin{multline} \label{KVninS} 
	K_{V,n}^{\alpha,\frac{1}{2}}(x,y)  =
	e^{-n(V(y)-g_{1,+}(y)+g_2(y) - \varphi_+(y) + \ell)}
	e^{n (g_{1,+}(x)+\varphi_+(x))} \\
	\frac{1}{2\pi i (x-y)}
			\begin{pmatrix} - e^{n \varphi_+(y)} & y^{\beta} e^{-n \varphi_+(y)} & 0 \end{pmatrix}
		S_+^{-1}(y) S_+(x) \begin{pmatrix} e^{-n \varphi_+(x)} \\ 
		x^{-\beta} e^{n \varphi_+(x)} \\ 0 \end{pmatrix}.
		\end{multline}
The scalar prefactor simplifies because of \eqref{phi} and we find
\begin{multline} \label{KVninS2} 
	K_{V,n}^{\alpha,\frac{1}{2}}(x,y)  =
	e^{\frac{1}{2} n (V(x) + g_2(x) - V(y) - g_2(y))} \\
	\frac{1}{2\pi i (x-y)}
			\begin{pmatrix} - e^{n \varphi_+(y)} & y^{\beta} e^{-n \varphi_+(y)} & 0 \end{pmatrix}
		S_+^{-1}(y) S_+(x) \begin{pmatrix} e^{-n \varphi_+(x)} \\ 
		x^{-\beta} e^{n \varphi_+(x)} \\ 0 \end{pmatrix}.
		\end{multline}
Now we take $x,y \in (0,r_n)$. We write $S = RP$ as in \eqref{Rdef} 
and $P$ is given by \eqref{Pansatz}. Then
\begin{multline} \label{Sphi1}
	S_+(x) \begin{pmatrix} e^{-n \varphi_+(x)} \\
		x^{-\beta} e^{n \varphi_+(x)} \\
		0 \end{pmatrix} 
		= R(x) E_{in}(x) \Phi_{\alpha,+}(n^3 f(x))
		\begin{pmatrix} 
	e^{-n (\varphi_+(x) - \frac{4}{3} \varphi_{1,+}(x) 
		- \frac{2}{3}\varphi_{2,+}(x))} \\ 
	e^{n (\varphi_+(x) - \frac{2}{3} \varphi_{1,+}(x) + 
		\frac{2}{3} \varphi_{2,+}(x))} \\ 0 
		\end{pmatrix} \\
	 = (-1)^n e^{\frac{1}{3} n (\varphi_{1,+}(x) + 2 \varphi_{2,+}(x))}
	 R(x) E_{in}(x) \Phi_{\alpha,+}(n^3 f(x))
		\begin{pmatrix} 1 \\ 1 \\ 0 \end{pmatrix}
		\end{multline}
where we used \eqref{Pansatz} and the formula \eqref{defD0} for $D_0$
for the first identity, and the fact
that $\varphi_+(x) = \varphi_{1,+}(x) - \pi i$ by \eqref{varphi1}
for the second identity in \eqref{Sphi1}. 
In a similar way,
\begin{multline} \label{Sphi2}
	\begin{pmatrix} - e^{n \varphi_+(y)} & y^{\beta} e^{-n \varphi_+(y)} & 0 \end{pmatrix}
	S_+^{-1}(y) \\
	= (-1)^n e^{-\frac{1}{3} n (\varphi_{1,+}(y) + 2 \varphi_{2,+}(y))}
	 \begin{pmatrix} -1 & 1 & 0 \end{pmatrix}
	 \Phi^{-1}_{\alpha,+}(n^3 f(y)) E_{in}^{-1}(y) R^{-1}(y).
		\end{multline}

We insert \eqref{Sphi1} and \eqref{Sphi2} into \eqref{KVninS}. Then
\eqref{KVninPhi} follows, since we also have 
\[ \tfrac{1}{2} V(x) + g_2(x) + \tfrac{1}{3} \varphi_{1,+}(x) + 
\tfrac{2}{3} \varphi_{2,+}(x) =
	\tfrac{2}{3} V(x) + \tfrac{1}{6} \ell \]
because of \eqref{phi0}, \eqref{varphi1}, and \eqref{varphi2}.	
\end{proof}

\subsection{Proof of Theorem \ref{mainThm}}

In view of \eqref{lemKVninPhi} and \eqref{Kformula} the
strategy for the proof of Theorem \ref{mainThm} is now clear. 

\begin{proof}[Proof of Theorem \ref{mainThm}]
We put $ x_n = \frac{x}{(c_V n)^3}$ and $y_n = \frac{y}{(c_V n)^3}$
with $x,y > 0$ fixed and $c_V$ is the constant defined
in Theorem \ref{mainThm}. 
It follows from Proposition \ref{f1f2ana} that $c_V^3 = f'(0)$
and therefore $ \Phi_{\alpha,+}(n^3 f(x_n)) \to \Phi_{\alpha,+}(x)$
and $\Phi_{\alpha,+}(n^3 f(y_n)) \to \Phi_{\alpha,+}(y)$
as $n \to \infty$, uniformly for $x,y$ in compact subsets of $(0,\infty)$.

We need to show that 
\begin{equation} \label{ERtoshow} 
	E_{in}^{-1}(y_n) R^{-1}(y_n) R(x_n) E_{in}(x_n) \to \mathbb I 
	\end{equation}
as $n \to \infty$, and then \eqref{Kalth} follows indeed
from  \eqref{KVninPhi} and \eqref{Kformula}. The proof
of \eqref{ERtoshow} (in a more precise form) is in Lemma \ref{lem:RestEstimate}
below. To establish this lemma, we first need two other
lemmas that we state and prove separately.
\end{proof}

We assume $x_n = \frac{x}{(c_V n)^3}$ and $y_n = \frac{y}{(c_V n)^3}$.

\begin{lemma} \label{lem:Restimate}
We have uniformly for $x,y$ in compact subsets of $[0,\infty)$,
\begin{align} \label{RinvRlimit}
R^{-1}(y_n) R(x_n) = 
	\mathbb{I} + \mathcal{O}\left((x-y) n^{-\frac{5}{2}}\right)
\end{align}
as $n \to \infty$.
 \end{lemma}
\begin{proof} 
Cauchy's integral formula yields that 
\begin{align*}
R(x_n) - R(y_n) & = \frac{1}{2\pi i} \oint_{|z|= r_n} 
	\left(	\frac{R(z)-\mathbb{I}}{z-x_n} - \frac{R(z)-\mathbb{I}}{z-y_n} \right) dz \\
	& = \frac{x_n-y_n}{2\pi i} \oint_{|z|=r_n}
		\frac{R(z) - \mathbb{I}}{(z-x_n)(z-y_n)} dz
\end{align*}
By \eqref{RisI} we have $R(z) - \mathbb{I} = \mathcal{O}(n^{-1})$
and it follows by easy estimation that
$ R(x_n) - R(y_n) = \mathcal{O}\left((x-y) n^{-\frac{5}{2}} \right)$
as $n \to\infty$. As a result, since $R^{-1}(y_n)$ remains bounded,
\begin{align*}
R^{-1}(y_n) R(x_n) & = \mathbb{I} + R^{-1}(y_n) (R(x_n)-R(y_n)) 
 =  \mathbb{I} + \mathcal{O}\left((x-y) n^{-\frac{5}{2}}\right).
\end{align*}
and the constant implied in the $\mathcal{O}$ term is uniform
for $x$ and $y$ in compact subsets of $[0,\infty)$. 
\end{proof}

\begin{lemma} \label{lem:Eestimate}
We have uniformly for $x,y$ in compact subsets of $[0,\infty)$,
\begin{align} \label{EinvElimit}
E_{in}^{-1}(y_n) E_{in}(x_n) = 
	\mathbb{I} + \mathcal{O}\left((x-y) n^{-\frac{1}{2}}\right)
\end{align}
as $n \to \infty$.
 \end{lemma}
\begin{proof}
We note that by  Definition \ref{defEinEout}
and \eqref{defEin1}, \eqref{E2in}, \eqref{E3in}, \eqref{E4in}
\begin{equation} \label{Einfactors} 
	E_{in}(z) = \frac{\sqrt{3}}{2\pi} n^{2\beta} 
	\left(\frac{f(z)}{z} \right)^{\frac{2\beta}{3}}
	\left(\mathbb I - B_n^{(3)}(z) \right) 
	\left(\mathbb I - B_n^{(2)}(z) \right)	
	\left( \mathbb I - B_n^{(1)}(z) \right)	
	E_n(z)
\end{equation}
with
\begin{equation} \label{Bnk} 
	B_n^{(k)}(z) = \frac{A_n^{(k)}(z) - A_n^{(k)}(0)}{n^{3k}z},
	\qquad \text{ for } k =1,2,3. 
	\end{equation}
We are going to show that for every $C > 0$ and
for $|z|, |s| \leq Cn^{-\frac{5}{2}}$ we have as $n \to \infty$,  
\begin{align} \label{Anktoshow}
	E_n^{-1}(z) A_n^{(k)}(s) E_n(z) & = O\left(n^{\frac{5(k-1)}{2}}\right),\\
	\label{Bnktoshow}
	E_n^{-1}(z) B_n^{(k)}(s) E_n(z) & = O\left(n^{-\frac{k}{2}}\right),
	\end{align}
for $k=1,2,3$.

We start by noting that for $z,s = \mathcal{O}(n^{-\frac{5}{2}})$ 
one has by Lemma \ref{En1-En2Lem} that $E_n^{-1}(z) E_n(s) = \mathcal{O}(1)$,
and then \eqref{Anktoshow} with $k=1$ follows immediately from
this and the definition \eqref{defAn1} of $A_n^{(1)}$.

Note that by Cauchy's formula
\begin{equation} \label{An1formula} 
	E_n^{-1}(z) 
	\frac{A_n^{(1)}(s) - A_n^{(1)}(0)}{s} E_n(z)
		= \frac{1}{2\pi i} \oint_{|t|= 2Cn^{-\frac{5}{2}}}
			\frac{E_n^{-1}(z) A_n^{(1)}(t) E_n(z)}{t(t-s)} dt 
			\end{equation}
which we can estimate because of \eqref{Anktoshow} with $k=1$ to give
\begin{equation} \label{An1estimate} 
	E_n^{-1}(z) 
	\frac{A_n^{(1)}(s) - A_n^{(1)}(0)}{s} E_n(z)
	=  \frac{\mathcal{O}\left(1\right)}{n^{-\frac{5}{2}}}
		= \mathcal{O}\left(n^{\frac{5}{2}} \right)  
	\end{equation}
whenever $|z|, |s| \leq  Cn^{-\frac{5}{2}}$. 
Then \eqref{Bnktoshow} with $k=1$ follows because of \eqref{Bnk}.

Then by \eqref{An2}, \eqref{Bnk}, and \eqref{Anktoshow} with $k=1$
\begin{align*} E_n^{-1}(z) A_n^{(2)}(s) E_n(z)
	& = n^3 \left( E_n^{-1}(z) A_n^{(1)}(0) E_n(z) \right)
	\left(E_n^{-1}(z) B_n^{(1)}(s) E_n(z) \right) \\
	&= n^3 \cdot \mathcal{O}(1) \cdot 
		\mathcal{O}\left(n^{-\frac{1}{2}}\right)
		= \mathcal{O}\left(n^{\frac{5}{2}} \right) \end{align*}
for  $|z|, |s| \leq  Cn^{-\frac{5}{2}}$, which 
proves \eqref{Anktoshow} with $k=2$.

From formula \eqref{An3} for $A_n^{(3)}$ we find in a similar way
\begin{align*} E_n^{-1}(z) A_n^{(3)}(s) E_n(z)
	& = - n^6  \left(E_n^{-1}(z) B_n^{(1)}(s) E_n(z) \right)
	\left( E_n^{-1}(z) A_n^{(1)}(s) E_n(z) \right)
	\left(E_n^{-1}(z) B_n^{(1)}(s) E_n(z) \right) \\
	&= n^6 \cdot \mathcal{O}\left(n^{-\frac{1}{2}}\right) 
	\cdot \mathcal{O}(1) \cdot 
		\mathcal{O}\left(n^{-\frac{1}{2}}\right)
		= \mathcal{O}\left(n^{5} \right) \end{align*}
for  $|z|, |s| \leq  Cn^{-\frac{5}{2}}$. This proves \eqref{Anktoshow}
with $k=3$.

Finally,  \eqref{Bnktoshow} for $k=2,3$ follows from
\eqref{Anktoshow} with $k=2,3$ in the same way as we obtained
\eqref{Bnktoshow} for $k=1$ with the identity \eqref{An1formula} 
and the estimate \eqref{An1estimate}.

Having \eqref{Bnktoshow} we go back to \eqref{Einfactors} which
we rewrite and estimate for $z = \mathcal{O}(n^{-3})$ as
\begin{align*}
	E_{in}(z) &  = \frac{\sqrt{3}}{2\pi} n^{2\beta} 
	\left(\frac{f(z)}{z} \right)^{\frac{2\beta}{3}}
	E_n(z) \left(\mathbb I - E_n^{-1}(z) B_n^{(3)}(z) E_n(z) \right) \\
	& \qquad \times
	\left(\mathbb I - E_n^{-1}(z) B_n^{(2)}(z) E_n(z) \right)	
	\left( \mathbb I - E_n^{-1}(z) B_n^{(1)}(z) E_n(z) \right)	\\
	& = \frac{\sqrt{3}}{2\pi} n^{2\beta} 
	\left(\frac{f(z)}{z} \right)^{\frac{2\beta}{3}}
	E_n(z) \left( \mathbb I + \mathcal{O}(n^{-\frac{3}{2}}) \right)
	\left( \mathbb I + \mathcal{O}(n^{-1}) \right)
	\left( \mathbb I + \mathcal{O}(n^{-\frac{1}{2}}) \right) \\
	& = \frac{\sqrt{3}}{2\pi} n^{2\beta} 
	\left(\frac{f(z)}{z} \right)^{\frac{2\beta}{3}}
	E_n(z)	\left( \mathbb I + \mathcal{O}(n^{-\frac{1}{2}}) \right).
\end{align*}
Since all factors are analytic and invertible we in fact have
\begin{equation} \label{EinCn} 
	E_{in}(z) = \frac{\sqrt{3}}{2\pi} n^{2\beta} 
	\left(\frac{f(z)}{z} \right)^{\frac{2\beta}{3}}
	E_n(z) \left(\mathbb{I} + C_n(z) \right), \end{equation}
where $C_n$ is analytic and 
\begin{equation} \label{Cnestimate} 
	C_n(z) = \mathcal{O}(n^{-\frac{1}{2}}), \qquad \text{for } 
	z = \mathcal{O}(n^{-3}).
	\end{equation}

We use \eqref{EinCn} for $z=x_n$  and $z= y_n$.
 Since $f$ is analytic with $f(0) = 0$, $f'(0) > 0$, 
 it is easy to check that
\begin{equation} \label{fcancels} 
	\left(\frac{f(x_n)}{x_n} \right)^{\frac{2\beta}{3}}
	\left( \frac{f(y_n)}{y_n} \right)^{- \frac{2\beta}{3}}
		= 1 + \mathcal O\left((x-y) n^{-3}\right) \end{equation}
as $n \to \infty$.
In view of \eqref{Enestimate} we have $E_n^{-1}(y_n) E_{n}(x_n) = 
	\mathbb{I} + \mathcal{O} \left((x-y) n^{-\frac{1}{2}}\right)
	$
and we obtain from this,  \eqref{EinCn}, \eqref{Cnestimate}, 
and \eqref{fcancels} 
\begin{align*} 
	E_{in}^{-1}(y_n) E_{in}(x_n)
	& = 	\left( \mathbb{I} + C_n(y_n) \right)^{-1}
	\left( \mathbb{I} + \mathcal{O}	\left((x-y) n^{-\frac{1}{2}}\right)
	\right) \left( \mathbb{I} + C_n(x_n) \right) \\
	& = \left( \mathbb{I} + C_n(y_n) \right)^{-1}
	 \left( \mathbb{I} + C_n(x_n) \right)
	 	+  \mathcal{O}	\left((x-y) n^{-\frac{1}{2}}\right) \\
	 & = \mathbb{I} +
	 	\left( \mathbb{I} + C_n(y_n) \right)^{-1}
	 \left(C_n(x_n) - C_n(y_n) \right)
	 	 	+  \mathcal{O}	\left((x-y) n^{-\frac{1}{2}}\right)
	 \end{align*}
From the fact that $C_n(z)$ is analytic with 
$C_n(z) = O\left(n^{-\frac{1}{2}} \right)$ for $z = \mathcal{O}(n^{-3})$ it follows (with an argument based on Cauchy's formula), that
\[  C_n(x_n) - C_n(y_n) = \mathcal{O}\left((x-y) n^{-\frac{1}{2}} \right) \]
and \eqref{EinvElimit} follows.
\end{proof}

\begin{lemma} \label{lem:RestEstimate}
As $n \to \infty$
\begin{equation} \label{RestEstimate}
	E_{in}^{-1}(y_n) R^{-1}(y_n) R(x_n) E_{in}(x_n)
= \mathbb I + \mathcal{O}\left((x-y) n^{-\frac{1}{2}} \right)
\end{equation}
uniformly for $x,y$ in compact subsets of $[0,\infty)$.
\end{lemma}
\begin{proof}
By \eqref{RinvRlimit} and \eqref{EinvElimit} we have
\begin{align} \nonumber
E_{in}^{-1}(y_n)  R^{-1}(y_n) R(x_n) E_{in}(x_n) 
	& =  E_{in}^{-1}(y_n) \left(\mathbb{I}+\mathcal{O}((x-y) n^{-5/2})\right) E_{in}(x_n) \\
	& = \mathbb{I} + \mathcal{O}\left((x-y)n^{-\frac{1}{2}} \right) 
	+ \label{ERestimate1}
	E_{in}^{-1}(y_n) 
	\mathcal{O}\left((x-y) n^{- \frac{5}{2}} \right) E_{in}(x_n).
	\end{align}
From \eqref{EinCn} and \eqref{fcancels} we obtain 
\begin{align} \nonumber
	E_{in}^{-1}(y_n) 
	\mathcal{O}\left((x-y) n^{-\frac{5}{2}})\right) E_{in}(x_n)
 	&	= \left(\mathbb{I} + C_n(y_n) \right)^{-1} 
		E_n^{-1}(y_n) 
		\mathcal{O}\left((x-y) n^{- \frac{5}{2}}\right)
		E_n(x_n) \left( \mathbb{I} + C_n(x_n) \right) \\
		\label{ERestimate2}
	& = E_n^{-1}(y_n) 
	\mathcal{O}\left((x-y) n^{- \frac{5}{2}}\right)
	E_n(x_n) \end{align}
where in the second step we used \eqref{Cnestimate}. Finally recall that $E_n(x_n) = \mathcal{O}(n)$ 
and $E_n^{-1}(y_n) = \mathcal{O}(n)$ by Lemma \ref{En1-En2Lem}.
Combining this with \eqref{ERestimate1}, \eqref{ERestimate2}, we find
\eqref{RestEstimate}, and the lemma follows. 
\end{proof}

\section*{Acknowledgements}

ABJK is supported by long term structural funding - Methusalem grant
of the Flemish Government, and by FWO Flanders projects G.0864.16 and EOS G0G9118N. 
LDM is supported by a PhD fellowship of the Research Foundation Flanders (FWO).

\phantomsection

\end{document}

%% file: FiguurS.tex
%
%
\resizebox{10.5cm}{6cm}{%
\begin{tikzpicture}[>=latex]
	\draw[-] (-7,0)--(7,0);
	\draw[-] (-3,-4)--(-3,4);
	\draw[fill] (-3,0) circle (0.1cm);
	\draw[fill] (3,0) circle (0.1cm);
	
	\node[above] at (-3.25,0) {\large 0};	
	\node[above] at (3.1,0.05) {\large $q$};	
	\node[above] at (-3.5,1.8) {\large $\Delta_1^+$};	
	\node[below] at (-3.5,-1.8) {\large $\Delta_1^-$};	
	\node[above] at (-5,0.05) {\large $\Delta_2$};	
	\node[above] at (0,0.05) {\large $\Delta_1$};
	
	\draw[->, ultra thick] (-3,0) -- (-3,1.2);
	\draw[-, ultra thick] (-3,0) -- (-3,2);
	
	\draw[->, ultra thick] (-3,0) -- (-3,-1.2);
	\draw[-, ultra thick] (-3,0) -- (-3,-2);
	
	\draw[->, ultra thick] (3,0) -- (5,0);
	\draw[->, ultra thick] (-3,0) -- (0,0);
	\draw[-, ultra thick] (-3,0) -- (7,0);
	
	\draw[->, ultra thick] (-7,0) -- (-4.5,0);
	\draw[-, ultra thick] (-7,0) -- (-3,0);
	
	\draw[-, ultra thick] (-3,2) to [out=0, in=135] (3,0);
	\draw[-, ultra thick] (-3,-2) to [out=0, in=225] (3,0);

	\draw[->, ultra thick] (-0.9,1.88) to (-0.7,1.85);
	\draw[->, ultra thick] (-0.9,-1.88) to (-0.7,-1.85);
	
\end{tikzpicture}
}
%

%% file: FiguurP.tex
%
%
\resizebox{6cm}{6cm}{%
\begin{tikzpicture}[>=latex]
	\draw[-] (-5,0)--(5,0);
	\draw[-] (0,-5)--(0,5);
	\draw[fill] circle (0.15cm);
	\draw[fill] (1,0) circle (0.15cm);
	\draw[fill] (3,0) circle (0.15cm);
	\draw[ultra thick] (0,0) circle (1 cm);
	\draw[ultra thick] (0,0) circle (3cm);
	\node[above] at (-0.5,0) {\huge 0};
	\node[above] at (1.35,0) {\huge $r_n$};
	\node[above] at (3.45,0) {\huge $R$};
	\node[right] at (0,1.5) {\huge $\Delta_1^{+}$};
	\node[below] at (4,0) {\huge $\Delta_1$};
	\node[left] at (0,-1.5) {\huge $\Delta_1^{-}$};
	\node[above] at (-4,0) {\huge $\Delta_2$};	
	
	\draw[->, ultra thick] (1.5,0.0) -- (2.0,0.0);
	\draw[-, ultra thick] (0,0) -- (3,0);
	
	\draw[->, ultra thick] (-2.0,0) -- (-1.5,0);
	\draw[-, ultra thick] (-3,0) -- (0,0);
		
	\draw[->, ultra thick] (0,0) -- (0,0.8);
	\draw[-, ultra thick] (0,0) -- (0,1);
	
	\draw[->, ultra thick] (0,0) -- (0,-0.8);
	\draw[-, ultra thick] (0,0) -- (0,-1);
	
	\draw[->, ultra thick] (2.16,2.1) to (2.06,2.2);
	\draw[->, ultra thick] (0.7,0.7) to (0.6,0.8);
	
\end{tikzpicture}
}
%

%% file: FiguurR.tex
%
%
\resizebox{10.5cm}{6cm}{%
\begin{tikzpicture}[>=latex]
	\draw[-] (-7,0)--(7,0);
	\draw[-] (-4,-4)--(-4,4);
	\draw[fill] (-4,0) circle (0.1cm);
	\draw[fill] (3,0) circle (0.1cm);
	\draw[ultra thick] (-4,0) circle (1.5cm);
	\draw[ultra thick] (-4,0) circle (0.5cm);
	\draw[ultra thick] (3,0) circle (1.5cm);
	\node[above] at (-4.2,-0.075) {\large 0};	
	\node[above] at (2.8,0) {\large $q$};	
	
	\draw[->, ultra thick] (-4,0.5) -- (-4,1.2);
	\draw[-, ultra thick] (-4,0.5) -- (-4,3);
	
	\draw[->, ultra thick] (-4,-0.5) -- (-4,-1.2);
	\draw[-, ultra thick] (-4,-0.5) -- (-4,-3);
	
	\draw[->, ultra thick] (4.5,0) -- (5.5,0);
	\draw[-, ultra thick] (4.5,0) -- (7,0);
	
	\draw[-, ultra thick] (-4,3) to [out=0, in=135] (1.939,1.061);
	\draw[-, ultra thick] (-4,-3) to [out=0, in=225] (1.939,-1.061);

	\draw[->, ultra thick] (-2.86,0.96) to (-2.96,1.08);
	\draw[->, ultra thick] (-3.65,0.35) to (-3.75,0.47);	
	\draw[->, ultra thick] (4.14,0.96) to (4.04,1.08);
	
	\draw[->, ultra thick] (-1.5,2.85) to (-1.3,2.82);
	\draw[->, ultra thick] (-1.5,-2.85) to (-1.3,-2.82);
\end{tikzpicture}
}
